\input amstex
\documentstyle{amsppt} 

\magnification=1200 \font\ftitle=cmbx10 scaled\magstep2 \line{}
\vskip4.5truecm \nologo \NoBlackBoxes

\topmatter \centerline{\ftitle A Theory of Dimension}\bigskip
\centerline{{\smc  Roberto Longo$^{1*}$} and {\smc John E.
Roberts}\footnote"$^*$"{Supported in part by MURST  and
CNR-GNAFA\hfill\break
 E-mail  longo\@tovvx1.ccd.utovrm.it, roberts\@mat.utovrm.it }}
\bigskip \centerline{ Dipartimento di Matematica, Universit\`a di
Roma ``Tor Vergata''} \centerline {via della Ricerca Scientifica,
I--00133 Roma, Italy} \medskip\centerline{$^1$Centro Linceo
Interdisciplinare, Accademia Nazionale dei Lincei} \centerline{via
della Lungara 10, I--00165 Roma, Italy} \vskip 2truecm

\centerline{{\it Dedicated to R.V. Kadison on the occasion of his
seventieth birthday}} \vskip 4truecm

{\it Abstract.} In which a theory of dimension related to the Jones
index  and based on the notion of conjugation is developed. An
elementary proof  of the additivity and multiplicativity of the
dimension is given and there is  an associated trace. Applications
are given to a class of endomorphisms of factors and to the theory
of subfactors. An important role is played by a notion of 
amenability inspired by the work of Popa.

\vfill\eject

scaled\magstep1 \newcount\REFcount \REFcount=1
\def\numref{\number\REFcount} \def\addref{\global\advance\REFcount
by 1} \def\wdef#1#2{\expandafter\xdef\csname#1\endcsname{#2}}
\def\wdch#1#2#3{\ifundef{#1#2}\wdef{#1#2}{#3}
    \else\write16{!!doubly defined#1,#2}\fi}
\def\wval#1{\csname#1\endcsname}
\def\ifundef#1{\expandafter\ifx\csname#1\endcsname\relax}
\def\ref(#1){\wdef{q#1}{yes}\ifundef{r#1}$\diamondsuit$#1
  \write16{!!ref #1 was never defined!!}\else\wval{r#1}\fi}
\def\inputreferences{
    \def\REF(##1)##2\endREF{\wdch{r}{##1}{\numref}\addref}
    \REFERENCES} \def\references{
    \def\REF(##1)##2\endREF{
        \ifundef{q##1}\write16{!!ref. [##1] was never quoted!!}\fi
        \item{[\ref(##1)]}##2}
   \par\REFERENCES}  
	\def\REFERENCES{
 	\REF(BaaSka){\smc Baaj, S., Skandalis, G., } {\it Unitaries
multiplicatifs et dualit\'e pour les produits croises de
$C^*$--algebres}, Ann. Sci. Ec.  Norm. Sup. {\bf 26} (1993),
425--488.\endREF
 \REF(Cho){\smc Choda, M.}, {\it Extension algebras of $II_1$
factors via $^*$-endomorphisms}, preprint 1994.\endREF
 \REF(Con){\smc Connes, A., } {\it Une classification des facteurs
de type III}, Ann. Sci. Ec. Norm. Sup. {\bf 6}, (1973), 133--252.
 \endREF
 \REF(DHR1){\smc Doplicher, S., Haag, R., Roberts, J.E.,} {\it Local 
observables and particle statistics I}, Commun. Math. Phys.  {\bf
23} (1971), 199--230.
 \endREF
 \REF(DHR2){\smc Doplicher, S., Haag, R., Roberts, J.E., } {\it
Local  observables and particle statistics II}, Commun. Math. Phys. 
{\bf 35} (1974), 49--85.
 \endREF
  \REF(DR1){\smc Doplicher, S., Roberts, J.E.}, {\it A new duality
theory for compact groups}, Invent. Math. {\bf 98} (1989),
157--218.\endREF \REF(DR2){\smc Doplicher, S., Roberts, J.E.}, {\it
Duals of compact Lie groups realized in the Cuntz algebras and their
actions on C$^*$--algebras},  J. Funct. Anal. {\bf 74} (1987),
96--120.\endREF
 \REF(FrKer){\smc Fr\" ohlich J., Kerler T.,}{\it Quantum groups,
quantum categories and quantum field theory}, Lect. Notes Math. {\bf
1542}, Springer Verlag, New York-Heidelberg-Berlin 1993.\endREF 
\REF(Kos1){\smc Kosaki, H., } {\it Extension of Jones' theory on
index to arbitrary subfactors}, J. Funct. Anal. {\bf 66} (1986),
123--140.\endREF
 \REF(KosLong){\smc Kosaki, H., Longo, R.,} {\it A remark on the
minimal index of subfactors}, J. Funct. Anal. {\bf 107} (1992),
458--470.\endREF
 \REF(Jones){\smc Jones, V.R.F., } {\it Index for subfactors},
Invent. Math. {\bf 72} (1983), 1--25. \endREF \REF(GLR){\smc Ghez,
P., Lima, R., Roberts, J.E., } {\it $W^*$-Categories}, Pac. J. Math.
{\bf 120} (1985), 79--109. \endREF \REF(KW){\smc Kawakami, S.,
Watatani, Y.,} {\it The multiplicativity  of the minimal index of
simple C$^*$-algebras, } Proc. Amer. Mat. Soc. (to appear).\endREF
 \REF(Long1){\smc Longo, R., } {\it Index of subfactors and
statistics  of quantum fields. I}, Commun. Math. Phys. {\bf 126}
(1989), 217--247.
 \endREF
 \REF(Long2){\smc Longo, R., } {\it Index of subfactors and
statistics of quantum fields. $II$ Correspondences, braid group
statistics and Jones polynomial}, Commun. Math. Phys. {\bf 130}
(1990), 285--309. \endREF \REF(Long3){\smc Longo, R.,}$\,$ {\it
Minimal index and braided subfactors},$\,$  J.$\,$ Funct.$\,$
Anal.$\,$ {\bf 109}$\,$ (1992),$\,$ 98--112. \endREF 
\REF(Long4){\smc Longo, R., } {\it Solution of the factorial
Stone--Weierstrass conjecture}, Invent. Math. {\bf 76} (1984),
145--155. \endREF
 \REF(Long6){\smc Longo, R.,} {\it A duality for Hopf algebras and
for subfactors. I}, Commun. Math. Phys. {\bf 159} (1994),
133--150.\endREF
  \REF(McL){\smc Mac Lane, S., } {\it Categories for the working
mathematician},
  Springer Verlag, Berlin-Heidelberg-New York 1971.\endREF
 \REF(Ocn){\smc Ocneanu, A., } {\it Quantized groups, string
algebras, and Galois theory}, in ``Operator Algebras and
Applications'' vol.$\!$ II, London Math. Soc. Lecture Note Series
{\bf 136}, Cambridge Univ. Press, 1988.\endREF
 \REF(PiPo){\smc Pimsner, M., Popa, M., } {\it Entropy and index for
subfactors}, Ann. Sci. Ec. Norm. Sup. {\bf 19} (1986),
57--106.\endREF 
 \REF(Popa1){\smc Popa, S., } {\it Classification of subfactors:
reduction to commuting squares}, Invent. Math. {\bf 101} (1990),
19--43.\endREF \REF(Popa2){\smc Popa, S., } {\it Classification of
amenable subfactor of type $II$}, Acta Math. {\bf 172} (1994),
352--445.\endREF \REF(Popa3){\smc Popa, S., } {\it Classification of
subfactors and of their endomorphisms}, CBMS lectures, preprint
1994.\endREF \REF(Rob1){\smc Roberts, J.E.,} {\it New light on the
mathematical structure of  algebraic field theory}, Proc. of Symp.
in Pure Math. {\bf 38} (1982), 523--550. \endREF
 \REF(SaRi){\smc Saavedra Rivano, N.,} {\it Cat\'egories
Tannakiennes}, Lect.$\!$ Notes
 in Math.$\!$ {\bf 265} Springer Verlag, New
York--Heidelberg--Berlin 1972.\endREF
 \REF(Take2){\smc Takesaki, M., } {\it Conditional expectation in von
Neumann algebras}, J. Funct. Anal. {\bf 9} (1970), 306--321 \endREF
 \REF(Yam){\smc Yamagami S.,} {\it Non-existence of minimal actions
for compact quantum groups}, preprint 1994 \endREF
 \REF(Wen1){\smc Wenzl, H., } {\it Hecke algebras of type $A_n$ and
subfactors}, Invent. Math. {\bf 92} (1988), 349--383.\endREF 
\REF(Wor1){\smc Woronowicz, S.L.}, {\it Compact matrix
pseudogroups}, Commun.  Math. Phys. {\bf 111} (1987), 613--665.
\endREF}

\inputreferences

\NoBlackBoxes \beginsection{1. Introduction}

The starting point for this paper is the operation of conjugation on 
finite dimensional objects in a  strict tensor $C^*$--category. The
terminology derives from group  theory where to every unitary
representation of a group there corresponds a  conjugate
representation defined on the conjugate Hilbert space. It includes
as  a special case the notion  of conjugate endomorphisms of a
$C^*$--algebra which was effectively first introduced in
[\ref(DHR2)], Thm. 3.3. If we look on a $C^*$--algebra as a
$C^*$--category with a single object then an endomorphism becomes an 
endofunctor and the conjugate endomorphism is just the left (and
indeed the right)  adjoint endofunctor. In this sense, the concept
is just a special case of the concept  of adjoint functor. The
notion of conjugate objects in a strict tensor $C^*$--category is 
also not really new either since it is implicit, for example, in the
notion of dual  object in [\ref(SaRi)].

When conjugates are understood in the sense of adjoint functors then 
the objects in question are necessarily finite dimensional in a
sense which will be made precise. In the case of von Neumann
algebras it corresponds to the finiteness of the index first
introduced by Jones [\ref(Jones)].  The relation between the index
and the square of the dimension  was first established in
[\ref(Long1)] whereas the concept  of dimension goes  back to
[\ref(DHR1)]  where it was referred to as the statistical dimension 
and  defined using the permutation symmetry on the strict tensor 
$C^*$--category.  This had the great merit that  the existence of
conjugates for finite dimensional objects could be established. 
However, although it was clear that the  dimensions of objects with
conjugates were independent of the choice of permutation symmetry 
it was not clear that the dimension could be defined in the absence
of any  permutation symmetry and hence made sense in any strict
tensor $C^*$--category.   Again there is a related notion of the
rank of an object with a dual defined for symmetric tensor
categories [\ref(SaRi)]  whose definition again involves the
permutation symmetry and which, in general,  depends on the choice
of that symmetry.

In Section 2, the properties of conjugation are developed for finite
dimensional  objects in a strict tensor $C^*$--category. In
particular we show how there is a related  notion of conjugation on
arrows and associated left and right scalar products on the 
category. We introduce the notions of left and right inverses for
objects and prove an  analogue of the Pimsner--Popa inequality
[\ref(PiPo)], already established in the particular context of
[\ref(DHR1)]. 

In Section 3, the dimension is introduced in the special case that
the tensor unit is  irreducible and shown to be additive on direct
sums and multiplicative on tensor products. It is  therefore a
Perron--Frobenius eigenvalue. It can also be characterized by a
property of  minimality related to the concept of minimal index. As
a consequence a strict tensor  $^*$--functor cannot increase
dimensions. We introduce the Jones projections  associated with
conjugation and conclude that the range of the dimension has the
same a priori restrictions as the square root of the Jones index. 
Associated with the dimension there is a canonical trace defined on
the full subcategory  of finite dimensional objects.

In Section 4, we briefly consider the status of conjugates and
dimension in a braided tensor $C^*$--category showing that the
braiding gives rise to a central element $\kappa$  and that the
braiding is automatically compatible with conjugates.

In Section 5, we study model endomorphisms canonically associated
with an object in a strict tensor $C^*$--category. The starting
point is the graded $C^*$--algebra $\Cal O_\rho$  associated with an
object $\rho$ in a strict tensor $C^*$--category introduced in
[\ref(DR1)].  When $d(\rho) < \infty$, the canonical trace defines a
trace on the grade zero subalgebra  which hence extends to a unique
gauge invariant state on $\Cal O_\rho$. This state  is shown to be a
KMS state with respect to gauge transformations. Taking a weak
closure in the  corresponding GNS--representation gives a von
Neumann algebra $\Cal M_\rho$ associated intrinsically with $\rho$.
$\Cal M_\rho$, like $\Cal O_\rho$, comes equipped with a  canonical
endomorphism $\hat \rho$. An important role is played here by a
notion of  amenability inspired by Popa's notion of a strongly
amenable inclusion of factors.  Roughly speaking amenability means
that if we pass from $\rho$ to $\hat\rho$  there are no new
intertwining operators. In fact, the theory is developed in greater
generality, replacing the standard left inverse which defines the
canonical trace by an arbitrary faithful left inverse $\varphi$ to
get a von Neumann algebra $\Cal M_\varphi$ equipped with a canonical
endomorphism. We illustrate our notion with examples. 

In Section 6,  we present applications to subfactors and the duality
theory for finite dimensional Hopf algebras. We introduce the notion
of  abstract $Q$--system, a categorical version of the $Q$--systems 
introduced in [\ref(Long6)]. We show that there is a 1--1
correspondence between isomorphism classes of amenable $Q$--systems
and of standard  AFD inclusions.

Section 7 has quite a different goal: the formalism of tensor
$C^*$--categories  is not sufficiently general to cover all
interesting applications. Instead, we need  to use
2--$C^*$--categories. This generalization presents no difficulties
and has  been deferred to Section 7 merely to avoid concepts that
are superficially technical. Here we present examples of
2--$C^*$--categories to illustrate the range of the resulting 
theory. We then generalize the basic notions of conjugation and
dimension, again providing  some key examples to illustrate the
utility of this approach. Finally, it is shown  in an appendix how a
2--$C^*$--category can completed so that it has subobjects and 
(finite) direct sums.

This paper by no means exhausts the basic themes of conjugation and
dimension so we  would like to conclude this introduction with some
comments on work in progress. In the  first place, although, in
defining dimension, we have restricted ourselves to the case that
the  endomorphisms of the tensor unit $\iota$ reduce to the complex
numbers, the  definitions used and  examples presented in
[\ref(DR1)] already establish the need for a  more general notion.
In the course of this work, it has on occasions proved irksome to
have to add the hypothesis that a von Neumann algebra is a factor
just to be able to talk about the dimension of an endomorphism.
However, although a notion of dimension can in fact be developed in
more generality on the basis of the ideas presented here, there may
still be room for improvement so details will be deferred to a
sequel. Again our notion of conjugation in terms of adjoint
functors  is restricted to the finite dimensional case despite the
fact that the notion of  conjugate representation makes sense in
infinite dimensions and the fact that  conjugates always exist for
normal endomorphisms of von Neumann algebras [\ref(Long2)]. We
propose to give a general definition covering the infinite
dimensional case elsewhere.  Although the results on inclusions we
have obtained in Section 6 using the canonical endomorphism are
relatively complete in themselves, it would still be of interest to
treat  the inclusion directly as a $\circ$--unit  in an appropriate
2--$C^*$--category. It would furthermore be of interest  to look at
duality at the level of abstract $Q$--systems.

\beginsection{2. Conjugation in Tensor C$^*$--Categories}

In this section we will work with strict tensor C$^*$--categories,
called strict {\it  monoidal} C$^*$--categories in [\ref(DR1)], \S
1, where the formal definition may  be found\footnote{It is
convenient to consider only {\it strict} tensor categories because
the generality achieved by dropping this restriction is minimal
[\ref(McL)].}. The change in name is justified since the set of
arrows between two objects  is automatically a linear space and
opportune because of the popularity and familiarity of  the term
tensor as opposed to monoidal. Briefly, such a category $\Cal T$
with objects  denoted by $\rho,\, \sigma,\, \tau,$ etc. has a space
of arrows between $\rho$ and $\sigma$,  $(\rho,\sigma)$, which is a
complex Banach space. Composition of arrows denoted  by $\circ$ is
bilinear and the adjoint $^*$ is an involutive contravariant functor
acting as  the identity on objects. The norm satisfies the
C$^*$--property: $\|R^*\circ R\|=\|R\|^2.$  This makes $\Cal T$ into
a C$^*$-category. To be a strict tensor C$^*$--category, we need  an
associative bilinear bifunctor $\otimes : \Cal T \times \Cal T \to
\Cal T$ with a unit $\iota$ and commuting with $^*$. To save space
and render expressions more legible, the  functor $\otimes$ applied
to a   pair of objects $\rho$ and $\sigma$ will be denoted simply
$\rho\sigma$ rather than $\rho\otimes\sigma$ as would be more
conventional. 

We begin by defining the notion of conjugation for finite dimensional
objects of an arbitrary strict tensor C$^*$--category  $\Cal T.$
This will eventually lead us to a notion of dimension related to 
the Jones index. We define a full subcategory $\Cal T_f$ of $\Cal T$
as follows:  an object $\rho$ of $\Cal T$ lies in $\Cal T_f$ if
there is an object $\bar \rho$ of $\Cal T$ and $R \in (\iota, \bar
\rho \rho)$ and $\bar R \in (\iota, \rho \bar \rho)$ such that
$$\bar R^*\otimes 1_\rho \circ 1_\rho\otimes R=1_\rho ; \quad
R^*\otimes 1_{\bar \rho} \,\circ \, 1_{\bar\rho}\,\otimes \bar
R=1_{\bar \rho}.$$ Here $\iota$ denotes the tensor unit in $\Cal T.$
The definition is symmetric in  $\rho$ and $\bar\rho$ so that
$\bar\rho$ is also an object of $\Cal T_f.$ The  above equations
will be referred to as the {\it conjugate equations} and  $\bar\rho$
will be said to be a {\it conjugate} for $\rho$. In essence, the
above condition means  that $\rho$ is finite dimensional. However,
one should bear in mind that $\rho$  may fail to be an object of
$\Cal T_f$ because its conjugate, although existing in a larger 
$C^*$--category containing $\Cal T$ as a full tensor
${}^*$--subcategory, fails  to exist in $\Cal T$. For this reason,
it is advisable to choose $\Cal T$ from  the outset so that it
contains all the necessary objects. In particular, it  is wise to
take $\Cal T$ to have subobjects and (finite) direct sums. 

Actually, as was pointed out in [\ref(DHR2)], the basic mathematical
concept underlying the notion of conjugation is that of adjoint
functor. If $\rho$ and $\bar \rho$  are objects in a tensor
$C^*$--category, then the functor of tensoring on the left by 
$1_\rho$ has the functor of tensoring on the left by $1_{\bar\rho}$
as a right adjoint precisely when we can find $R$ and $\bar R$ as
above. The condition for adjunction takes on such a simple form
because a tensor C$^*$--category has a tensor unit and a
$^*$--operation and because of the particular nature of the functors
involved. The symmetry between $\rho$ and $\bar \rho$ shows that
$1_{\bar \rho} \,\otimes$ is also a left adjoint for
$1_\rho\,\otimes$. The same conditions are also equivalent to 
$\otimes\,1_{\bar \rho}$ being a left (or right) adjoint to
$\otimes\, 1_\rho$. From the better known alternative equivalent
definition of adjoint functors it follows,  as can easily be
verified, that $\bar\rho$ is a conjugate for $\rho$ if and only if
there are isomorphisms $i_{\sigma , \tau}:(\rho\sigma,\tau) \to
(\sigma,\bar \rho \tau)$ natural in $\sigma$ and $\tau$.

To make the basic construction explicit, we state it in the form of
a lemma: \medskip \proclaim{Lemma 2.1} If $R$, $\bar R$ satisfy the
conjugate equations for $\rho$,  then for any $\rho_1$ and $\rho_2$,
the maps $$S \mapsto 1_{\bar\rho}\,\otimes\,S\,\circ R\,\otimes
1_{\rho_1}, \quad S\in  (\rho \rho_1,\rho_2),$$ $$S'\mapsto \bar
R^*\,\otimes\, 1_{\rho_2}\circ\, 1_{\rho}\,\otimes\, S', \quad S'\in
(\rho_{1},\bar \rho \rho_{2}).$$ are inverses of one another.
Similarly, $$T \mapsto T\,\otimes\, 1_{\bar
\rho}\,\circ\,1_{\rho_1}\,\otimes \bar R, \quad T\in
(\rho_1\rho,\rho_2),$$ $$T' \mapsto 1_{\rho_2}\,\otimes\,
R^*\,\circ\, T'\,\otimes 1_\rho, \quad  T' \in
(\rho_1,\rho_2\bar\rho)$$ are inverses of one another.\endproclaim

\demo{Proof} $$1_{\bar \rho}\otimes \bar R^*\otimes 1_{\rho_2}\circ
1_{\bar \rho \rho}\otimes S' \circ R\otimes 1_{\rho_1}=
1_{\bar\rho}\otimes\bar R^*\otimes 1_{\rho_2}\circ  R\otimes
1_{\bar\rho \rho_2}\circ S'= 1_{\bar \rho \rho_2} \circ S'=S',$$
$$\bar R^*\otimes 1_{\rho_2}\circ 1_{\rho \bar \rho}\otimes  S\circ
1_{\rho}\otimes R\otimes 1_{\rho_1}=S\circ \bar R^* \otimes 1_{\rho
\rho_1}\circ 1_{\rho}\otimes R\otimes 1_{\rho_1}=S\circ 1_{\rho
\rho_1}=S.$$

In the special case of group representations, this is in essence the 
Frobenius Reciprocity Theorem.
 
From the above lemma it follows in particular that any $R_1 \in
(\iota, \rho_1 \rho)$ is of the form $R_1=X\,\otimes\,
1_{\rho}\,\circ\, R$ for  a unique $X \in (\bar \rho, \rho_1)$ and
$R_1$ also defines a conjugate for  $\rho$ if and only if $X$ is
invertible. The corresponding $\bar R_1 \in
 (\iota, \rho \rho_1)$ is given by $\bar R_1=1_{\rho}\,\otimes\,
X^{*-1} \,\circ\,\bar R$. Similarly, any $R' \in (\iota,\bar \rho
\rho')$ is  of the form $R'=1_{\bar \rho}\,\otimes\, Y\,\circ\, R$
for a unique $Y \in (\rho,\rho')$ and $R'$ defines a conjugate for
$\rho'$ if  and only if $Y$ is invertible. The corresponding $\bar
R' \in (\iota, \rho'\bar \rho)$ is given by $\bar
R'=Y^{*-1}\,\otimes\, 1_{\bar \rho} \,\circ\, \bar R.$ Another
obvious consequence of the above isomorphisms is that the spaces
$(\rho, \rho), (\iota, \bar \rho \rho)$ and $(\bar  \rho, \bar
\rho)$ are isomorphic. In particular $\rho$ is irreducible  if and
only if $\bar \rho$ is irreducible.  \medskip We now show that
slightly weaker conditions suffice to define a conjugate.

\proclaim{Lemma 2.2} Two objects $\rho$ and $\bar \rho$ in a strict
tensor  C$^*$--category $\Cal T$  are conjugate if and  only if
there are $R \in (\iota, \bar \rho \rho)$ and $\bar R \in (\iota,
\rho \bar \rho)$ such that $\bar R^* \otimes 1_\rho \circ 1_\rho 
\otimes R$ is  invertible and the positive linear mapping $\omega $
on $(\bar \rho, \bar \rho) \to (\iota,\iota)$ defined by $$\omega
(T) = R^* \circ T \otimes 1_\rho \circ R, \ \ T \in  (\bar \rho,\bar
\rho),$$ is faithful. \endproclaim

\demo{Proof} If $R$ defines a conjugate for $\rho$ then $T \in (\bar 
\rho,\bar \rho)$ and $T \otimes 1_\rho \circ R = 0$ imply $T=0$ so
that  $\omega$ is faithful. Conversely, given $R$ and $\bar R$ as in
the  statement of the lemma, let $X=\bar R^* \otimes 1_\rho \circ
1_\rho \otimes  R \in (\rho, \rho),$ then replacing $\bar R$ by
$X^{-1 *} \otimes 1_\rho  \circ \bar R,$ we see that we may assume
that $X=1_\rho.$ But now $$ R = 1_{\bar \rho}   \otimes\bar R^*
\otimes 1_\rho \circ 1_{\bar \rho \rho} \otimes  R\circ R =
1_{\bar\rho} \otimes \bar R^* \otimes 1_\rho \circ R\otimes
1_{\bar\rho \rho}\circ R$$ and since $\omega$ is faithful, we
conclude that $1_{\bar \rho} \otimes  \bar R^*  \circ R \otimes
1_{\bar \rho}=1_{\bar \rho}$ and $\rho$  and $\bar \rho$ are
conjugates, completing the proof.\hfill$\square$ \medskip

In the category $\Cal T_f$ we can define a conjugation on arrows. We
pick for each object $\rho$ with a conjugate an $R_\rho \in  (\iota,
\bar \rho \rho)$ defining this conjugate and set for $T \in  (\rho,
\rho')$ $$T^\bullet:=1_{\bar \rho'} \otimes \bar R_{\rho}^* \circ
1_{\bar \rho'}  \otimes T^* \otimes  1_{\bar \rho} \circ R_{\rho'}
\otimes 1_{\bar \rho}\in (\bar \rho, \bar \rho').$$ This mapping is
just the composition of the mappings $$(\rho, \rho')\to (\iota,\rho'
\bar \rho)\overset* \to \rightarrow  (\rho' \bar \rho, \iota) \to
(\bar \rho, \bar \rho'),$$ where the first and last mappings are
defined as discussed above using the properties of adjoint functors.
A simple computation shows that $T^\bullet$ is the unique element of
$(\bar \rho, \bar \rho'),$ satisfying either $$T^\bullet \otimes
1_{\rho} \circ R_\rho=  1_{\bar \rho'} \otimes T^* \circ  R_{\rho'}
\quad \text {or} \quad  1_{\rho'} \otimes T^{\bullet *}  \circ \bar
R_{\rho'}= T \otimes 1_{\bar \rho} \circ \bar R_\rho.$$ At the same
time our choice $\rho \mapsto R_\rho$ serves to define an 
$(\iota,\iota)$--valued scalar product on  $\Cal T_f$ setting for
$S,T\in (\rho,\rho')$ $$(S,T)_{\ell} :=R_{\rho}{}^*\circ
1_{\bar\rho}\otimes(S^*\circ T)\circ R_{\rho}.$$ We obviously have
\item {a)} $(S,S)_{\ell} \geq 0,$\par \item {b)} $(S, X \circ
T)_{\ell} = (X^* \circ S, T)_{\ell},$\par \noindent whenever the
above expressions are well defined. \proclaim{Lemma 2.3} The mapping
$T \mapsto T^\bullet$ is antilinear and has the  following
properties \smallskip \item {$a$)} $1^\bullet_\rho = 1_{\bar
\rho},$\par \item {$b$)} $S^\bullet \circ T^\bullet = (S \circ
T)^\bullet$ \smallskip Furthermore if $(\cdot,\cdot)_\ell$ is
tracial, i.e. if $(S,T)_{\ell}=(T^*,S^*)_{\ell}$ then \item {$c$)}
$T^{\bullet *} = T^{*\bullet}.$ \endproclaim \demo{Proof} a) is
trivial and b) follows from the above formula and  $(S \circ T)^* =
T^* \circ S^*. $ To prove c) we note that if $S,T\in (\rho,\rho')$
then $$(S,T)_{\ell}=R_{\rho}{}^*\circ 1_{\bar\rho}\otimes S^*\circ
T^{* \bullet}\otimes 1_{\rho'} \circ R_{\rho'}$$ whereas
$$(T^*,S^*)_{\ell}=R_{\rho}{}^*\circ T^{\bullet *}\otimes
1_\rho\circ 1_{\bar\rho'}\otimes S^*\circ
R_{\rho'}=R_{\rho}{}^*\circ 1_{\bar\rho}\otimes S^*\circ T^{\bullet
*}\otimes 1_{\rho'} \circ R_{\rho'}.$$ Since, by Lemma 2.1,
$1_{\bar\rho}\otimes S\circ R_{\rho}$ is an arbitary element of 
$(\iota,\bar\rho\rho')$,  $c$) follows.\hfill$\square$\medskip
\proclaim{Theorem 2.4} The subcategory $\Cal T_f$ is closed under
conjugates and  tensor products. If $\Cal T$ has (finite) direct
sums and subobjects then $\Cal T_f$  is closed under direct sums
and  subobjects. \endproclaim \demo{Proof} $\Cal T_f$ is obviously
closed under conjugates and to see that it is closed under tensor
products let $R_i\in(\iota,\bar\rho_i\rho_i)$ and $\bar R_i\in
(\iota,\rho_i\bar\rho_i)$ define a conjugate for $\rho_i, i=1,2,$
then $$R :=1_{\bar\rho_2}\otimes R_1\otimes 1_{\rho_2}\circ
R_2;\quad \bar R := 1_{\rho_1}\otimes \bar R_2\otimes
1_{\bar\rho_1}\circ\bar R_1$$ solve the conjugate equations for
$\rho_1\rho_2$. To check closure under direct sums,  we suppose
$W_i\in (\rho_i,\rho)$ and $\bar W_i\in (\bar\rho_i,\bar\rho),
i=1,2$  to be isometries realizing $\rho$ and $\bar\rho$ as direct
sums and set  $$R=\sum_i \bar W_i\otimes W_i\circ R_i; \quad \bar
R=\sum_i W_i\otimes  \bar W_i\circ \bar R_i.$$ A simple computation
again shows that $R$ and $\bar R$ solve the conjugate equations for
$\rho$. It remains to check closure under subobjects. Suppose that
$\rho$ is an object of $\Cal T_f$ and $W\in (\rho',\rho)$ is an
isometry. Let $R,\bar R$ solve the conjugate equations for  $\rho$,
then, by Lemma 2.3,  $$1_{\bar\rho}\otimes (W\circ W^*)\circ
R=F\otimes 1_{\rho}\circ R,$$  where $F\in (\bar\rho,\bar\rho)$ is a
(not necessarily orthogonal) projection.  Since $\Cal T$ is closed
under subobjects there is a corresponding subobject $\bar\rho'$ of
$\bar\rho$ and hence an $\Bar X\in (\bar\rho',\bar\rho)$ and  $\bar
Y\in (\bar\rho,\bar\rho')$ such that $\bar Y\circ \bar
X=1_{\bar\rho'}$ and $\bar X\circ \bar Y=F$. Define $$R'=\bar
Y\otimes W^*\circ R;\quad \bar R'=W^*\otimes \bar X^*\circ \bar R.$$
A computation now shows that $R'$ and $\bar R'$ solve the conjugate
equations for $\rho'$, completing the proof. \hfill$\square$\medskip

We define a left inverse $\varphi$ for an object $\rho$ in a strict
tensor $C^*$--category $\Cal T$ to be a set $$\varphi_{\sigma,
\tau}: (\rho \sigma, \rho \tau) \to (\sigma, \tau)$$ for objects
$\sigma, \tau$ of $\Cal T$ of linear mappings which are natural  in
$\sigma$ and $\tau,$ i.e. given $S \in (\sigma, \sigma')$ and $T \in
(\tau, \tau')$ we have $$\varphi_{\sigma', \tau'} (1_\rho \otimes T
\circ X \circ 1_{\rho} \otimes S^*) = T \circ \varphi_{\sigma, \tau}
(X) \circ S^*, \ X \in (\rho \sigma, \rho \tau),$$ and furthermore
satisfy $$\varphi_{\sigma \pi, \tau \pi} (X \otimes 1_\pi) = 
\varphi_{\sigma, \tau} (X)  \otimes 1_\pi, \ X \in (\rho \sigma,
\rho \tau)$$ for each object $\pi $ of $\Cal T.$ We will say that
$\varphi$ is  positive if  $\varphi_{\sigma, \sigma}$ is positive
for each $\sigma$ and normalized if $$\varphi_{\iota, \iota}(1_\rho)
= 1_\iota$$ When $\varphi $ is positive, we say that $\varphi$ is
faithful if  $\varphi_{\sigma, \sigma}$ is faithful for each object
$\sigma.$\par By way of explanation of the name left inverse, this
term has already been used in connection with endomorphisms of a
$C^*$--algebra $\Cal A$. When $\varphi$ is a left inverse of an
endomorphism $\rho$ and we define for each pair of endomorphisms
$\sigma, \tau$ of $\Cal A$ and $X \in (\rho\sigma,\rho\tau)$
$$\varphi_{\sigma,\tau}(X)=\varphi(X) \in (\sigma,\tau),$$ then we
get a left inverse for the object $\rho$ of End$\Cal A$. We may 
think of a left inverse for an object in a tensor $C^*$--category as
being  some kind of virtual object: $\varphi_{\sigma,\tau}$
corresponds to tensoring  on the left by this virtual object looked
at as a map from  $(\rho\sigma,\rho\tau)$ to $(\sigma,\tau)$. If
$\rho $ has a conjugate then we may always construct non--zero left
inverses; for suppose $R, S \in (\iota, \bar \rho \rho)$ and define
$$\varphi_{\sigma, \tau} (X): = S^* \otimes 1_\tau \circ 1_{\bar
\rho} \otimes X \circ R \otimes 1_\sigma, \quad X \in (\rho \sigma,
\rho \tau),$$ then $\varphi$ is obviously a left inverse which is
positive if $S=R$ and is normalized if $S^* \circ R = 1_\iota.$ Note
that when $S=R$ defines a conjugate for $\rho$ then we get a
faithful left inverse. \bigskip \proclaim{Lemma 2.5} Suppose $\rho$
has a conjugate $\bar \rho$ then every left  inverse is of the above
form. More precisely if $R \in (\iota, \bar \rho  \rho)$ defines a
conjugate for $\rho$ then left inverses for $\rho$ are in $1-1$
correspondence with the elements $Y$ of $(\bar \rho, \bar \rho)$  by
taking $S^*= R^* \circ Y \otimes 1_\rho$ in the above formula and a
left  inverse is positive if and only if $Y$ is positive.\par
\medskip \endproclaim\par\demo{Proof} Let $X \in (\rho \sigma, \rho
\tau)$ then trivially, $$\align  X &= 1_\rho \otimes R^*\otimes
1_\tau \circ \bar R \otimes 1_{\rho \tau}  \circ X \circ \bar R^*
\otimes 1_{\rho \sigma} \circ 1_{\rho} \otimes R \otimes 1_{\sigma}\\
{}&= 1_\rho \otimes R^* \otimes 1_\tau \circ 1_{\rho\bar\rho}
\otimes X \circ (\bar R \circ \bar R^* )  \otimes  1_{\rho\sigma}
\circ 1_\rho \otimes R\otimes 1_\sigma. \endalign$$ Hence, if
$\varphi$ is a left inverse for $\rho$ we have  $$\varphi_{\sigma,
\tau} (X) = R^*  \otimes 1_{\tau} \circ 1_{\bar \rho} \otimes X
\circ \varphi_{\bar \rho, \bar \rho} (\bar R \circ \bar R^*)\otimes 
1_{\rho \sigma} \circ R \otimes 1_\sigma.$$ Since $\varphi_{\bar
\rho, \bar \rho} (\bar R \circ \bar R^*) \in (\bar \rho,\bar \rho)$
we see that every left inverse is of the stated form and is positive
if and only if $\varphi_{\bar \rho, \bar \rho}  (\bar R \circ \bar
R^*)$ is positive.\hfill$\square$ \medskip A left inverse $\varphi$
is thus uniquely determined by $\varphi_{\bar\rho,\bar\rho} (\bar
R\circ\bar R^*)$. Notice that, alternatively, every left inverse for
$\rho$ is of the form $$\varphi_{\sigma, \tau} (X) = S^*  \otimes
1_{\tau} \circ 1_{\bar \rho} \otimes X \circ R \otimes 1_{\sigma},
\quad X \in {(\rho \sigma, \rho \tau)},$$ where $S$ is a uniquely
determined but arbitrary element of ${(\iota, \bar \rho \rho)}.$
Hence by Lemma 2.1, a left inverse $\varphi$ is also uniquely
determined by $\varphi_{\iota,\iota}$. \medskip

Before giving the basic inequality on left inverses associated with
solutions  of the conjugate equations, we need a preliminary lemma.
\proclaim{Lemma 2.6} Let $S\in (\iota,\sigma)$ then $$S\circ S^*
\leq (S^*\circ S)\otimes 1_\sigma;\quad S\circ S^* \leq 1_\sigma
\otimes (S^*\circ S)$$.\endproclaim \demo{Proof.} $S^*\circ S\circ
S^*=(S^*\circ S)\otimes 1_\iota\circ 1_\iota\otimes S^* =S^*\circ
(S^*\circ S)\otimes 1_\sigma$. We now have a positive element $X
:=S\circ S^*$  and a positive central element $Z := (S^*\circ
S)\otimes 1_\sigma$ in the C$^*$-algebra  $(\sigma,\sigma)$ such
that $X^2=XZ$. Hence $X \leq Z$, since $\pi(X) \leq \pi(Z)$  for
each irreducible representation $\pi$. The second inequality is
proved similarly.\hfill$\square$ \medskip \proclaim{Lemma 2.7} Let
$R$ and $\bar R$ solve the conjugate equations for  $\rho$ and let
$\varphi$ be the associated left inverse, then  $$X^*\circ X \leq
\bar R^*\circ \bar R\otimes 1_{\rho\sigma}\circ 1_\rho\otimes
\varphi_{\sigma,\sigma}(X^*\circ X),\quad X\in
(\rho\sigma,\rho\tau).$$ \endproclaim \demo{Proof}  Set $Y:= 1_{\bar
\rho}  \otimes X\circ R \otimes 1_{\sigma}$ then $\bar R^* \otimes
1_{\rho \tau}
 \circ 1_\rho \otimes Y = X.$ Hence $X^* \circ X = 1_\rho \otimes
Y^* \circ (\bar R \circ \bar R^*) \otimes 1_{\rho \tau}\circ 1_\rho
\otimes Y.$  Applying Lemma 2.6 with $S=\bar R\in
(\iota,\rho\bar\rho)$ $$X^* \circ X \leq (\bar R^* \circ \bar
R)\otimes 1_{\rho\sigma}\circ  1_\rho \otimes (Y^* \circ Y) =  (\bar
R^* \circ \bar R)\otimes 1_{\rho\sigma}\circ 1_\rho \otimes \varphi_
{\sigma,\sigma}(X^* \circ X).$$\hfill$\square$ \medskip Actually
this inequality is the best possible in the sense that $\bar
R^*\circ\bar R$ is the smallest element of $(\iota,\iota)$ for which
this inequality holds.  This is seen by taking $X^* \circ X=
 \bar R \circ \bar R^* \in (\rho \bar \rho, \rho \bar \rho).$
 
 As a consequence of this inequality there are various finiteness
properties of the objects 
 $\rho$ of $\Cal T_f$. For example, the image of $(\rho,\rho)$ under
a representation is finite dimensional whenever the image of the
centre is finite dimensional.    \medskip

Now just as we have introduced the notion of a left inverse of an
object $\rho$  in a tensor $C^*$-category, so we can introduce the
notion of a right inverse replacing tensoring on the left by $\rho$
by tensoring on the right. Thus a  right inverse $\psi$ for $\rho$
is a set of mappings $$\psi_{\sigma,\tau}: (\sigma\rho,\tau\rho) \to
(\sigma,\tau),$$ natural in $\sigma$ and $\tau$, i.e. such that
given $S\in (\sigma,\sigma')$  and $T\in (\tau,\tau')$ we have
$$\psi_{\sigma',\tau'}(T\otimes 1_\rho\circ X\circ S^*\otimes
1_\rho)=T\circ \psi_{\sigma,\tau}(X)\circ S^*, \quad X\in
(\sigma\rho,\tau\rho),$$ and satisfying
$$\psi_{\pi\sigma,\pi\tau}(1_\pi \otimes X)=1_\pi
\otimes\psi_{\sigma,\tau}(X),  \quad X\in (\sigma\rho,\tau\rho).$$
Since this notion results from that of left inverse by dualizing
with  respect to $\otimes$, all the above results have their
counterparts for  right inverses. In particular, if $R, \bar R$
solve the conjugate equations for  $\rho$ there is an associated
right inverse  given by $$\psi_{\sigma,\tau}(X)=1_\tau \otimes \bar
R^*\circ X\otimes 1_{\bar\rho} \circ 1_\sigma\otimes \bar R, \quad
X\in (\sigma\rho,\tau\rho).$$

In the same way, we have another $(\iota,\iota)$--valued scalar
product  $(\cdot,\cdot)_{r}$ on $\Cal T_f$ associated with a choice
$\rho\mapsto R_\rho$  for each object $\rho$ of $\Cal T_f$, namely
if $S,T\in (\rho,\rho')$, $$(S,T)_{r}:=\bar R_\rho{}^*\circ
(S^*\circ T)\otimes 1_{\bar\rho}\circ \bar R_\rho.$$

Objects $\rho$ of the form $\bar\sigma\sigma$ will be termed {\it
canonical}  since they are an abstract analogue of the canonical
endomorphisms introduced in  [\ref(Long4)]. They have some special
properties which are well known in the theory of  adjoint functors.
In fact, let $R\in (\iota,\bar\sigma\sigma)$ and $\bar R\in 
(\iota,\sigma\bar\sigma)$ satisfy the conjugate equations then $R\in
(\iota,\rho)$  and $M:=1_{\bar\sigma}\otimes\bar R^*\otimes
1_\sigma\in (\rho^2,\rho).$ We have  $$M\circ M\otimes 1_\rho=M\circ
1_\rho\otimes M;\quad M\circ R\otimes 1_\rho=1_\rho;  \quad M\circ
1_\rho\otimes R=1_\rho$$ and since these are the associative and
unit laws, we can, following the standard  practice in category
theory, say that they make $\rho$ into an algebra in the  tensor
$C^*$--category $\Cal T$ in question. Since ${}^*$ is a duality for
$\Cal T$,  this is equivalent to saying that $M^*$ and $R^*$ make
$\rho$ into a cogebra  in $\Cal T$. We also have a further relation
$$1_\rho\otimes M\circ M^*\otimes 1_\rho=M^*\circ M$$ which is not
part of the theory of adjoint functors as it involves both $M$ and 
$M^*$. We shall, however, see in Section 7, that, in our context, it
is actually a  consequence of the other equations.

Note that if $F : \Cal T_1 \to \Cal T_2$ is a strict tensor
$^*$-functor then $F$ maps solutions of the conjugate equations into
solutions for the image object and intertwines the  associated left
and right inverses. Hence $F$ maps conjugates  onto conjugates,
canonical objects onto canonical objects and induces a strict tensor
$^*$-functor from $\Cal T_{1f}$ to  $\Cal T_{2f}$.

\beginsection{ 3. Dimension and Trace}

We turn now to the related notions of dimension and trace which are
central to this  paper. In fact, in view of their importance, we
want, throughout this section, to restrict ourselves to developing 
this theory in the  important special case that  $(\iota,\iota)=
\Bbb C$, i.e. that $\iota$ is irreducible. In the language of von
Neumann algebras this corresponds to defining the index just for the
case of subfactors. 

We note the following simple corollary of Lemma 2.2. \proclaim{Lemma
3.1} Two irreducible objects $\rho$ and $\bar \rho$ are conjugate if
and only if there are $R \in (\iota, \bar \rho \rho)$ and $\bar R
\in (\iota, \rho \bar \rho)$ such that $\bar R^* \otimes 1_\rho
\circ 1_\rho \otimes R \not= 0.$ \endproclaim The next result lends
support to the interpretation of the objects of $\Cal T_f$  as being
finite dimensional. \proclaim{Lemma 3.2} When $\iota$ is irreducible
and $\rho$ is an object of $\Cal T_f$  then $(\rho,\rho)$ is finite
dimensional. In particular if $\Cal T$ has subobjects, $\rho$ is
just a finite direct sum of irreducibles. \endproclaim \demo{Proof}
Let $X_1, X_2, \dots , X_n$ be positive elements of $(\rho,\rho)$ of
norm 1 with $\sum_i X_i \leq 1_\rho$, then Lemma 2.7 shows that
$n\leq \bar R^*\circ \bar  R\circ R^*\circ R$, so that $(\rho,\rho)$
is finite dimensional.\hfill$\square$ \medskip

We study further the set of solutions of the conjugate equations for
given $\rho$ and $\bar\rho$ defining equivalence classes and
normalizing solutions. On this basis, we can define a well behaved
notion of dimension for the corresponding left inverses. We  shall
then show that there exists a privileged class of normalized
solutions called the  standard solutions. This privileged class will
be characterized in a variety of ways and serves  to define a
dimension for the objects themselves.

The pertinent notion of equivalence is characterized in the
following lemma. \proclaim{Lemma 3.3} Two solutions $R,\,\bar R$ and
$R',\,\bar R'$ of the conjugate  equations for $\rho,\,\bar\rho$
define the same left inverse for $\rho$ if and only if  there is a
unitary $U\in (\bar\rho,\bar\rho)$ such that $R'=U\otimes
1_\rho\circ R$. \endproclaim \demo{Proof} Since a left inverse
$\varphi$ is uniquely determined by $\varphi_{\iota, \iota}, R$ and
$R'=U\otimes 1_\rho\circ R$ define the same left inverse if and only
if  $$R^*\circ(U^*\circ U)\otimes X\circ R=R^*\circ
1_{\bar\rho}\otimes X\circ R,\,\,X\in (\rho,\rho),$$ i.e. if and
only if $$R^*\circ 1_{\bar\rho}\otimes X\circ(U^*\circ U)^{\bullet
*}\circ R,\,\,X\in (\rho,\rho),$$ where ${}^\bullet$ is the map from
$(\bar\rho,\bar\rho)$ to $(\rho,\rho)$ defined by $R$. Since the
left inverse is faithful, we conclude that $(U^*\circ U)^{\bullet
*}=1_\rho$  so by Lemma 2.3, $U^*\circ U=1_{\bar\rho}$ as
required.\hfill$\square$\medskip

We can, of course, equally well define an equivalence class by
requiring $\bar R$ and $\bar R'$ to define the same left inverse for
$\bar\rho$. However, this will define the same equivalence class if
and only if the map ${}^\bullet$ commutes with the adjoint hence, by
Lemma 2.3c if and only if $\varphi_{\iota,\iota}$ is tracial.

On the basis of Lemma 3.3, we make the following definition which
gives an analogue of the square root of the index relative to a
conditional expectation in the theory of subfactors [\ref(Kos1)].
\proclaim{Definition 3.4} Let $\varphi$ be the left inverse of $\rho$
defined by a solution $R,\,\bar R$ of the conjugate equations, then
we define the dimension $d(\varphi)$ of $\varphi$ by
$$d(\varphi)=\|R\|\,\|\bar R\|.$$ \endproclaim\medskip

Note that if $\mu > 0$, then $d(\mu\varphi)=d(\varphi),$ so the
dimension is independent of the normalization of $\varphi$. For this
reason, it is convenient to restrict oneself to {\it normalized}
solutions $R,\,\bar R$, i.e. those for which $\|R\|=\|\bar R\|$. In
this case, we have $$d(\varphi)=R^*\circ R=\bar R^*\circ \bar R=
d(\bar\varphi),$$ where $\bar\varphi$ denotes the left inverse of
$\bar\rho$ defined by $\bar R$. \proclaim{Lemma 3.5} Let $R,\,\bar
R$ be normalized solutions of the conjugate equations for 
$\rho,\,\bar\rho$ and $\varphi$ the corresponding left inverse. Then
$d(\varphi)=1$ if  and only if $R$ is unitary.\endproclaim
\demo{Proof} Clearly, if $R$ is unitary, $d(\varphi)=1$. Conversely,
if $d(\varphi)=1$ then computing $X^*\circ X$ for $X=\bar R\otimes
1_\rho - 1_\rho\otimes R$, we see that  $\bar R\otimes
1_\rho=1_\rho\otimes R$. Similarly $R\otimes
1_{\bar\rho}=1_{\bar\rho}\otimes\bar R.$ Thus $$R\circ
R^*=1_{\bar\rho\rho}\otimes R^*\circ R\otimes
1_{\bar\rho\rho}=1_{\bar\rho}\otimes\bar R ^*\otimes 1_\rho\circ
1_{\bar\rho}\otimes\bar R\otimes 1_\rho=1_{\bar\rho\rho},$$ and $R$
is unitary.\hfill$\square$\medskip

We now show that our notions behave well under direct sums and
tensor products. Let  $W_i \in (\rho_i,\rho)$ and $\bar W_i \in
(\bar \rho_i,  \bar \rho)$ be isometries expressing $\rho$ and
$\bar\rho$ as  a direct sum of irreducibles and let  $R_i \in
(\iota, \bar \rho_i \rho_i)$ and $\bar R_i  \in
(\iota,\rho_i\bar\rho_i)$ be solutions of the conjugate equations
then  $$R = \sum_i \bar W_i \otimes W_i \circ R_i,\quad  \bar R =
\sum_i W_i \otimes \bar W_i \circ \bar R_i$$  solve the conjugate
equations for $\rho,\,\bar\rho$ and are said to be a direct sum of
the original solutions $R_i,\,\bar R_i$. They are normalized if the
original solutions were normalized. This notion is compatible with
equivalence classes and we have:  $$R^* \circ R = \sum_i R^*_i \circ
R_i,$$ and hence for the corresponding left inverses:
$$d(\varphi)=\sum_i d(\varphi_i).$$
 Turning to the product, a trivial computation gives the following
result. \proclaim{Lemma 3.6} Let $R_i\in(\iota,\bar\rho_i\rho_i)$ and
$\bar R_i\in (\iota,\rho_i\bar\rho_i)$ define a conjugate for
$\rho_i, i=1,2,$ then $$R :=1_{\bar\rho_2}\otimes R_1\otimes
1_{\rho_2}\circ R_2;\quad \bar R := 1_{\rho_1}\otimes \bar
R_2\otimes 1_{\bar\rho_1}\circ\bar R_1$$ define a conjugate for
$\rho_1\rho_2$. If the $R_i$ are normalized, so is $R$.
\endproclaim\medskip

It is clear that this notion of product is compatible with
equivalence classes and that  for the corresponding left inverses we
have, with an obvious notation,
$$d(\varphi)=d(\varphi_1)d(\varphi_2).$$

If $\rho$ is a zero object, i.e. if $(\rho,\rho)=0$ then
$d(\varphi)=0.$ When $R,\bar R$ define a conjugate for a non--zero
object $\rho$, the  defining equations show that $$ 1 \leq \|\bar
R^* \otimes 1_{\rho}\|\,\|1_{\rho} \otimes R\| \leq \|\bar R\|\,\|
R\|,$$ so that $d(\varphi)\geq 1$.

We now establish a connection with the work of Jones [\ref(Jones)].
If $R \in (\iota,\bar \rho \rho)$ and $\bar R \in
(\iota,\rho\bar\rho)$ satisfy  the conjugate equations, we let  $E$
and $\bar E$  be the range projections of $R$  and $\bar R$, 
$$E:=\|R\|^{-2}R\circ R^*;\quad\bar E:=\|\bar R\|^{-2}\bar
R\circ\bar R^*.$$ These will be called the associated {\it Jones
projections\/}. Note that they are unique up to a unitary
equivalence within an equivalence class of solutions, where the
unitary in question can be chosen from the subalgebras
$(\bar\rho,\bar\rho)\otimes 1_\rho$ and  $1_\rho\otimes
(\bar\rho,\bar\rho)$, respectively. Since $$1_\rho\otimes(R\circ
R^*)\circ(\bar R\circ\bar R^*)\otimes 1_\rho\circ 1_\rho\otimes 
(R\circ R^*)=1_\rho\otimes(R\circ R^*),$$ we have $$1_\rho\otimes
E\circ\bar E\otimes 1_\rho\circ 1_\rho\otimes
E=\|R\|^{-2}\cdot\|\bar  R\|^{-2}\, 1_\rho\otimes E$$ and similarly
$$1_{\bar\rho}\otimes\bar E\circ E\otimes 1_{\bar\rho}\circ
1_{\bar\rho}\otimes\bar  E=\|R\|^{-2}\,\cdot\|\bar
R\|^{-2}\,1_{\bar\rho}\otimes\bar E.$$ These are the Jones relations
and by looking at the associated representations of  the braid
group, we may conclude as in [\ref(Jones)] that
$$\|R\|^{-2}\cdot\|\bar R\|^{-2}\,\in \{4\text{cos}^2{\pi\over
k}:k\in \Bbb N,\,\,k\geq 3\}\cup  [4,\infty).$$ 

We note that if $\varphi$ is a normalized faithful left inverse,
then by Lemma 2.7, $d(\varphi)$ is the smallest constant $c$ such
that, for all objects $\sigma,\,\tau$, we  have $$X^*\circ X\leq
c^2\,\,1_\rho\otimes\varphi_{\sigma,\sigma}(X^*\circ X),\quad X\in
(\rho\sigma,\rho\tau).$$ In the special case that $\Cal
T=\text{End}\Cal A$ for a unital $C^*$--algebra $\Cal A$, we can
regard $\varphi$ as a left inverse for the endomorphism $\rho$ by
setting  $$\varphi(X):=(R^*R)^{-1}\,R^*\bar\rho(X)R$$ and $\Cal
E:=\rho\varphi$ is a conditional expectation from $\Cal A$ onto
$\rho(\Cal A)$. Setting $Y=\bar\rho(X)R,$ cf. Lemma 2.6, we have
$$X^*X=\rho(Y^*)\bar R\bar R^*\rho(Y)\leq\rho(Y^*Y)\bar R^*\bar
R=\Cal E(X^*X)\,R^*R\, \bar R^*\bar R.$$ Hence $d(\varphi)^{-2}$ is
the Pimsner--Popa bound of $\Cal E$, in other words, Ind$\Cal
E=d(\varphi)^2.$

If $\varphi$ is a linear mapping of $\Cal A$ into $\Cal A$ such that
$$\varphi(\rho(A)B\rho(C))=A\varphi(B)C,\quad A,B,C\in \Cal A,$$
then for each  pair $\sigma,\tau$ of endomorphisms of $\Cal A$,
$$\varphi(X)\sigma(A)=\varphi(X\rho\sigma(A))=\varphi(\rho\tau(A)X)=\tau(A)\varphi(X)
,\quad X\in (\rho\sigma,\rho\tau).$$ Thus we get an associated left
inverse of $\rho$ in End$\Cal A$ by taking  $\varphi_{\sigma,\tau}$
to be the restriction of $\varphi$ to $(\rho\sigma,\rho\tau)$.

Thus for each object $\rho$, there are two sets of Jones
projections, those of the form $E$  associated with faithful right
inverses, and those of the form $\bar E$ associated with faithful
left inverses. For each faithful left or right inverse, there is a
unique unitary  equivalence class of Jones projections. Each Jones
projection $E$ determines a solution  of the conjugate equations up
to a scalar.

The vital step towards introducing a dimension for the objects of
$\Cal T_f$ rather than  their left inverses is to choose a special
equivalence class of solutions of the conjugate equations. We first
take this step for irreducible objects by picking the unique
equivalence class of normalized solutions.

In the general case, we write $\rho$ as a direct sum of irreducibles
and say  that $R,\,\bar R$ is a {\it standard} solution if $R$ is of
the form $$R = \sum_i \bar W_i \otimes W_i \circ R_i,$$ where $R_i
\in (\iota, \bar \rho_i \rho_i)$ is a normalized solution for an 
irreducible object $\rho_i$  and  $W_i \in (\rho_i,\rho)$ and $\bar
W_i \in (\bar \rho_i,  \bar \rho)$ are isometries expressing $\rho$
and $\bar\rho$ as  a direct sum of irreducibles. As we know we then
have  $$\bar R = \sum_i W_i \otimes \bar W_i \circ \bar R_i$$  and
$R,\,\bar R$ is a normalized solution.  With an obvious notation,
the corresponding left inverses are related by 
$$\varphi(X)=\sum_i\varphi_i(W_i^*\circ X\circ W_i).$$  It is easy
to see that $\varphi$ is independent of the choice of the $W_i$  so
that we get in this way a single well determined equivalence class
of  solutions of the conjugate equations. This left inverse will be
called the {\it standard} left inverse, a  notion introduced in
[\ref(DR1)] in the context of a symmetric  tensor $C^*$--category. 

\proclaim {Lemma 3.7} If $\rho \mapsto R_\rho$ is standard then the
corresponding  scalar product $(\cdot.\cdot)$, defined as in Section
2, is independent of the choice of  standard solutions and is
tracial. \endproclaim \demo{Proof} The first statement follows from
the uniqueness of the standard left inverse; to prove the second,
we  use Lemma 2.3c. Now when the $R_\rho$ are standard, $R_{\rho}$
and $R_{\rho'}$ have the form: $$R_\rho = \sum_i \bar W_i \otimes
W_i \circ R_{\rho_i}; \quad R_{\rho'}= \sum_j \bar W'_j \otimes 
W'_j \circ R_{\rho'_j};$$ where, without loss of generality, we may
suppose that, when  $\rho_i$\ and\ $\rho'_j$ are equivalent, they
are actually equal. If $T\in(\rho, \rho')$ then $$T^\bullet\otimes
1_\rho\circ \sum_i\bar W_i\otimes W_i\circ R_{\rho_i} =
1_{\bar\rho'}\otimes T^*\circ\sum_j\bar W'_j\otimes W'_j\circ
R_{\rho'_j}$$ Multiplying on the left by $\bar{W'}^*_j \otimes 
W^{*}_i,$ we conclude that $$\bar{W'}^*_j \circ T^\bullet \circ \bar
W_i = ({W'}^*_j \circ T \circ W_i)^\bullet= W^*_i \circ T^* \circ W'_j
= \bar{W'}*_j \circ  T^{*\bullet *} \circ \bar W_i.$$ Thus
$T^{\bullet *} = T^{*\bullet} $ as required.\hfill$\square$ \medskip

There is a notion of trace is implicit in the above scalar product.
To make it explicit, given $T \in (\rho, \rho),$ we  write $tr (T):=
(1_\rho, T).$ It is easy to  compute the trace of a projection $E.$ 
In fact, if we write $E= W \circ W^*,$  where $W \in (\sigma, \rho)$
is an isometry then  $$tr (E)= (1_\rho, E)= (W^*, W^*)= (W, W)= tr
(1_\sigma).$$  Thus the trace of $E$ is just the dimension of the
standard left inverse of $\sigma$. This shows that the trace on our
full subcategory is entirely determined  by the dimensions of the
standard left inverses.

We now define $d(\rho)$ for an object $\rho$ to be the dimension of
its standard left inverse $\varphi$.  We then have
$d(\rho)=d(\bar\rho)$ so that $d(\rho)$ depends only on the unitary
equivalence class of $\rho$ and is additive on direct sums. Since 
it is additive, it follows  from Theorem 2.4 that an object of
dimension $<2$ is necessarily irreducible. 

Let us now look at a few simple examples. If we take our strict
tensor C$^*$--category to be a category of finite  dimensional
Hilbert spaces then every object is a finite direct sum of
1--dimensional objects and the dimension is the usual one. More
generally, if we look at a category of finite-dimensional,
continuous, unitary representations of a compact group then we can
take $$R\lambda = \lambda \sum_i \bar e_i \otimes
 e_i; \quad \bar R\lambda = \lambda \sum_i e_i\otimes \bar e_i;
\quad \lambda \in \Bbb C;$$ where $e_i$ is an orthonormal basis of
the underlying Hilbert space and  $\bar e_i$ is the conjugate basis
of the conjugate Hilbert space. Thus we again have the usual notion
of dimension.

Turning to compact quantum groups in the sense of Woronowicz, we
again get a strict  tensor C$^*$--category with $(\iota,\iota)=\Bbb
C$ by taking unitary representations on  objects of a strict tensor
C$^*$--category of Hilbert spaces. For a unitary representation  on
a finite--dimensional Hilbert space, we get solutions of the
conjugate equations by taking  $$R\lambda = \lambda \sum_i f_i
\otimes
 J^{-1}f_i; \quad \bar R\lambda = \lambda \sum_i e_i\otimes Je_i;
\quad \lambda \in \Bbb C;$$ where $e_i$ is again an orthonormal
basis of the underlying Hilbert space and  $f_i$ is an orthonormal
basis of the Hilbert space of the conjugate[\ref(Wor1)]. However $J$
is  now just an invertible antilinear intertwiner and is not
necessarily antiunitary. Since  $$R^*\circ
R=\text{Tr}J^{-1*}J^{-1};\quad \bar R^*\circ \bar R=\text{Tr}J^*J,$$
the dimension will no longer coincide with the Hilbert space
dimension in general.

Let us compute a simple example, the case of $S U_q (2),$ with $- 1
< q < 1$. A few remarks will clarify the general structure. It is
known that the fusion rules for the irreducible representations are
the same as in the classical case of the dual of $S U (2)$
[\ref(Wor1)]. The irreducible representations of $S U_q (2)$ are
indexed by $\{\rho_i, i=0, 1,2 \dots \},$ $\,\rho_0 = id$ and
$\rho_1$ the  fundamental representation, and we have
$$\rho_i\otimes\rho_j = \bigoplus\limits^{i+j}_{k=|i-j|} \rho_k.$$
This means the principal graph of tensoring by $\rho_1$ is
$A_\infty.$ Now $d(\rho_0)=1$ and we shall see later that the
dimension is multiplicative  on tensor products. Hence it suffices
to compute $d := d(\rho_1)$ since the other dimensions  can be
computed from the recurrence relation
$$d(\rho_{n+1})=d(\rho_n)d-d(\rho_{n-1}).$$  The task of computing
$d := d(\rho_1)$ is facilitated by the fact that $\rho_1$ is
self-conjugate. It can be easily checked that with the matrix
realization of Eq.$\!$ (1.34)  of [\ref(Wor1)], we can take the
antilinear intertwiner $J$ to be given by  $$Je_1=q^{-1}e_2;\quad
Je_2=e_1;$$ where $e_1,\,e_2$ are the canonical basis of $\Bbb C^2$.
Thus we conclude that  $d(\rho_1)=q+q^{-1}$.

Before being able to show that the dimension is multiplicative, we
need to examine the relationship between  (faithful) left and right
inverses.  There are in fact two  natural ways of passing from a
left inverse to a right inverse.  One way is to take a solution
$R,\,\bar R$ of the conjugate equations and to  pass from the left
inverse defined by $R$ to the right inverse defined by $\bar R$.
Obviously, this procedure depends only on the equivalence class of
the solutions  and is compatible with direct sums and products. The
alternative procedure is to  start with a left inverse $\varphi$ and
look for the right inverse $\psi$ such that 
$$\varphi_{\iota,\iota}=\psi_{\iota,\iota}.$$ This condition is
again compatible with taking direct sums. However, it is less
obvious that it is compatible with taking direct products.

\proclaim{Theorem 3.8} Let $\varphi$ and $\psi$ be left and right
inverses of $\rho$ and  $\varphi'$ and $\psi'$ left and right
inverses of $\rho'$ such that
$$\varphi_{\iota,\iota}=\psi_{\iota,\iota}\quad
\varphi'_{\iota,\iota}=\psi'_{\iota,\iota}$$  then the corresponding
product left and right inverses for $\rho\rho'$, $\varphi'\varphi$
and  $\psi\psi'$ satisfy
$$(\varphi'\varphi)_{\iota,\iota}=(\psi\psi')_{\iota,\iota}.$$
\endproclaim \demo{Proof} Let $X\in (\rho\rho')$, then
$$\varphi'_{\iota,\iota}\varphi_{\rho',\rho'}(X)=\psi'_{\iota,\iota}\varphi_{\rho',\rho'}(X),$$
since $\varphi_{\rho',\rho'}(X)\in (\rho',\rho')$. Now any right
inverse commutes with any left inverse, so the above equals
$$\varphi_{\iota,\iota}\psi'_{\rho,\rho}(X)=\psi_{\iota,\iota}\psi'_{\rho,\rho}(X),$$
again since $\psi'_{\rho,\rho}(X)\in
(\rho,\rho).$\hfill$\square$\medskip

It is clear that the standard left inverse agrees with the standard
right inverse of $\rho$  on $(\rho,\rho)$ because it is trivially
true on irreducibles and the equality is compatible with direct
sums. Alternatively, one may observe that the standard right inverse
may equally well be used to define the dimension and the trace. What
is important is the following  characterization of the standard
solutions. 

\proclaim{Lemma 3.9} Let $R\in (\iota,\bar\rho\rho)$ and $\bar R\in
(\iota,\rho\bar\rho)$  solve the conjugate equations for
$\rho,\bar\rho$ then the corresponding left and right  inverses
$\varphi$ and $\psi$ for $\rho$, $$\varphi_{\sigma,\tau}(X) =
R^*\otimes 1_\tau\circ 1_{\bar\rho}\otimes  X\circ R\otimes
1_\sigma, \quad X \in (\rho\sigma,\rho\tau);$$
$$\psi_{\sigma.\tau}(X) = 1_\tau\otimes\bar R^*\circ X\otimes
1_{\bar\rho} \circ 1_\sigma\otimes\bar R, \quad X \in 
(\sigma\rho,\tau\rho);$$  satisfy
$\varphi_{\iota,\iota}=\psi_{\iota,\iota}$ if and only if R is
standard.

\endproclaim \par\demo{Proof} We know that
$\varphi_{\iota,\iota}=\psi_{\iota,\iota}$ when $R$  is standard.
However any other choice of $R$ and $\bar R$ has the form  $R'=
1_{\bar\rho}\otimes Y\circ R$ and $\bar R'=Y^{*-1}\otimes
1_{\bar\rho} \circ\bar R$, where $Y\in (\rho,\rho)$ is invertible.
Thus the corresponding  left and right inverses $\varphi'$ and
$\psi'$ satisfy $\varphi'_{\iota,\iota} =\psi'_{\iota,\iota}$ if and
only if $$\text{tr}Y^*XY=\text{tr}Y^{-1}XY^{*-1}, \quad X\in
(\rho,\rho).$$ This is the case if and only if
$\text{tr}XYY^*=\text{tr}X(YY^*)^{-1}, X \in (\rho,\rho)$ and taking
$X$ to be an eigenprojection of the positive invertible  operator
$YY^*$, we see that $YY^*=I$ so that $Y$ is unitary and $R'$ is
standard.\hfill$\square$\medskip

The multiplicativity of the dimension now follows as a corollary.

\proclaim{Corollary 3.10}\quad If $\rho_1$ and $\rho_2$ have
conjugates, then \smallskip  \item {$a$)}
$d(\rho_1\rho_2)=d(\rho_1)d(\rho_2)$\par \item {$b$)} If $\varphi_1$
and $\varphi_2$ are standard left inverses for $\rho_1$  and
$\rho_2$, then $\varphi$ is a standard left inverse for
$\rho_1\rho_2$,
 where
$\varphi_{\sigma,\tau}=\varphi_{2\,\sigma,\tau}\varphi_{1\,\rho_2\sigma,
\rho_2\tau}.$\par \item {$c$)} $(S_1\otimes S_2, T_1\otimes
T_2)=(S_1,T_1)(S_2,T_2), \quad  S_i,T_i\in (\rho_i,\rho'_i),\quad
i=1,2.$\endproclaim We have therefore given a simple general proof
of the multiplicativity of the  dimension and hence of the  minimal
index. Of course, there are other proofs known in  different
contexts, e.g. [\ref(KosLong)],\ref(Long3)],  but we must here draw
attention to the simple proof given
 in [\ref(KW)] for inclusions of $C^*$-algebras with centre $\Bbb
CI.$ To recognize the relation with our proof, we suppose that our
category is End$\Cal A$, where  $\Cal A$ is a $C^*$-algebra with
centre $\Bbb CI$. If the endomorphism $\rho$ of  $\Cal A$ has
conjugates defined by $R,\bar R$, we have a corresponding 
conditional expectation $\Cal E$ of $\Cal A$ onto $\rho(\Cal A)$
defined by $$\Cal E(A) := a\rho(R^*\bar\rho(A)R),\quad A\in \Cal
A,$$ where $a=\| R\|^{-2}.$ Now $$a^{-1}\bar R^*\Cal E(\bar RA)=\bar
R^*\rho(R^*\rho(\bar RA)R)=A,$$ so Proposition 2 of [\ref(KW)] tells
us that $\Cal E$ is minimal if and only if $$a^{-1}\bar R^*X\bar
R=c\Cal E(X),\quad X\in(\rho,\rho),$$ for some constant $c > 0$. But
this just means that
$$\psi_{\iota,\iota}(X)=ca^2\varphi_{\iota,\iota}(X),\quad
X\in(\rho,\rho).$$ where $\varphi$ and $\psi$ are the left and right
inverses of $\rho$ defined  by $R,\bar R$. Hence by Lemma 3.7,
$c^{1/4}a^{1/2}R$ is standard and $d(\rho)=c^{1/2}a\|
R\|^2=c^{1/2}.$ So Index$\Cal E=d(\rho)^2$. In particular, we see 
that the proof of multiplicativity in [\ref(KW)] rests on the same
criterion as our  proof.

We now give an alternative method of characterizing standard
solutions and the dimension  based on the concept of minimality and
related to using minimal conditional  expectations in the theory of
the index. \medskip \proclaim{Theorem 3.11}\quad Let $R, \bar R$
define a conjugate $\bar\rho$ for $\rho$ then $$R^*\circ R\circ\bar
R^*\circ\bar R\geq d(\rho)^2,$$ with equality holding if and only if
some multiple of R is standard. \endproclaim \demo{Proof} The
inequality itself results by taking $X=1_\rho$ in Lemma 2.7.  Suppose
$R,\bar R$ are standard, then if $R',\bar R'$ also define 
$\bar\rho$ as a conjugate for $\rho$, $R'=1_{\bar\rho}\otimes X\circ
R$ and  $\bar R'=X^{*-1}\otimes 1_{\bar\rho}\circ\bar R$, where
$X\in(\rho,\rho)$ is  invertible. Hence $${R'}^*\circ R'=R^*\circ
1_{\bar\rho}\otimes(X^*\circ X)\circ R=\text{tr}X^*\circ X,$$ $$\bar
{R'}^*\circ\bar R'=\bar R^*\circ(X^{-1}\circ X^{*-1})\otimes
1_{\bar\rho} \circ\bar R=\text{tr}(X^*\circ X)^{-1}.$$ Now writing
$X^*\circ X=\sum_i \mu_i E_i$ and $(X^*\circ X)^{-1}=\sum_i
\mu_i^{-1} E_i,$ where the $E_i$ are mutually orthogonal projections
in $(\rho,\rho)$, we get $${R'}^*\circ R'\circ\bar{R'}^*\circ\bar R'=
\sum_{i,j}\mu_i\mu_j^{-1}\,\text{tr}(E_i) \text{tr}(E_j).$$ Since
$\mu_i\mu_j^{-1}+\mu_j\mu_i^{-1}\geq 2,$ we deduce that $${R'}^*\circ
R'\circ\bar{R'}^*\circ\bar R'\geq\sum_{i,j}\text{tr}E_i\,\text{tr}E_j
=(\text{tr}1_\rho)^2=d(\rho)^2.$$ Thus we have rederived the
required inequality and we see that equality holds only if each
$\mu_i=1$, i.e. only if $X^*\circ X=I$, i.e. only if $R'$ is
standard.\hfill$\square$ \medskip

Let us now consider the behaviour of the dimension under a functor
in the light  of the above theorem. If $F : \Cal T_1\to\Cal T_2$ is
a strict tensor $^*$-functor
 where the tensor unit is irreducible in $\Cal T_1$ and $\Cal T_2$
then $F(R)$ defines  a conjugate for $F(\rho)$ whenever
$R\in(\iota,\bar\rho\rho)$ defines a conjugate  for $\rho$. It
follows that $$d(F(\rho))\leq d(\rho).$$ In special cases, we have
equality here. This is trivially the case when  $d(\rho)=1$ but it
is also true when $F(\rho)$ is irreducible or, more  generally, when
$F$ induces an isomorphism of $(\rho,\rho)$ and 
$(F(\rho),F(\rho))$. Note that the mapping of
$(\rho\sigma,\rho\tau)$ into  $(F(\rho\sigma),F(\rho\tau))$ is
injective whenever $\rho$ has a conjugate.

More interesting is the fact that equality holds whenever there are
just a finite  number of equivalence classes of irreducibles. Let
$\rho_i, i\in I,$ be a maximal set  of pairwise inequivalent
irreducibles and set $d_i := d(\rho_i)$ then since the dimension is
additive and multiplicative, $$d(\rho)d_i=\sum_j m^\rho_{ij}\,d_j,$$
where $m^\rho_{ij}$ is the dimension of $(\rho_j,\rho\rho_i)$. Thus
$d(\rho)$ is  a Perron-Frobenius eigenvalue of the matrix
$m^\rho_{ij}$. Since a $^*$-functor  preserves direct sums, we
conclude that, whenever there is a unique Perron-Frobenius
eigenvalue, as is the case when $I$ is finite, $F(d(\rho))=d(\rho)$.
In particular, when $I$ is finite there can only be a tensor
$^*$-functor into the category  of Hilbert spaces if the dimension
takes on integral values.

If we consider the matrix $m^\rho_{ij}$ as an operator on
$\ell^2(I)$ then 
$$\lim_{n\to\infty}\text{dim}(\rho^n,\rho^n)^{1\over 2n}=\|m^\rho\|
\leq d(\rho).$$  This may be proved as the corresponding result in
the theory of subfactors [\ref(Popa2)]. We conclude from this that
if $F$ is a faithful tensor $^*$--functor then  $$\|m^\rho\|\leq
\|m^{F(\rho)}\|\leq d(F(\rho))\leq d(\rho).$$ In particular, if
$\|m^\rho\|=d(\rho)$, these inequalities are all equalities. The
question of when we have equality $\|m^\rho\|=d(\rho)$ is also of
interest in  the discussion of amenability given in \S 5. Of course,
by the  Perron--Frobenius Theorem, this equality must hold whenever
the category in  question has a finite set of equivalence classes of
irreducibles. However, there are a few other very elementary but
relevant remarks. Obviously, $m^{\rho\oplus\sigma}=
m^\rho+m^\sigma$. Consequently the equality
$d(\rho\oplus\sigma)=\|m^{\rho\oplus \sigma}\|$ implies both
$d(\rho)=\|m^\rho\|$ and $d(\sigma)=\|m^\sigma\|$. Similarly, 
$m^{\rho\sigma}=m^\sigma m^\rho$ and we can draw the analogous
conclusion provided   neither $\rho$ nor $\sigma$ is a zero object.
Finally, $m^{\bar\rho}=m^{\rho *}$ and  we deduce that
$d(\bar\rho\rho)=\|m^{\bar\rho\rho}\|$ if and only if
$d(\rho)=\|m^\rho\|$.

If we define a dimension function to be a  function from the set of
objects with conjugates to $\Bbb R_+$, depending only  on the
equivalence class and additive on direct sums and multiplicative on
tensor  products and taking equal values on conjugate objects, then
$\rho \mapsto  d(F(\rho))$ is a dimension function, whenever F is a
strict tensor $^*$-functor.  Furthermore, if $f$ is any dimension
function  $f(\rho)$ is a Perron-Frobenius eigenvalue of the matrix
$m^\rho_{ij}$. Hence a  necessary condition for the existence of a 
tensor $^*$-functor into the  category of Hilbert spaces is for
there to be an integral-valued dimension function. Thus the
representation categories considered by Woronowicz in [\ref(Wor1)]
trivially have integral-valued dimension functions even when the
intrinsic dimensions of the  objects themselves are not integral. If
there is to be a  right regular object,  i.e. an object  $\rho$ with
the property that $\rho\sigma$ is equivalent to a finite multiple
of  $\sigma$ for all $\sigma$ with conjugates then we again need an
integral-valued dimension function. The condition even becomes
sufficient  if we are allowed to choose $\rho$ from a suitable
closure of the category under  direct sums.  

In the next section devoted to braided tensor $C^*$-categories,  we
shall meet further conditions which guarantee that $d(F(\rho))=d(
\rho).$ Note, however, that in any case the set of objects $\rho$
such that  $F(d(\rho))=d(\rho)$ is closed under finite direct sums,
conjugates and  subobjects.

\beginsection{ 4. Braided Tensor $C^*$--Categories}

We now consider the question of conjugates in strict symmetric, or
braided tensor $C^*$--categories. We begin by proving that if $\rho$
has a conjugate then a left inverse $\varphi$ for $\rho$ is uniquely
determined by the value it  takes on $\varepsilon (\rho,\rho).$
\proclaim{Lemma 4.1} Let $R\in (\iota,\bar\rho\rho)$ define a
conjugate for  $\rho$ and let $\varepsilon$ be a braiding then
$$\varphi_{\bar\rho,\bar\rho}(\bar R\circ\bar R^*)=
1_{\bar\rho}\otimes \bar R^*\circ\varepsilon(\rho,\bar\rho)\otimes
1_{\bar\rho} \circ\varphi_{\rho,\rho}(\varepsilon(\rho,\rho))\otimes
1_{\bar\rho\bar\rho}\circ \bar R\otimes 1_{\bar\rho.} $$ In
particular, $\varphi$ is uniquely determined by the value it takes on
$\varepsilon(\rho,\rho).$\endproclaim \demo{Proof}
$\varphi_{\bar\rho,\bar\rho}(\bar R\circ\bar R^*)= 
\varphi_{\bar\rho,\bar\rho} (1_{\rho\bar\rho}\otimes\bar R^*
\circ\bar R\otimes 1_{\rho\bar\rho})=1_{\bar\rho}\otimes\bar
R^*\circ\varphi_ {\bar\rho,\iota}(\bar R\otimes 1_\rho)\otimes
1_{\bar\rho}.$ $$\bar R\otimes
1_\rho=\varepsilon(\rho,\rho\bar\rho)\circ 1_\rho\otimes \bar
R=1_\rho\otimes\varepsilon(\rho,\bar\rho)\circ\varepsilon(\rho,\rho)
\otimes 1_{\bar\rho}\circ 1_\rho\otimes\bar R$$ Applying $\varphi$
to this equation and substituting leads to the desired result and
$\varphi$ is uniquely determined by the value it takes on 
$\varepsilon(\rho,\rho)$ since we have seen in the course of proving
Lemma 2.5 that is uniquely determined by the value  it takes on
$\bar R\circ\bar R^*.$\hfill$\square$

To render the statement of results as simple as possible, we shall
suppose in what follows that our category has conjugates. We now
show, using the standard left inverse, that a braiding $\varepsilon$
gives rise to a central element of the category.

\proclaim{Theorem 4.2} Let $\varepsilon$ be a braiding and $\varphi$
the standard left inverse for $\rho$ and set
$$\kappa(\rho):=d(\rho)\varphi_{\rho,\rho}(\varepsilon(\rho,\rho)) $$
then $\rho\mapsto \kappa(\rho)\in (\rho,\rho)$ defines an element of
the  centre of the category i.e. if $T\in (\rho,\sigma),$ $$T\circ
\kappa(\rho)=\kappa(\sigma)\circ T.$$ Furthermore, we have
$$\kappa(\rho_1\rho_2)=\kappa(\rho_1)\otimes
\kappa(\rho_2)\circ\varepsilon(\rho_2,\rho_1)
\circ\varepsilon(\rho_1,\rho_2).$$ so that $\kappa$ is
multiplicative if and only if the braiding is in fact a symmetry.
\endproclaim \demo{Proof} We have $$\kappa(\rho)=\sum_i R_i^*\otimes
1_\rho\circ 1_{\bar\rho_i}\otimes W^*_i\otimes  1_\rho\circ
1_{\bar\rho_i}\otimes\varepsilon (\rho,\rho)\circ 1_{\bar\rho_i}
\otimes W_i\otimes 1_\rho\circ R_i\otimes 1_\rho.$$ where
$R=\sum_i\bar W_i\otimes W_i\circ R_i$ is the usual expression for a
standard arrow $R$ defining a conjugate for $\rho.$ A simple
computation now gives $W^*_j\circ \kappa(\rho)=\kappa(\rho_j)\circ
W^*_j$ and $\kappa(\rho)\circ W_j=W_j\circ \kappa(\rho_j).$ But now
if we pick $T\in (\rho,\sigma)$ and isometries
$V_k\in(\sigma_k,\sigma)$ with $\sum_kV_k\circ V^*_k=1_\sigma$ then
$$ T\circ \kappa(\rho) =\sum_{i,k}V_k\circ V^*_k\circ T\circ
W_i\circ W_i^*\circ \kappa (\rho).$$ Since we may suppose that
$\rho_i $ and $\sigma_k$ are either inequivalent or equal so that
$V^*_k\circ T\circ W_i$ is just a multiple of $1_{\rho_i},$ $\kappa$
is indeed central. The formula for $\kappa(\rho_1\rho_2)$ follows
from Corollary 3.9b and the following computations where, to be
concise, we write $\varepsilon$ for $\varepsilon (\rho_1\rho_2,
\rho_1\rho_2)$ and $\varepsilon_{ij}$ for
$\varepsilon(\rho_i,\rho_j)$ and omit indicating the objects
involved when acting with the left inverses $\varphi_1$ and
$\varphi_2$ of $\rho_1$ and $\rho_2$ $$\varepsilon
=1_{\rho_1}\otimes\varepsilon_{12}\otimes 1_{\rho_2}\circ\varepsilon
_{11}\otimes 1_{\rho_2\rho_2}\circ 1_{\rho_1\rho_1}\otimes
\varepsilon_{22} \circ 1_{\rho_1}\otimes\varepsilon_{21}\otimes
1_{\rho_2}.$$ thus $$ \aligned \varphi_1(\varepsilon)& =
\varepsilon_{12} \otimes 1_{\rho_2} \circ\varphi_1
(\varepsilon_{11})\otimes 1_{\rho_2 \rho_2}\circ 1_{\rho_1} \otimes
\varepsilon_{22}\circ\varepsilon_{21}\otimes 1_{\rho_2}\\ & =
1_{\rho_2}\otimes \varphi_1 (\varepsilon_{11})\otimes 1_{\rho_2}
\circ\varepsilon_{12}\otimes 1_{\rho_2}\circ 1_{\rho_1}\otimes
\varepsilon_{22}\circ \varepsilon_{21}\otimes 1_{\rho_2}\\ & =
1_{\rho_2}\otimes\varphi_1(\varepsilon_{11})\otimes 1_{\rho_2}\circ
\varepsilon (\rho_1\rho_2, \rho_2)\circ \varepsilon_{21}\otimes
1_{\rho_2}\\ & = 1_{\rho_2}\otimes\varphi_1(\varepsilon_{11})\otimes
1_{\rho_2}\circ
1_{\rho_2}\otimes\varepsilon_{21}\circ\varepsilon(\rho_1\rho_2,
\rho_2)\\ & = 1_{\rho_2}\otimes \varphi_1(\varepsilon_{11})\otimes
1_{\rho_2}\circ 1_{\rho_2}\otimes
\varepsilon_{21}\circ\varepsilon_{22}\otimes 1_{\rho_1} \circ
1_{\rho_2}\otimes \varepsilon_{12}\endaligned  $$ Applying
$\varphi_2$ gives $$ \varphi_2\varphi_1(\varepsilon) =
\varphi_1(\varepsilon_{11})\otimes \varphi_2 (\varepsilon_{22})\circ
\varepsilon_{21}\circ\varepsilon_{12}. $$ yielding the required
result.\hfill$\square$ \smallskip Note that since $\rho\mapsto
\kappa(\rho)$ is central the value $\in\Bbb C$ it takes on
irreducible objects depends only on the equivalence class of the 
object. Furthemore, as for any central element. $$ \kappa(\rho) =
\sum_n \mu_n E_n $$ where the $E_n$ are the minimal central
projections in $(\rho, \rho)$ and the eigenvalues $\mu_n$\ of\
$\kappa(\rho)$ are simply the values of $\kappa$ on the 
corresponding irreducibles.\par If we introduce a conjugation
$^\bullet$ then by Lemma 2.3 $\rho\mapsto \kappa(\bar\rho)^\bullet
:= \kappa^\bullet (\rho)$ is again central and obviously when $\rho$
is irreducible the scalar $\kappa^\bullet(\rho)$ is just the complex
conjugate of the scalar $\kappa(\bar\rho).$ Thus this involution on
the centre of our category is independent of the choice of
conjugation.\par We now prove the following result which shows that
a braiding is in a sense automatically compatible with conjugates
and that $\kappa(\rho)$ plays another role. \proclaim{Lemma 4.3} Let
$R\in (\iota,\bar\rho \rho)$ define a conjugate for  $\rho$\ then\
$\varepsilon(\bar\rho, \rho)\circ R$ defines a conjugate for 
$\bar\rho.$ $R$ standard implies $\varepsilon(\bar\rho, \rho)\circ
R$ standard  for all $\rho$ if and only if $\kappa$ is unitary.
Furthermore, when $R$ is standard $$ \aligned
&\varepsilon^{-1*}(\bar\rho, \rho)\circ R = \kappa(\rho)^*\otimes
1_{\bar\rho} \circ \bar R,\\ &\varepsilon (\bar\rho, \rho)\circ R =
\kappa^\bullet(\rho)^{-1*}\otimes 1_{\bar\rho} \circ \bar
R.\endaligned $$\endproclaim \demo{Proof} If $R$ defines a conjugate
then $\varepsilon(\bar\rho, \rho) \circ R\ne 0$ and hence trivially
defines a conjugate for $\bar\rho$\  when\ $\rho$ is irreducible.
Furthermore the question of whether $\varepsilon(\bar\rho,
\rho)\circ R$ defines a conjugate for $\bar\rho$ is independent of
the choice of $R.$ Thus we may choose $R$ to be standard writing $R
= \sum_i\bar W_i\otimes W_i \circ R_i$ as in the definition of 
standard. But then $$ \varepsilon (\bar\rho,\rho)\circ R = \sum_i
W_i \otimes \bar W_i \circ  \varepsilon (\bar \rho_i, \rho_i)\circ
R_i $$ which again obviously defines a conjugate for $\bar\rho.$
Since $\varepsilon^{-1*}(\bar\rho, \rho)\circ R$ defines a conjugate
for $\bar\rho,$ we know that $$ \varepsilon^{-1*}(\bar\rho,
\rho)\circ R = X^* \otimes 1_{\bar\rho}\circ \bar R, $$ where
$X\in(\rho, \rho)$ is invertible. Now $$ \aligned X &= (\bar
R^*\circ X\otimes 1_{\bar\rho})\otimes 1_\rho \circ 1_\rho \otimes R
= (R^*\circ \varepsilon^{-1}(\rho, \bar\rho))\otimes 1_\rho\circ
1_\rho\otimes R\\ &= R^*\otimes 1_\rho\circ
1_{\bar\rho}\otimes\varepsilon (\rho, \rho)\otimes
\varepsilon^{-1}(\rho, \bar\rho\rho)\circ 1_\rho \otimes R =
R^*\otimes 1_\rho  \circ 1_{\bar\rho}\otimes\varepsilon(\rho,
\rho)\circ R\otimes 1_\rho =  \kappa(\rho).\endaligned $$ Now
$\varepsilon(\bar\rho, \rho)^{-1 *}\circ R = \kappa(\rho)^*\otimes
1_{\bar\rho}\circ \bar R$ so $R = \varepsilon(\bar\rho, \rho)^*
\otimes 1_{\bar\rho}\circ \bar R = 1_{\bar\rho}\otimes \kappa(\rho)^*
\circ \varepsilon^*(\rho, \bar\rho)\circ \bar R$ thus reversing the
roles of $\rho$\ and $\bar\rho$ we conclude that $$
\varepsilon^*(\bar\rho, \rho)\circ R = 1_\rho\otimes \kappa (\bar
\rho)^{-1*} \circ\bar R = \kappa^\bullet(\rho)^{-1}\otimes
1_{\bar\rho} \circ \bar R, $$ yielding the required result since
passing from $\varepsilon$\ to\ $\varepsilon^*$ means passing from
$\kappa$\ to\ $\kappa^*.$\par Now $\varepsilon(\bar\rho, \rho)\circ
R$ will be standard whenever $R$ is standard if this is the case for
all irreducibles $\rho.$ But this means that the scalars
$\kappa^\bullet(\rho)^{-1*}$ must be a phase so that $\kappa$ is
unitary. Conversely, if $\kappa$ is unitary and $R$ is standard then
so is $\varepsilon(\bar\rho, \rho)\circ R.$\hfill$\square$\smallskip
Note that the above computations show that if the braiding is
unitary then $\kappa$ is unitary, too.\par Now the definition of
conjugates used for strict symmetric tensor  $C^*$-categories in
connection with the duality theory of compact groups [8] required
that $\bar R = \varepsilon(\bar\rho, \rho)\circ R.$ This condition 
was however relaxed in $\S$7 of [8] where it was simply required
that there exist $R\in(\iota, \bar\rho\rho)$\ and\ $\tilde
R\in(\iota,\rho\bar\rho)$ such that $$ R^*\otimes 1_\rho \circ
1_{\bar\rho}\otimes \varepsilon(\rho, \rho)\circ R \otimes 1_\rho\
\text{and}\ \tilde R^*\otimes 1_{\bar\rho} \circ
1_\rho\otimes\varepsilon (\bar\rho, \bar\rho)\circ \tilde R\otimes
1_{\bar\rho} $$   are invertible. This condition is however
equivalent to requiring that $\rho$ has a conjugate in the sense
employed here. For suppose that $R\in (\iota,  \bar\rho\rho)$
defines a conjugate for $\rho$ and is standard and set $\tilde
R:=\bar R$ then the above expressions equal $\kappa(\rho)$\ and\
$\kappa(\bar\rho)$ respectively and are invertible. On the other
hand it was show in $\S$ 7 of [8] that a new symmetry $\hat
\varepsilon$ can be introduced so  that $\rho$\ and\ $\bar\rho$
becomes conjugates in the stricter sense employed  there and hence,
a fortiori in the sense employed here.\par We next present a
sufficient condition for a normalized left inverse to be standard;
this will be of use shortly. \proclaim{Proposition 4.4} Let
$\varphi$ be a normalized left inverse of an object $\rho$ with
conjugates, then if $\kappa$ is unitary and $\varphi_{\rho, \rho}
(\varepsilon(\rho,\rho))$ is a multiple of $1_{\rho}, \varphi$ is
standard. \endproclaim \demo{Proof} Let $W_i \in(\rho_i, \rho)$ be
isometries with $\rho_i$ irreducible and $\sum_iW_i\circ W_i^* =
1_\rho,$ then there are normalized left inverses $\varphi_i$\ of\
$\rho_i$ such that

$$ \varphi_{i\sigma, \tau}(S)\varphi_{\iota,\iota}(W_i\circ W_i^*) =
\varphi_{\sigma,\tau}(W_i\otimes  1_\tau\circ S\circ W_i^*\otimes
1_\sigma),\quad S\in (\rho_i\sigma,\rho_i\tau). $$

Hence

$$ \varphi_{i\rho_i,
\rho_i}(\varepsilon(\rho_i.\rho_i))\varphi_{\iota,\iota}(W_i\circ
W_i^*) = \varphi_{\rho_i,\rho_i}(W_i\otimes
1_{\rho_i}\circ\varepsilon(\rho_i,\rho_i)\circ W^*_i \otimes
1_{\rho_i})=W^*_i\circ\varphi_{\rho,\rho}(\varepsilon(\rho,
\rho))\circ W_i, $$ and when $\varphi_{\rho, \rho}(\varepsilon(\rho,
\rho)) = \lambda 1_\rho,$ we conclude, since any normalized left
inverse of $\rho_i$ is standard, that

$$ d(\rho_i)\lambda 1_{\rho_i} = \kappa
(\rho_i)\varphi_{\iota,\iota}(W_i\circ W^*_i). $$ If $\kappa$ is
unitary, taking absolute values and summing over $i$ gives
$d(\rho)|\lambda| = 1.$ Hence $\varphi_{\iota,\iota}(W_i\circ W^*_i)
d (\rho) = d(\rho_i)$ and finally

$$ d(\rho)\lambda 1_\rho = d(\rho)\lambda \sum_i W_i \circ W^*_i =
\sum_i W_i \circ \kappa(\rho_i)\circ W^*_i = \kappa (\rho), $$
showing that $\varphi$ is indeed standard.\hfill$\square$\smallskip
Now the trace we have introduced in $\S$ 3 depends only on the
notion of  conjugate. But it does make essential use of the
$^*$-operation and to define it we were forced to restrict ourselves
to standard arrows defining conjugates. However, in a braided tensor
category there is another natural trace which is  not positive, in
general, and depends on the braiding. On the other hand, this trace
can be defined without reference to the $^*$-operation and can be 
defined using arrows defining conjugation without requiring $R$ to
be standard. We set for $S, T\in(\rho, \rho)$ $$ (S, T)_\varepsilon
= \bar R^* \circ \varepsilon (\bar\rho,\rho)\circ 1_{\bar
\rho}\otimes (S^* \circ T)\circ R, $$ where $R\in(\iota, \bar\rho,
\rho)$ defines a conjugate for $\rho.$ We see immediately, that this
expression is independent of the choice of $R.$ It is also easy to
check that it is a trace using the argument of Lemma 4.1.  However,
this can also be seen by relating the above expression to the trace
defined previously. In fact, by Lemma 4.3, $$\gather
\varepsilon^*(\rho, \bar\rho)\circ \bar R = 1_{\bar\rho}\otimes
\kappa (\rho)^{*-1} \circ R,\\ (S, T)_\varepsilon = (S, T \kappa
(\rho)^{-1})\endgather $$ and since $\rho \to \kappa (\rho)^{-1}$ is
in the centre of the category it follows at once that $$ (S,
T)_\varepsilon = (T^*, S^*)_\varepsilon $$ This trace has no reason
in general to be multiplicative however if  $\varepsilon$ is a
symmetry $$ (S\otimes S', T\otimes T')_\varepsilon = (S,
T)_\varepsilon (S', T')_ \varepsilon. $$

Our analysis of the invariant $\kappa$ leads to an interesting
condition  guaranteeing that a strict tensor $^*$-functor preserves
dimensions. \proclaim{Proposition 4.5} Let $F: \Cal I_1\to\Cal I_2$
be a strict tensor $^*$-functor, where $\Cal I_1$ has conjugates and
subobjects and $\Cal I_2$ has a braiding $\varepsilon$ for which
$\kappa$ is unitary. Then if $\rho$ is irreducible and
$\varepsilon(F(\rho), F(\rho))$\ is in $F(\rho^2, \rho^2)$\ then\
$d(F(\rho)) = d(\rho).$\endproclaim \demo{Proof}  Picking
$R\in(\iota, \bar\rho \rho)$ standard defining $\bar\rho$ as a
conjugate for $\rho$ and setting $\sigma := F(\rho),\, \bar\sigma :=
F(\bar\rho),$ we see that the  left inverse $\varphi$\ of\ $\sigma$
defined by $$ \varphi_{\lambda,\nu}(X) := d(\rho)^{-1} F(R)^*\otimes
1_\nu\circ  1_{\bar\sigma}\otimes X\circ F(R)\otimes 1_\lambda,\quad
X\in(\sigma\lambda, \sigma\nu) $$ has $\varphi_{\sigma,
\sigma}(\varepsilon(\sigma, \sigma))$ a scalar $c 1_\sigma$ since
this is the image under $F$ of an element of $(\rho, \rho)$ and
$\rho$ is irreducible. Thus by Lemma 4.4, $\varphi$ is a standard
left  inverse and consequently $F(R)$ is standard so that
$d(F(\rho)) = F(R)^*\circ F(R) = F(R^*\circ R) = d(\rho),$ as
required.\hfill$\square$ 

\beginsection{5. Amenable Objects}

As an application of our results, we construct  C$^*$--algebras
$\Cal O_\rho$  and von Neumann algebras $\Cal M_\rho$ equipped in a
canonical way with an  endomorphism $\hat\rho$ starting from
suitable objects $\rho$ in a  tensor C$^*$--category. In this way,
we get  examples of subfactors $\hat\rho(\Cal M_\rho)\subset \Cal
M_\rho$ and,  making contact with the work of Popa, we introduce the
notion of an amenable object in a tensor C$^*$--category.\par In
this section $\Cal T$ always denotes a {\it strict tensor
$C^*$-category}; we will consider objects with conjugates in the
sense of Section 2. By considering the full $C^*$-subcategory they
generate, we may also assume that all objects have conjugates.
Moreover it is always possible to add objects and  hence assume, as
we shall do, that $\Cal T$ has direct sums and subobjects (see
Appendix).  \par

We now fix an object $\rho$ of the  $C^*$--category $\Cal T$ and
review the construction  of the $C^*$--algebra $\Cal O_\rho$ given
in [\ref(DR1)]. We first form the  graded ${}^*$--algebra 
$$\dot{\Cal O}_\rho\equiv\bigoplus_{k\in\Bbb Z}\dot{\Cal
O}^{(k)}_\rho$$ where $$ \dot{\Cal
O}^{(k)}_\rho=\underset\longrightarrow\to\lim_r (\rho^r, \rho^{r+k})
$$ and the embedding $(\rho^r,\rho^s)\to (\rho^{r+1},\rho^{s+1})$ is
given by tensoring on the right with the identity of $\rho :
X\mapsto X\otimes 1_\rho$, which may, for the purposes of this
paper, be supposed to be injective. The $C^*$--algebra $\Cal O_\rho$
is defined as the  completion of  $\dot{\Cal O}_\rho$ with respect
to the unique $C^*$--norm making the expectation $\epsilon
:\dot{\Cal O}_\rho\to\dot{\Cal O}^{(0)}_\rho$,  $$\epsilon
\equiv\int_{\Bbb T}\alpha_t dt\,,\eqno (4.1)$$ continuous, where
$\alpha :\Bbb T\to\text{Aut} (\Cal O_\rho)$ is the {\it gauge}
action  corresponding to the grading, namely $$\alpha_t(X) =
t^kX,\qquad X\in\dot{\Cal O}^{(k)}_\rho.$$\par We now recall that
$\Cal{O}_\rho$ is equipped with a natural endomorphism $\hat\rho$
determined by $$\hat\rho(X) = 1_\rho \otimes X , \qquad X \in
(\rho^r,\rho^s), \quad r,s \in \Bbb N_0,$$ where we identify arrows
between powers of $\rho$ with their images in $\Cal O_\rho$. 
Furthermore, with this identification, we have 
$$(\rho^r,\rho^s)\subset (\hat\rho^r,\hat\rho^s),\quad r,s\in \Bbb
N_0.$$\par Now let $\varphi$ be a positive left inverse for $\rho$
then if $X\in (\rho\sigma,\rho\tau),$
$$\varphi_{\sigma,\tau}(X)^*\circ \varphi_{\sigma,\tau}(X)\leq
\|\varphi_{\iota,\iota}(1_\rho)\| 
\,\,\varphi_{\sigma,\sigma}(X^*\circ X)\leq
\|\varphi_{\iota,\iota}(1_\rho)\|\,\,\|X\|^2\,\,
\varphi_{\iota,\iota}(1_\rho)\otimes 1_\sigma,$$  so that
$$\|\varphi_{\sigma,\tau}(X)\|\leq\|\varphi_{\iota,\iota}(1_\rho)\|\,\,\|X\|.$$
Now, by definition, a left inverse commutes with tensoring on the
right by $1_\rho$ and it hence induces a linear mapping
$\hat\varphi$ of each $\dot{\Cal O}^{(k)}_\rho$  into itself and
extends by continuity to each $\Cal O^{(k)}_\rho$ and we obviously
have  $$\eqalign{\hat\varphi(S^*) & =\hat\varphi(S)^*,\cr
\hat\varphi(T\hat\rho(S)) & =\hat\varphi(T)S,\cr 
\hat\varphi(\hat\rho(S)T) & =S\hat\varphi(T)}.$$ The argument
leading to Theorem 4.12 of [\ref(DR1)] shows that $\hat\varphi$
extends to a  positive linear mapping of $\Cal O_\rho$, also denoted
by $\hat\varphi$ and satisfying  the above equations. Obviously,
$\hat\varphi(I)$ is the image in $\Cal O_\rho$ of 
$\varphi_{\iota,\iota}(1_\rho)$. If $\varphi$ is associated with a
solution $R,\bar R$ of the conjugate equations for $\rho$, then
Lemma 2.7 shows that $$\|X^*X\|\leq\|\bar R^*\circ \bar
R\|\,\,\|\hat\varphi(X^*X)\|,\quad X\in (\rho^r, \rho^r).$$ The
analogous inequality therefore holds for $X\in \Cal O^{(0)}_\rho$ so
that $\hat\varphi$ is faithful on $\Cal O^{(0)}_\rho$.\par From now
on we suppose that $\varphi$ is normalized, then iterating
$\hat\varphi$ gives an $(\iota,\iota)$--valued state $\omega$ on
$\Cal O^{(0)}_\rho$, i.e. a unit-preserving positive linear mapping 
$\omega : \Cal O^{(0)}_\rho \to (\iota,\iota)$, with the property
that  $$\omega(S)=\hat\varphi^r(S),\quad S\in (\rho^r,\rho^r).$$ By
composing this state with the conditional expectation $\epsilon$, we
extend it to a  gauge--invariant $(\iota,\iota)$--valued state, also
denoted by $\omega$ on $\Cal O_\rho$.\par

At this point we make the simplifying    assumption that, as in
Section 3,  $$(\iota,\iota) = \Bbb C,$$ since the following
construction would otherwise have to depend on the choice of an 
auxiliary faithful state on the Abelian C$^*$--algebra
$(\iota,\iota)$. Letting $(\Cal H_\omega, \xi_\omega, \pi_\omega)$
be the Gelfand-Naimark-Segal triple
 associated with  $\omega$, we denote by $$\Cal M_\varphi \equiv
\pi_\omega ( \Cal O_\rho)^{''}$$ the von Neumann algebra generated
by $\pi_\omega (\Cal O_\rho).$

\proclaim{Lemma 5.1} There is a unique isometry $V$ determined by
$$V\pi_\omega(X)\xi_\omega=\pi_\omega(\hat\rho(X))\xi_\omega,\quad
X\in \Cal O_\rho,$$  and we have 
$$V\pi_\omega(X)=\pi_\omega(\hat\rho(X))V;\quad \hat\varphi(X)=
V^*\pi_\omega(X)V,\quad X\in \Cal O_\rho.$$ Hence $\hat\rho$ extends
to a normal endomorphism of $\Cal M_\varphi$ and $\hat\varphi$ to a
normal left inverse of $\hat\rho$ on $\Cal M_\varphi.$ 
Consequently, $\hat\rho\hat\varphi$ is a normal conditional
expectation of $\Cal M_\varphi$ onto $\hat\rho(\Cal M_\varphi)$. If
$\hat\varphi$ is faithful on  $\Cal O^{(0)}_\rho$, $\pi_\omega$ is
faithful and $\hat\rho\hat\varphi$ is faithful. \endproclaim 

\demo {Proof} The existence and uniqueness of the isometry $V$ is
evident. Since  $\hat\rho$ is an endomorphism,
$V\pi_\omega(X)=\pi_\omega(\hat\rho(X))V.$ Thus, if $X_1,X_2\in \Cal
O_\rho$,  $$(\pi_\omega(X_1)\xi_\omega,V^*\pi_\omega(X)V\pi_\omega
(X_2)\xi_\omega)=(\pi_\omega(X_1)\xi_\omega,\hat\varphi(X)\pi_\omega(X_2)\xi_\omega),$$ 
showing that $\hat\varphi(X)= V^*\pi_\omega(X)V,$ as required. If
$\hat\varphi$ is faithful on  $\Cal O^{(0)}_\rho, \pi_\omega$ is
faithful on $(\rho^r,\rho^r), r\in \Bbb N$ and hence on  $\Cal
O_\rho.$ Since $\omega$ is gauge-invariant, $\pi_\omega$ is faithful
by the uniqueness  of the norm on $\dot{\Cal
O_\rho}$.\hfill$\square$\bigskip

Note that since $\omega$ is gauge-invariant, we have an action of
the gauge group on  $\Cal M_\varphi$ compatible with that on $\Cal
O_\rho$ and that  $$\pi_\omega(\rho^r,\rho^s)\subset
(\hat\rho^r,\hat\rho^s).$$

A case of particular interest is when $\varphi$ is the standard left
inverse of $\rho$ and we will then write $\Cal M_\rho$ in place of
$\Cal M_\varphi$ since the von Neumann  algebra is uniquely
determined by $\rho.$ Furthermore, in this case $\omega\vert\Cal
O^{(0)}_\rho$ is  tracial.

\proclaim{Lemma 5.2} $\omega$ is a faithful state of $\Cal
M_\varphi$. In particular,   $\pi_\omega$ is faithful and
$\xi_\omega$ is a separating vector for $\Cal M_\varphi$. 
\endproclaim

\demo{Proof} We shall show that $\omega$ is a KMS state of $\Cal
O_\rho$ for a certain  1--parameter group of automorphisms. To this
end, let $e^{-H}\in (\rho,\rho)$ be the  positive invertible
operator such that $$\varphi(X)={\text tr}(X\circ e^{-H}),\quad X\in
(\rho,\rho).$$ If $e^{-H_n}$ denotes the $n$--th tensor power of
$e^{-H}$, then $$\varphi^n(X)={\text tr}(X\circ e^{-H_n}),\quad X\in
(\rho^n,\rho^n).$$ Now there is a unique continuous 1--parameter
group $\sigma_t$ of automorphisms of  $\Cal O_\rho$ such that
$$\sigma_t(X)=e^{itH_s}\circ X\circ e^{-itH_r},\quad X\in
(\rho^r,\rho^s).$$ To show that $\omega$ is a KMS state for this
1--parameter group, it suffices to pick $X,Y\in (\rho^r,\rho^s)$ and
compute: $$\omega(Y^*\circ \sigma_t(X))={\text tr}(Y^*\circ
e^{itH_s}\circ X\circ e^{-itH_r}e^{-H_r})$$  and performing an
analytic continuation, we get $$\omega(Y^*\circ
\sigma_{t+i}(X))=\omega(\sigma_t(X)\circ
Y^*).$$\hfill$\square$\medskip

Note that in the special case of $\Cal M_\rho$ our modular
automorphisms $\sigma_t$ are just gauge transformations, more
precisely, we have  $\sigma_t=\tilde\alpha_s$ with 
$s=d(\rho)^{-it}$, where $\tilde\alpha_s$ denotes  the normal
extension of $\alpha_s$ to $\Cal M_\rho$. 

We remark that $t\mapsto\sigma_t$ may also be defined directly in
terms of the category  $\Cal T_\rho$. In fact, a state on a
$C^*$--category can be understood as a state $\omega_\sigma$ on
$(\sigma,\sigma)$ for each object $\sigma$ of the category. When  we
have a faithful state on a $W^*$-category, there is a corresponding
modular group of automorphisms of the category[\ref(GLR)]. $\omega$
defines a faithful state on $\Cal T_\rho$  and the corresponding
modular group is a group of automorphisms of the tensor category 
inducing the automorphisms $\sigma_t$ of $\Cal M_\varphi.$

We define $\Cal M_\varphi^{(k)}$ as the $k$-spectral subspace of
$\tilde\alpha$, namely $\Cal M_\varphi^{(k)}$ is the weak closure of
$\Cal O_\rho^{(k)}$. To avoid certain special cases which will be of
no interest in what follows we assume that $\Cal O_\rho^{(k)}$ is
non-zero for some $k\neq 0$ and that $d(\rho)>1.$

\proclaim{Corollary 5.3} $\Cal M_\rho$ is a factor if  the
homogeneous part $\Cal M_\rho^{(0)}$ is a factor, i.e. if $\omega$
is an extremal trace of  ${\Cal O}_\rho^{(0)}$. In this case $\Cal
M_\rho$ is a factor of type $III_\mu$ with $\mu = d(\rho)^{-n}$
where $n=\text {min}\{k\in\Bbb N |\, \exists\, r\,\text{s.t.}\,
(\rho^r,\rho^{r+k})\neq 0\}$. \endproclaim \demo{proof} The center
$Z(\Cal M_\rho)$ of $\Cal M_\rho$ is pointwise invariant under the
modular group $\sigma$ hence it belongs to $\Cal M_\rho^{(0)}$
showing that if $\Cal M_\rho^{(0)}$ is a factor then $\Cal M_\rho$
is a factor too. Since the centralizer  $\Cal M_\rho^{(0)}$ of
$\omega$ is a factor, the Connes invariant[\ref(Con)] $S(\Cal
M_\rho)$ is equal to the exponential of the spectrum of the modular
group $\sigma$ of $\omega$, proving the result.$\hfill
\square$\medskip

\proclaim{Lemma 5.4} Let $X\in \Cal M_\varphi$ then 
$\hat\varphi^n(X)\rightarrow \omega(X)I$ in the  s--topology as
$n\rightarrow\infty$. Consequently, the following conditions are 
equivalent: \item{$a)$} $\hat\varphi(X)=X$, \item{$b)$}
$\hat\rho(X)=X$, \item{$c)$} $X\in (\iota,\iota).$

Furthermore,  $\hat\rho^n(X)\rightarrow\omega (X)I$ in the
$\sigma$--topology as $n\rightarrow\infty$. Thus $\omega$ is the
unique
 normal $\hat\rho$--invariant state and $\hat\varphi$ is the unique
left inverse of 
 $\hat\rho$ having a normal invariant state. More generally, if
$X,\,Y,\,A\in \Cal M_\varphi$ and $A$ has grade $k\geq 0$, then
$$\hat\rho^n(X)A\hat\rho^n(Y)\rightarrow\omega(X\hat\rho^k(Y))A$$ in
the $\sigma$--topology as $n\rightarrow\infty$. \endproclaim
\demo{Proof} If $X\in (\rho^r,\rho^r)$ then
$\hat\varphi^r(X)=\omega(X)I$. Given  $X\in \Cal M^{(0)}_\varphi$
and $\varepsilon>0$, we may pick $r\in \Bbb N$ and  $Y\in
(\rho^r,\rho^r)$ such that $\|Y\|\leq\|X\|$ and
$\|(Y-X)\xi_\omega\|<\varepsilon$. Hence
$\|\hat\varphi^n(Y-X)\xi_\omega\|=\|(\omega(Y)I-\hat\varphi^n(X))\xi_\omega\|<\varepsilon$,
for $n\geq r$, and  $|\omega(Y-X)|<\varepsilon$. Now
$(\iota,\rho^k)$ is finite dimensional and is mapped into itself by
$\hat\varphi$. Any eigenvalue has modulus $\leq 1$ and if we know
that there are no eigenvalues on the unit circle for $k\not=0$ then
$\hat\varphi^n(X)\rightarrow 0=\omega(X)I$ as $n\rightarrow\infty$
for $X\in (\iota,\rho^k)$, hence for $X\in \cup_r 
(\rho^r,\rho^{r+k})$ and, arguing as above, even for $X\in \Cal
M^{(k)}_\varphi$. However, $\hat\varphi(X)=\lambda X$ with
$|\lambda|=1$ implies $\hat\varphi(XX^*)\geq
\hat\varphi(X)\hat\varphi(X)^*=XX^*$. Thus, evalueting in the state
$\omega$, we see that $XX^*=\omega(XX^*)I$ and this implies that
$X=0$ since $(\iota,\rho^k)$ would otherwise contain a unitary
operator which is impossible as we are assuming $d(\rho)>1$. Thus 
$\hat\varphi^n(X)\rightarrow \omega(X)I$ in the s--topology, as
required. To show that $\hat\rho^n(X)\rightarrow\omega (X)I$ in the
$\sigma$--topology, it suffices to show that
$\hat\rho^n(X)\xi_\omega$ tends  weakly to $\omega(X)\xi_\omega$.
However if $Y\in\Cal
M_\varphi,\,\,\omega(Y\hat\rho^n(X))=\omega(\hat\varphi^n(Y)X)
\rightarrow\omega(Y)\omega(X)$ as required. $\omega$ is now the
unique normal $\hat\rho$--invariant state. Thus the associated
modular group commutes with $\hat\rho$ and leaves $\hat\rho(\Cal
M_\varphi)$ globally stable. Hence, by a result of Takesaki
[\ref(Take2)], $\hat\rho\hat \varphi$ is the unique normal
conditional expectation leaving $\omega$ invariant or, equivalently,
since any left inverse of $\hat\rho$ is normal and any state it
leaves invariant is $\hat\rho$--invariant, $\hat\varphi$ is the
unique left inverse with a normal invariant state.

To complete the proof, it suffices, since $\omega$ is faithful to
show that for each  $\ell\in \Bbb Z,\,\,r\in\Bbb N$ and
$S\in(\rho^r,\rho^{r+\ell})$,
$$\omega(S\hat\rho^n(X)A\hat\rho^n(Y))\rightarrow\omega(SA)\omega(X\hat\rho^k(Y))$$ 
But for $n\geq r$ and $\ell\geq 0$, we have 
$$\omega(S\hat\rho^n(X)A\hat\rho^n(Y))=\omega(\hat\rho^{n+\ell}(X)SA\hat\rho^n(Y))
=\omega(\hat\rho^\ell(X)\hat\varphi^n(SA)Y)$$ whose limit as
$n\rightarrow\infty$ is $\omega(\hat\rho^\ell(X)Y)\omega(SA)$. If 
$\ell\leq 0$, a similar computation gives a limit
$\omega(X\hat\rho^{-\ell}(Y))\omega(SA)$. Noting that $\omega(SA)=0$
unless $\ell=-k$, we have the required result.
\par\hfill$\square$\medskip
 
 The last result shows that our systems $(\Cal
M_\varphi,\,\hat\rho)$ have some 
 asymptotic Abelianness and a cluster property. These can be used,
as usual, to derive 
 further interesting properties.
 \proclaim{Lemma 5.5} Let $r\in \Bbb N_0$ and $E,\,E'$ be two
projections in 
 $\Cal M^{(0)}_\varphi$. If there is an $X\in \Cal M_\varphi$ with
$E\hat\rho^r(X)E'\neq 0,$ then there is a $Y\in \Cal
M^{(0)}_\varphi$ such that $E\hat\rho^r(Y)E'\neq 0$. \endproclaim
\demo{Proof} We may suppose that $X$ has grade $k$. Then letting $F$
be the left support of $\hat\rho^r(X)E'$, $F\in\Cal M^{(0)}_\varphi$
and $EF\neq 0$. Now $E\hat\rho^{r+n}(X^*)F\neq 0$  for some $n$ since
$$E\hat\rho^{r+n}(X^*)F\hat\rho^{r+n}(X)\rightarrow
EF\omega(X^*X)\neq 0.$$ But then
$E\hat\rho^r(\hat\rho^n(X^*)X)E'\neq 0,$  completing the
proof.\hfill$\square$\medskip \proclaim{Corollary 5.6} Let $\tau$
denote the restriction of $\hat\rho$ to $\Cal M^{(0)}_\varphi$ then
$$(\tau^r,\tau^r)=(\hat\rho^r,\hat\rho^r)\cap\Cal
M^{(0)}_\varphi,\quad r\in\Bbb N_0.$$ \endproclaim \demo{Proof} Let
$E$ be a projection in $(\tau^r,\tau^r)$ then
$E\hat\rho^r(Y)(I-E)=0$ for all $Y\in \Cal M^{(0)}_\varphi$. Hence
by Lemma 5.5, $E\hat\rho^r(X)(I-E)=0$ for all $X\in\Cal M_\varphi$,
i.e. $E\in (\hat\rho^r,\hat\rho^r)$.\hfill$\square$\medskip

Note that for $r=0$ this corollary relates the centres of the two
von Neumann algebras. We shall see later, that, with slightly
strengthened hypotheses, cf. Corollary 5.12, 
$(\hat\rho^r,\hat\rho^r)$ is in the zero grade part so that the
spaces of intertwiners  will in fact coincide. We also give at this
stage a result of a similar nature which we will need  later on.
\proclaim{Proposition 5.7} Let $\Cal M^{(0)}_\omega$ denote the
fixed-point algebra of $\Cal M^{(0)}_\varphi$ under the modular
group of $\omega$ and $\mu$ the restriction of  $\hat\rho$ to $\Cal
M^{(0)}_\omega$. then 
$$(\mu^r,\mu^r)=(\hat\rho^r,\hat\rho^r)\cap\Cal M^{(0)}_\omega,\quad
r\in\Bbb N_0.$$ \endproclaim \demo{Proof}We must just prove the
analogue of Lemma 5.5. The only point to note is that  for the
analogous proof to go through instead of taking $X$ to be of grade
$k$ we must take  it to transform like a character under the modular
group. This can be done since the finite--dimensional spaces
$(\rho^r,\rho^s)$ are invariant under the modular group and their
union is a total set in $\Cal M_\varphi$.\hfill$\square$\medskip

If we consider a pair $(\Cal M,\rho)$ consisting of a factor $\Cal
M$ and a normal endomorphism $\rho$ and if $\Cal M$ admits a unique
normal $\rho$--invariant state  and this state is faithful then we
may define the {\it Jones index} of $\rho$ to be $d(\varphi)^2$, 
where $\varphi$ is the unique left inverse of $\rho$ leaving the
state invariant. We have just seen in Lemma 5.4 how this situation
arises in our context. When the state in question  is a trace, this
definition coincides with the usual definition of the Jones index
for the inclusion $\rho(\Cal M)\subset\Cal M$. This index continues
to be multiplicative in the sense that if endomorphisms
$\rho,\,\sigma,$ and $\rho\sigma$ of $\Cal M$ each satisfy the above
conditions then the Jones index of $\rho\sigma$ is a product of the
Jones indices of $\rho$  and $\sigma$.

To illustrate the above construction and motivate what follows, we
consider a class of examples. We first take our object to be a
finite-dimensional Hilbert space $H$ of dimension $d >1$ then $\Cal
O_H$ is just the  Cuntz algebra. The von Neumann algebra $\Cal
M^{(0)}_H$ has a faithful trace state and is an inductive limit of
full matrix algebras. It is therefore the hyperfinite factor of type
$II_1$.  Since $(\Bbb C,H)\not=0, \Cal M_H$ is the hyperfinite
factor of type $III_{1\over d}.$

More generally, we can consider a faithful left inverse $\Phi$ of
$H$ and construct the von  Neumann algebra $\Cal M_\Phi$. This is a
factor by Lemma 5.4 since any element commuting  with $H$ is
necessarily invariant under $\rho$. However, we have the following
stronger result.

\proclaim{Lemma 5.8} Let $\theta\in (H^2,H^2)$ denote the unitary
flip and $\sigma :=\sigma_H$ the endomorphism of $\Cal M_\Phi$
generated by $H$ then an $X\in \Cal M_\Phi$ such that
$[X,\sigma^r(\theta)]=0,\,\, r\in\Bbb N_0$ is a multiple of the
unit. \endproclaim \demo{Proof} Since $\{X\in\Cal M_\Phi :
[X,\sigma^r(\theta)]=0,\,\, r\in \Bbb N_0\}$ is  globally invariant
under gauge transformations and weak operator closed, it suffices
to  prove the result when $\alpha_\lambda(X)=\lambda^k\,X$. First
suppose $k=0$ and pick  $Y,Y'\in (H^r,H^r)$ then 
$$\omega(Y^*XY')=\omega(\theta(s,1)Y^*XY'\theta(s,1)^*)=
\omega(\sigma(Y^*)\theta(s,1)
X\theta(s,1)^*\sigma(Y'))=\omega(\sigma(Y^*)X\sigma(Y')),$$ provided
$s\geq r$. Here $\theta(s,1)$ permutes $H^s$ and $H$. Hence 
$$(Y\xi_\omega,XY'\xi_\omega)=(Y\xi_\omega,V^*XVY'\xi_\omega)=(Y\xi_\omega,
\hat\Phi(X)Y'\xi_\omega),$$ for $Y,Y'\in\cup_r(H^r,H^r)$. Thus
$X=\hat\Phi(X)$, so that $X\in \Bbb CI$ by Lemma 5.4. Now suppose
$k\not= 0$, then $X^*X$ and $XX^*$ have grade zero and  are hence
equal to $\|X\|^2I$. Now for every element $Y$ of grade $k$, 
$$\theta(r+k,1)Y\theta(r,1)^*\rightarrow\hat\rho(Y)$$ in the
s-topology as $r\rightarrow\infty,$ since this is trivially true on
the s--dense set $\Cal O^{(k)}_\rho$ and $\omega$ is invariant under
the inner endomorphisms generated by the $\theta(s,1)$.  Taking
$Y=X$ and multiplying on the left by $X^*$, we deduce that
$$\|X\|^2\,\theta(r+k,1)\theta(r,1)^*\rightarrow X^*\hat\rho(X),$$
as $r\rightarrow\infty$. Applying $\hat\Phi$, we get
$$\|X\|^2\,\hat\Phi(\theta(r+k,1)\theta(r,1)^*)\rightarrow\hat\Phi(X^*)X,$$
as $r\rightarrow\infty$. But $\hat\Phi(X)$ is again of grade $k$ and
commutes with each  $\sigma^r(\theta)$. It is hence just a multiple
of $X$: $\hat\Phi(X)=\mu X$. Thus
$$\|X\|^2\,\hat\Phi(\theta(r+k,1)\theta(r,1)^*)\rightarrow\bar\mu\,\|X\|^2.$$ 
If $X\not=0$, we can compute $\mu$; it is just
$$\omega(\theta(r+k,1)\theta(r,1)^*)=\sum_j\lambda^{k+1}_j\not=0,$$
where the sum is taken over the eigenvalues of the density matrix
determining the left inverse $\Phi$, cf. Lemma 5.2. Writing
$\theta(r;k):=\theta(r+k,1)\theta(r,1)^*$ for short,
$\Phi(\theta(r;k))\rightarrow\mu I$  and this would imply that
$\omega(\Phi(\theta(r;k)^*\Phi(\theta(r;k)))$ tends to $\mu^2I$.
However, choosing an orthonormal basis $\psi_i$ of $H$ consisting of
eigenvectors of the density matrix determining $\Phi$, we have
$\hat\Phi(Y)=\sum_i\lambda_i\psi_i^*Y\psi_i.$ Thus 
$$\hat\Phi(\theta(r;k))=\sum_i\lambda_i\psi_i^*\hat\rho^r(\psi_i).$$
A further computation shows that
$$\omega(\hat\Phi(\theta(r;k)^*)\hat\Phi(\theta(r;k)))=\sum_i\lambda_i^3\not=\mu^2$$
and  this contradiction shows that $X=0$ and completes the
proof.\hfill$\square$\medskip

We now take our object $\rho$ to be a continuous unitary
representation of a compact group  $G$ on $H$. Now $\Cal O_\rho$ can
be identified with the fixed points of $\Cal O_H$ under  the action
of $G$ induced by $\rho$ [\ref(DR2)]. Since the standard left
inverse of $\rho$ is just the restriction of the standard left
inverse of $H, \Cal M_\rho$ is just the von Neumann  subalgebra of
$\Cal M_H$ generated by $\Cal O_\rho$. As the state $\omega$ on
$\Cal M_H$ is invariant under the action of $G, \Cal M_\rho$ is the
fixed-point algebra of $\Cal M_H$ under  the action of $G$ induced
by $\rho$. More generally, any faithful left inverse $\varphi$ of
$\rho$ extends to a G--invariant faithful left inverse $\Phi$ of
$H$. Hence, arguing as above,  $\Cal M_\varphi$ is the fixed--point
algebra of $\Cal M_\Phi$ under the action of $G$ induced by $\rho$.

It follows from Lemma 5.8, that $\Cal M_\varphi{}'\cap \Cal
M_\Phi=\Bbb CI$  and that $\Cal M^{(0)}_\varphi$ and hence $\Cal
M_\varphi$ are factors but  this lemma also yields a stronger result.
\proclaim{Theorem 5.9} Let $\rho$ be a continuous unitary
representation of a compact group $G$ of finite dimension $d > 1$
then $$(\rho^r,\rho^s)=(\hat\rho^r,\hat\rho^s),\quad r,s\in\Bbb
N_0$$ on $\Cal M_\varphi$.\endproclaim \demo{Proof} Let $X\in
(\hat\rho^r,\hat\rho^s)$ and $\psi\in H^r,\,\,\psi'\in H^s$ then 
$\psi'{}^*X\psi\in\Cal M_\varphi{}'\cap\Cal M_\Phi$ so that
$\psi'{}^*X\psi\in\Bbb CI.$ Thus  $X\in (H^r,H^s)\cap\Cal
M_\varphi=(\rho^r,\rho^s)$ as required.\hfill$\square$\medskip

On the basis of the above result, it is clear that two systems of
the form $(\Cal M_\Phi,\hat\rho)$ are isomorphic if and only if the
density matrices defining the left  inverses $\Phi$ are unitarily
equivalent.

Now S. Popa [\ref(Popa1)] has introduced various notions of 
amenability for inclusions of factors with finite index, at the same
time providing a classification. In particular the graph
$\Gamma_{\Cal N,\Cal M}$ of an inclusion $\Cal N\subset\Cal M$ of
factors is said to be amenable if $\|\Gamma_{\Cal N,\Cal M}\|^2$ 
equals the minimal index of the inclusion. The analogue condition
for an object $\rho$ of  a  tensor $C^*$--category is the equality
$d(\rho)=\|m^\rho\|$ discussed at the end of \S 3. We prefer,
however, to take the conclusions of Theorem 5.9  as the definition
of an object $\rho$ being {\it amenable}. The analogous property for
$\hat\rho$  considered as an endomorphism of the $C^*$--algebra and
proved as Theorem 3.5 of  [\ref(DR2)] will be termed {\it
$C^*$--amenable}. A special object in a symmetric  tensor
$C^*$--category is $C^*$--amenable, Lemma 4.14 of [\ref(DR1)]. In
the same  way if $\varphi$ is a normalized left inverse for $\rho$,
we say that $\varphi$ is  amenable  if the corresponding property
holds when $\hat\rho$ denotes the induced  endomorphism  of $\Cal
M_\varphi$. This definition is closely related to one of the
consequences of what Popa terms strong amenability, namely that the
relative commutants in a Jones tunnel coincide with the relative 
commutants in the associated standard model $\Cal N^{st}\subset\Cal
M^{st}$, i.e. $$\Cal N_{i}^{st'}\cap\Cal N_{j}^{st}=\Cal
N_{i}^{'}\cap\Cal N_j.$$

Associated with an object $\rho$ in a tensor $C^*$--category $\Cal
T$ we have a
 tensor $C^*$--category, $\Cal T_\rho$, whose objects are the
integers $\Bbb N_0$ and where $(r,s)=(\rho^r,\rho^s)\quad r,s\in\Bbb
N_0,$ with the structure inherited from $\Cal T$. In practice the
objects $\rho^r$ will  usually be distinct, in particular when
$d(\rho) > 1$, and $\Cal T_\rho$ can be identified with the full
subcategory of $\Cal T$ whose objects are the powers of $\rho$. When
$\rho$ is an object of $\Cal T_f$, the  construction of $\Cal
M_\varphi$ provides a natural $^*$--functor $F_\varphi$ of $\Cal
T_\rho$  into the full tensor $C^*$--subcategory of $\text{End}(\Cal
M_\varphi)$ whose objects are the  
 powers of $\hat\rho$. $F_\varphi$ maps  $\rho^r$ to the endomorphism
$\hat\rho^r$ and an arrow $T$ into its embedding in $\Cal M_\varphi$
again denoted by $T$. Then $\varphi$  amenable just means that the
functor $F_\varphi$ is full, i.e. that $$(\hat\rho^r, \hat\rho^s) =
{(\rho^r,\rho^s)}$$ for all integers $r,s\geq 0$.

Note that we are assuming $(\iota,\iota)=\Bbb C$ hence, if $\varphi$
is  amenable, then $\Cal M_\varphi$ is a factor since $$Z(\Cal
M_\varphi) = (\hat\iota, \hat\iota) = (\iota,\iota) = \Bbb C.$$

In  our discussion of $\Cal O_\rho$ and $\Cal M_\varphi$ we have, up
till now, made no direct use of a conjugate $\bar\rho$ for $\rho$;
indirectly, of course, it has been used to provide  a standard left
inverse. The problem, of course, is that $\bar\rho$ cannot, in
general, be regarded as an endomorphism of $\Cal O_\rho$ or $\Cal
M_\varphi$ nor can solutions $R,\,\bar R$ of the conjugate
equations, in general, be interpreted as elements of $\Cal O_\rho$
or  $\Cal M_\varphi$. However, the following results  provide a
certain substitute.  \proclaim{Lemma 5.10} Let $\rho$, an object of
$\Cal T_f$, have a conjugate $\bar\rho$  that is a direct sum of
subobjects of tensor powers of $\rho$ and let $F :\Cal
T_\rho\to\hat{\Cal T}$ be a tensor $^*$--functor into a strict
tensor $C^*$--category  $\hat{\Cal T}$ with subobjects and (finite)
direct sums. Then $\hat\rho :=F(\rho)$ has a  conjugate in
$\hat{\Cal T}$ and for every left inverse $\varphi$ of $\rho$
associated with solutions of the conjugate equations there is a left
inverse $\hat\varphi$ of $\hat\rho$  such that, with an obvious
abuse of notation, $F\varphi=\hat\varphi F.$ If the units of $\Cal
T$ and $\hat{\Cal T}$ are irreducible then $d(\hat\rho)\leq
d(\rho)$. If, furthermore,  the map from $(\rho,\rho)$ to
$(\hat\rho,\hat\rho)$ induced by $F$ is surjective then
$d(\hat\rho)=d(\rho).$ \endproclaim \demo{Proof} $\bar\rho$ is a
direct sum of subobjects of tensor powers of $\rho$. Take a direct 
sum of objects representing the images of the corresponding
projections under $F$.  This can be done since $\hat{\Cal T}$ has
subobjects and direct sums. The resulting object should be a
conjugate for $F(\rho)$. This can be checked by taking solutions of
the conjugate equations for $\rho$, using them to construct
solutions of the conjugate equations for  $\hat\rho$ and performing
the necessary computations. Those who wish to avoid explicit
computations can take the following longer route. The obvious full
and faithful canonical  tensor $^*$--functor from $\Cal T_\rho$ to
$\Cal T$ extends to a full and faithful tensor $^*$--functor between
the canonical completions $C(\Cal T_\rho)$ and $C(\Cal T)$ under
subobjects and direct sums (see Appendix). Since $C(\Cal T_\rho)$
contains an object whose image is a conjugate for $\rho$, this
object must be a conjugate for $\rho$ considered as an object of
$C(\Cal T_\rho)$. Now $F$ extends to a tensor $^*$--functor $C(F) :
C(\Cal T_\rho)\to C(\hat{\Cal T})$ so $C(F)(\rho)$ has a conjugate
which can be taken to be in $\hat{\Cal T}$ as every object of
$C(\hat{\Cal T})$ is equivalent to an object of $\hat{\Cal T}$.  The
remaining assertions are now obvious.\hfill$\square$\bigskip

\proclaim{Proposition 5.11} Let $\rho$ have a conjugate $\bar\rho$
which is a direct sum of  subobjects of tensor powers of $\rho$ then
$\hat\rho$ has a conjugate $\bar{\hat\rho}$ in  End$\Cal M_\varphi$
and the index of the inclusion $\hat\rho(\Cal M_\varphi)\subset \Cal
M_\varphi$  with respect to the normal conditional expectation $\Cal
E :=\hat\rho\hat\varphi$ is  $$[\Cal M_\varphi : \hat\rho(\Cal
M_\varphi)]_{\Cal E}=d(\varphi)^2.$$  In particular, if $\Cal
M_\varphi$ is a factor, $d(\hat\rho)\leq d(\rho)$ and if $\varphi$
is amenable, $d(\hat\rho)=d(\rho).$  \endproclaim   \demo{Proof} We
apply Lemma 5.10 to the functor $F_\varphi$ from $\Cal T_\rho$ to
End$\Cal M_\varphi$ noting that the image category not only has
subobjects but also direct  sums since $\Cal M_{\varphi}$ is
properly infinite. It only remains to remark that the formula for
the index follows from [\ref(Long2)], since we have a solution of
the conjugate equations for $\hat\rho$  with associated left inverse
$d(\varphi)\hat\varphi$. \hfill$\square$\bigskip

We now want to give conditions implying that $\varphi$ is amenable.
One possibility is to follow the  path taken in [\ref(DR1)] when
proving (Lemma 4.14) that special objects in a  
 symmetric tensor $C^*$--category are $C^*$--amenable. We  will
first generalize the most important step, the Grading Lemma, Lemma
4.5 of 
 [\ref(DR1)].

\proclaim{Lemma 5.12} Let $\rho$ be a unital endomorphism of a
$C^*$--algebra $\Cal A$ and suppose that for each $m\in \Bbb N$
there is an isometry $S_m\in (\iota, \rho^{n(m)})$ such that
$$S_m^*\rho^m(S_m)=\lambda_m,$$ where $\lambda_m=\lambda_m^*$,
$\rho(\lambda_m)=\lambda_m$, $\|\lambda_m\|<1$  and $$\prod_m
{1+\|\lambda_m\| \over 1-\|\lambda_m\|}<+\infty.$$ Then $[\rho]^k,$
the norm closure of $\dot{[\rho]}^k=\cup_{r,r+k\geq
0}(\rho^r,\rho^{r+k})$,  $k\in \Bbb Z,$ determine a grading of the
$C^*$--algebra that they generate and the canonical map of $\Cal
O_\rho$ into $\Cal A$ is injective. \endproclaim \demo{Proof}  Let
$$X=X_0+(\sum_{k=-n}^{-m} + \sum_{k=m}^n)X_k$$ with
$X_k\in\dot{[\rho]}^k$ then, setting
$$X'=(I-\rho^r(\lambda_m))^{-1}(\rho^r(S_m^*)X\rho^r(S_m)-\rho^r(\lambda_m)X),$$
we have $$X'=X_0+(\sum_{k=-n}^{-m-1} + \sum_{k=m+1}^n)X'_k,$$ where 
$$X_k'=(I-\rho^r(\lambda_m))^{-1}(\rho^r(S_m^*)\rho^{r+k}(S_m)-\rho^r(\lambda_m))X_k\in\dot{[\rho]}^k$$
provided $r$ is chosen so that each $X_k\in (\rho^r,\rho^{r+k})$. We
conclude that $X=0$ implies $X_0=0$. But $X=0$ implies  $X^*X=0$ and
$(X^*X)_0=\sum_kX_k^*X_k$ so $X_k=0.$ Thus the linear span of the
$\dot{[\rho]}^k$ is a $\Bbb Z$--graded $^*$--algebra 
$\dot{[\rho]}$. Now $\|X'\|\leq {(1+\|\lambda_m\|)\over
(1-\|\lambda_m\|)} \|X\|.$ Thus $$\|X_0\|\leq \prod_m
{1+\|\lambda_m\| \over 1-\|\lambda_m\|}\|X\|<+\infty.$$ This
computation shows that the projection $m_0$ onto the zero grade part
of $\dot{[\rho]}$  is continuous for any $C^*$--norm on
$\dot{[\rho]}$ (cf. Lemma 2.10 of [\ref(DR2)]).  The argument of
Lemma 2.11 of [\ref(DR2)] now shows that $\dot{[\rho]}$ has a unique
$C^*$--norm extending the $C^*$--norm on $\dot{[\rho]}^0$. But
$\dot{[\rho]}^0$ being  an inductive limit of $C^*$--algebras admits
a unique $C^*$--norm. Hence $\dot{[\rho]}$  admits a unique
$C^*$--norm which must be the $C^*$--norm defined by both $\Cal
O_\rho$  and the ambient $C^*$-algebra $\Cal A$. Thus the canonical
map of $\Cal O_\rho$ into  $\Cal A$ is injective and the
$\dot{[\rho]}^k$ determine the grading on the image.
\hfill$\square$\bigskip We immediately have the following corollary.
\proclaim{Corollary 5.13} Let $\rho$ be an object in a tensor
$C^*$--category and suppose  that for each $m\in \Bbb N$ there exist
an isometry $S_m\in (\iota,\rho^{n(m)})$ such that  $$S_m^*\otimes
1_{\rho^m}\circ 1_{\rho^m}\otimes S_m :=\lambda_m,$$  satisfies the
conditions in Lemma 5.11. If $\hat\rho$ denotes the canonical 
endomorphism on $\Cal O_\rho$ then $[\hat\rho]^k=\Cal O^{(k)}_\rho.$
In the same  way if $\hat\rho$ denotes the canonical endomorphism on
$\Cal M_\varphi$, where  $\varphi$ is a normalized left inverse for
$\rho$, then $[\hat\rho]^k$ is $\sigma$--dense  in $\Cal
M^{(k)}_\varphi$. \endproclaim

In all applications of the above corollary to date, $\lambda_m=\pm
d^{-m}\,1_{\rho^m}.$ Apart from the case of special objects in a
symmetric tensor $C^*$--category (cf. Lemma 3.7  of [\ref(DR1)]),
the result applies to {\it real} objects of dimension greater than
one  in an arbitrary tensor $C^*$--category. These are objects
$\rho$ for which there is a solution $R\in (\iota,\rho^2), \,\bar
R=R$ of  the conjugate equations. In particular $\rho$ is
self-conjugate. To see that the corollary  can be applied, it
suffices to note that any power $\rho^m$ of a real object is real. 
The corresponding operator $R$ can be computed inductively using the
formula for products  given in the proof of Theorem 2.4:
$$R_m=1_\rho\otimes R_{m-1}\otimes 1_\rho\circ R;\quad R_1=R.$$ 
Thus $R_m^*\circ R_m=d^m$ where $d := R^*\circ R$. We can now take
$S_m=\sqrt{d^{-m}}R_m$  in Corollary 5.13.

It is easily checked that the direct sum of real objects is real and
that an object of the form  $\sigma\oplus\bar\sigma$ or
$\bar\sigma\sigma$ is real. If $\rho$ is self-conjugate and
irreducible then we can either choose $\bar R=R$ or $\bar R=-R$. The
second case we call  {\it pseudoreal} adopting the terminology used
for group representations. We may now give  a characterization of
real objects in terms of irreducibles. \proclaim{Theorem 5.14} An
object $\rho$ is real if and only if $(\sigma,\rho)$ and 
$(\bar\sigma,\rho)$ have the same dimension for any irreducible
$\sigma$ and  $(\sigma,\rho)$ has even dimension for any pseudoreal
irreducible $\sigma$. \endproclaim \demo{Proof} If $\rho$ satisfies
the conditions of the theorem, it is real by the above remarks.  On
the other hand, if $\rho$ is real, it is self-conjugate so the first
condition must be  satisfied. Now suppose $\sigma$ is a pseudoreal
irreducible, then consider the  antilinear map  $X\mapsto
X^{\bullet}$ on $(\sigma,\rho)$ where we agree to take $\bar R_\rho=
R_\rho$  and $\bar R_\sigma=-R_\sigma$. In this case, it is easily
checked that $X^{\bullet\bullet}=-X.$  This implies that
$(\sigma,\rho)$ has even dimensions.\hfill$\square$\bigskip This
characterization of real objects and the above formulae make it
clear that we  may restrict ouselves to standard solutions of the
conjugate equations when defining real  objects.

If we can take $\bar R=-R$ as in the case of a pseudoreal
irreducible then again Corollary 5.12  is applicable; the formulae
differ from the real case only by the presence of a minus sign  when
$m$ is odd. It is thus natural to ask whether Corollary 5.12 may be
applied to any  self-conjugate object of dimension greater than one.
We do not know the answer.

The next lemma provides a simpler and more general way of proving
$C^*$--amenability than that taken in Lemma 4.14 of [\ref(DR1)].
\proclaim{Lemma 5.15} Suppose that $[\hat\rho]^k=\Cal O^k_\rho$ and
that,  for some $n\in\Bbb N$, $(\iota,\rho^n)\not=0$ then $\rho$ is
$C^*$--amenable. \endproclaim \demo{Proof} If $S\in (\iota,\rho^n)$
is an isometry, so is $S^{\otimes m}\in (\iota,\rho^{mn}),
\,\,m\in\Bbb N.$ Suppose $X\in (\hat\rho^r,\hat\rho^s)$ then
$X=\hat\rho^s(V^*)X\hat\rho^r(V)$ for any isometry $V$. Since
$[\hat\rho]^k=\Cal O^k_\rho$, we  can find $X_m\in
(\rho^{r+k(m)n},\rho^{s+k(m)n})$ converging in norm to $X$. Picking
an  isometry $S_m\in (\iota,\rho^{k(m)n})$ and setting
$Y_m:=\hat\rho^s(S_{m}^*)X_m\hat\rho^r (S_m)\in (\rho^r,\rho^s)$ we
have $$\|Y_m - X\|\leq\|\hat\rho^s(S_{m}^*)(X_m -
X)\hat\rho^r(S_m)\|\leq\|X_m - X\|.$$ But $(\rho^r,\rho^s)$ is
closed so $X\in (\rho^r,\rho^s)$ and $\rho$ is $C^*$--amenable.
\hfill$\square$\medskip

We now give a criterion for $\varphi$ to be amenable in the form of
a theorem. \proclaim{Theorem 5.16} Let $\rho$ be real or irreducible
and pseudoreal then $\rho$ is  $C^*$--amenable.  Let $R\in
(\iota,\rho^2)$  be such that 
$$R^*XR=d(\varphi)\,\hat\varphi(X),\quad X\in (\rho,\rho),$$  then
$\varphi$ is amenable if and only if 
$$R^*XR=d(\varphi)\,\hat\varphi(X),\quad X\in
(\hat\rho,\hat\rho).$$  Hence, in particular, in this case
$d(\rho)=d(\hat\rho)$ and $\Cal M_\varphi$ is a factor. \endproclaim
\demo{Proof} The $C^*$--amenability of $\rho$ follows from Corollary
5.13 and  Lemma 5.15. Now  define $S_m$ as following Corollary 5.13.
This corollary shows that $(\hat\rho^r,\hat\rho^s)\subset \Cal
M_\rho^{(s-r)}$.  Hence, given $X\in (\hat\rho^r, \hat\rho^s)$ there
exists $X_m\in(\rho^{r+k(m)},\rho^{s+k(m)})$  uniformly bounded such
that $X_m\xi_\omega$ tends to $X\xi_\omega.$ Set
$Y_m=\hat\rho^s(S_{k(m)}^*)X_m\hat\rho^r(S_{k(m)})\in
(\rho^r,\rho^s)$ and we have  $$\leqalignno{\|(Y_m -
X)\xi_\omega\|^2 &=\|\hat\rho^s(S_{k(m)}^*)(X_m -
X)\hat\rho^r(S_{k(m)})\xi_\omega\|^2\cr &\leq\|(X_m -
X)\hat\rho^r(S_{k(m)})\xi_\omega\|^2 
=\omega(\hat\rho^r(S_{k(m)}^*)(X_m - X)^*(X_m -
X)\hat\rho^r(S_{k(m)})).\cr}$$  Now $(X_m - X)^*(X_m - X)\in
(\hat\rho^{r+k(m)},\hat\rho^{r+k(m)})$ and for any such  element $B$
we have $$\hat\rho^r(S_{k(m)}^*)B\hat\rho^r(S_{k(m)})
=\hat\psi^{k(m)}(B),$$  where $\hat\psi$ is the right inverse of 
$\hat\rho$ such that
$\hat\psi_{\iota,\iota}=\hat\varphi_{\iota,\iota}$. Although we do
not, a priori, know that $\Cal M_\varphi$ is a factor, we may still
conclude that the same relation holds between the powers of
$\hat\varphi$ and $\hat\psi$, since the computation in Theorem 3.8
is valid without assuming the tensor unit to be irreducible. Since,
$\omega\circ\hat\varphi=\omega$, we get
$$\omega(\hat\psi^{k(m)}(B))=\omega(\hat\varphi^r\hat\psi^{k(m)}(B))=\omega(\hat
\psi^{r+k(m)}(B))=\omega(\hat\varphi^{r+k(m)}(B))=\omega(B).$$  We
conclude that $\|(Y_m - X)\xi_\omega\|\leq\|(X_m - X)\xi_\omega\|$
so that $Y_m\xi_\omega$  converges to $X\xi_\omega$. Since
$\xi_\omega$ is separating and  $Y_m$ is uniformly bounded, $Y_m$
converges to $X$ in the strong operator topology. But
$(\rho^r,\rho^s)$ is closed in this topology, so $X\in
(\rho^r,\rho^s)$ and $\varphi$ is amenable. The converse is trivial.
\hfill$\square$\medskip

The following result is an immediate corollary. \proclaim{Corollary
5.17} Let $\rho$ be real or irreducible and pseudoreal then $\varphi$
is amenable if and only if  $(\rho,\rho)=(\hat\rho,\hat\rho)$.
\endproclaim 

In a way, the natural condition for $\varphi$ to be amenable would
be simply that $d(\rho)=d(\hat\rho)$. However, we do not know a
priori that $\Cal M_\varphi$ is a factor so that $d(\hat\rho)$ might
not be defined. Therefore, in the next result, we add a condition
which guarantees that $\Cal M_\varphi$ is a factor and which is
anyway automatically fulfilled whenever $\Cal M_\varphi$ is a
factor. \proclaim{Theorem 5.18} Let $\rho$ be real or irreducible
and pseudoreal and suppose  there is an $R\in
(\hat\iota,\hat\rho^2)$ such that
$$R^*XR=d(\varphi)\,\hat\varphi(X),\quad X\in(\hat\rho,\hat\rho),$$
then $\Cal M_\varphi$ is a factor and if $d(\rho)=d(\hat\rho)$ then
$\varphi$ is amenable. \endproclaim \demo{Proof} Let $Z$ be in the
centre of $\Cal M_\varphi$ then $R^*ZR=ZR^*R=Zd(\varphi)$.  Hence
$Z=\hat\varphi(Z)$ so that $\Cal M_\varphi$ is a factor by Lemma
5.4. Since  $d(\rho)=d(\hat\rho)$, a standard $R'\in (\iota,\rho^2)$
for $\rho$ will also be  standard for $\hat\rho$. Hence
$$R'{}^*XR'=R'{}^*\hat\rho(X)R',\quad X\in (\hat\rho,\hat\rho).$$
Now there is an invertible $Y\in (\rho,\rho)$ such that
$1_\rho\otimes Y\circ R'$  induces $d(\varphi)\hat\varphi$. Hence 
$$R'{}^*\hat\rho(Y^*XY)R'=d(\varphi)\,\hat\varphi(X),\quad X\in \Cal
M_\varphi,$$  so for $X\in (\hat\rho,\hat\rho)$ this is just
$R'{}^*Y^*XYR'$ and taking $Y\otimes 1_\rho \circ R'$ for $R$ in
Theorem 5.16, the result follows.\hfill$\square$\medskip

\proclaim{Corollary 5.19} Let $\rho$ be real or irreducible and
pseudoreal then if $\varphi$ is  not amenable, either $\Cal
M_\varphi$ is not a factor, or if it is a factor
$d(\rho)>d(\hat\rho)$. \endproclaim\medskip

\proclaim{Corollary 5.20} If $\rho$ is real or irreducible and
pseudoreal and $d(\rho)= \|m^\rho\|$ then $\varphi$ is amenable  if
and only if $\Cal M_\varphi$ is a factor. \endproclaim \demo{Proof}
$d(\rho)=\|m^\rho\|\leq\|m^{\hat\rho}\|\leq d(\hat\rho)\leq
d(\rho)$  so that we have equality
everywhere.\hfill$\square$\medskip

The next results provide us with a class of examples of amenable
objects. \proclaim{Proposition 5.21} Let $\rho$ be real or
irreducible and pseudoreal and let  $\Cal M_\varphi$ be a factor
then  $\hat\varphi$ is amenable. \endproclaim \demo{Proof} It is
obvious from Lemma 5.12 that $\Cal M_\varphi=\Cal M_{\hat\varphi}$
and $\hat{\hat\rho}=\hat\rho$.\hfill$\square$\medskip For our next
example, we would like to make contact with the work of Popa. To
facilitate this we give a result relating an endomorphism $\hat\rho$
and its restriction to the grade zero part. \proclaim{Proposition
5.22}Let $\rho$ be real or irreducible and pseudoreal and let $\tau$
denote the restriction of $\hat\rho$ to $\Cal O^{(0)}_\rho$ then
$$(\hat\rho^r,\hat\rho^r)=(\tau^r,\tau^r),\quad r\in \Bbb N_0,$$
hence $\|m^{\hat\rho}\|=\|m^\tau\|.$ Furthermore, if $\Cal O_\rho$
has trivial centre, $d(\hat\rho)=d(\tau)$. The analogous statement
holds if $\Cal O_\rho$ is replaced by $\Cal M_\varphi$. In the case
of $\Cal M_\varphi$, the Jones indices of $\hat\rho$ and $\tau$
coincide.  \endproclaim  \demo{Proof} To prove the equality
$(\hat\rho^r,\hat\rho^r)=(\tau^r,\tau^r)$, it suffices to consider
the case of $\Cal M_\varphi$ since this implies the result for $\Cal
O_\rho$. However, this case is a consequence of Corollary 5.6 and
Corollary 5.12. Now by the characterisation of left inverses in
Lemma 2.5, every left inverse $\chi$ of $\hat\rho$ preserves the
grading and defines a left inverse of $\tau$ by restriction. Since
$(\hat\rho,\hat\rho)=(\tau,\tau)$,  every left inverse of $\tau$
arises in this way. In particular, the standard left inverse arises
in this way. However, the dimensions of the faithful left inverses
cannot change on restriction  since the corresponding Jones
projections have grade zero. Hence $d(\hat\rho)=d(\tau)$. For the
same reason, the Jones indices of $\hat\rho$ and $\tau$ coincide. 
\hfill$\square$\medskip  \proclaim{Proposition 5.23} Let $\Cal
N\subset\Cal M$ be a proper inclusion  of properly infinite factors 
with  normal conditional expectation $\Cal E$ and $\Gamma_{\Cal
N,\Cal M}$ amenable  and ergodic then  $\varphi$ is amenable, where
$\varphi$ is the left inverse of the canonical endomorphism $\gamma
: \Cal M\to \Cal N\subset\Cal M$ induced by $\Cal E$ (i.e. the
composition of $\Cal E$ and its dual).\endproclaim \demo {Proof}
Saying that $\Gamma_{\Cal N,\Cal M}$ is ergodic in the sense of Popa
just means that  the completion of $\cup_n(\gamma^n,\gamma^n)$ in
the representation of the trace induced by $\Cal E$ is a factor of
type $II_1$ [\ref(Popa3)]. For our purposes, we do not need to
assume that $\Cal E$ induces a trace on the relative commutants in
the tunnel and interpret ergodicity to mean simply that  $\Cal
M^{(0)}_\varphi$ is a factor. But by Proposition 5.22, $\Cal
M_\varphi$ is then a factor. Since amenability of $\Gamma_{\Cal
N,\Cal M}$ in the sense of Popa means that 
$\|m^\gamma\|=\|\Gamma_{\Cal N,\Cal M}\|^2$ is the  minimal index of
the inclusion, i.e. $d(\gamma)=\|m^\gamma\|,$ and $\gamma$ is real,
the result  now follows from Corollary 5.20. \hfill$\square$\medskip
We therefore see that amenability in the sense of Popa is closely
related to our use of the term. Indeed, by using Popa's results we
obtain a partial converse to Corollary 5.20. Our argument involves
the endomorphism $\mu$, the restriction of $\hat\rho$ to $\Cal
M^{(0)}_\omega$. \proclaim{Lemma 5.24}The following conditions are
equivalent.  \item{$a)$}Some Jones projection of $\hat\rho$ lies in
$\Cal M^{(0)}_\omega$, \item{$b)$}Some Jones projection for
$\hat\varphi$ lies in $\Cal M^{(0)}_\omega$. \item{$c)$}Some Jones
projection for the standard left inverse of $\hat\rho$ lies in $\Cal
M^{(0)}_\omega$.

Under these conditions every normalized faithful left inverse of
$\mu$ has a unique extension to a faithful normalized left inverse
of $\hat\rho$ and some Jones projection for this extension lies in
$\Cal M^{(0)}_\omega$. Furthermore 
$\|m^\mu\|\leq\|m^{\hat\rho}\|,\,d(\mu)=d(\hat\rho)$ and the Jones
indices of $\mu$ and $\hat\rho$ coincide.  These equivalent
conditions are implied by the following: there is a solution
$R,\,\bar R\in (\hat\iota,\hat\rho^2)$ of the conjugate equations
such that both $R$ and $\bar R$ define $\hat\varphi$. In this case,
we say that $\varphi$ is selfconjugate.\endproclaim   \demo{Proof}We
first show the equivalence of the three conditions. Suppose
$R,\,\bar R\in (\iota,\hat\rho^2)$ solve the conjugate equations and
are such that $\bar R\bar R^*\in \Cal M^{(0)}_\omega$. Then
$\sigma_t(\bar R)=a^{it}\bar R,\,\, t\in \Bbb R$, for some $a\in
\Bbb R_+$. If $\chi$ is the corresponding faithful normalized left
inverse of $\hat\rho$, then $\chi(\Cal M^{(0)}_\omega)\subset\Cal
M^{(0)}_\omega$ and $\chi$ restricts to a faithful normalized left
inverse of $\mu$. By Lemma 2.5, every faithful normalized left
inverse of $\mu$ is of the form $A\mapsto \chi(X^*AX)$ for some
$X\in (\mu,\mu)$ with $\chi(X^*X)=1$. Now $(\mu,\mu)\subset
(\hat\rho,\hat\rho)$ by Proposition 5.7,  so the same formula
defines a faithful normalized left inverse of $\hat\rho$ and, since
$X^{*-1}\bar R\bar R^* X^{-1}\in \Cal M^{(0)}_\omega$, there is an
associated Jones projection in $\Cal M^{(0)}_\omega$. Since any
faithful normalized left inverse of $\hat\rho$ is uniquely
determined by its values on a Jones projection, every faithful
normalized left inverse of $\mu$ has a unique extension to a
faithful normalized left inverse of $\hat\rho$. Since both
$\hat\varphi$ and the standard left inverse leave $\Cal
M^{(0)}_\omega$ globally stable, $a)$ implies both $b)$ and $c)$ so
the conditions are obviously equivalent. Now any faithful normalized
left inverse of $\mu$ has the same dimension as its extension to a
normalized faithful left inverse of $\hat\rho$ since a Jones
projection for the latter lies in $\Cal M^{(0)}_\omega$. It is now
clear that $\|m^\mu\|\leq\|m^{\hat\rho}\|,\,d(\mu)= d(\hat\rho)$ and
that the Jones indices of $\mu$ and $\hat\rho$ coincide. Now suppose
that $R$ and $\bar R$ both define $\hat\varphi$. Then if $S\in
(\iota,\hat\rho^2)$,
$$\omega(SA)=\omega(\hat\rho^2(A)S)=\omega(A\hat\varphi^2(S)),\quad
A\in \Cal M_\varphi.$$ Hence $\sigma_i(S)=\hat\varphi^2(S)$.
However, since $R$ defines $\hat\varphi$, $\hat\varphi(\bar
R)=d(\varphi)^{-1}R$, and since $\bar R$ also defines $\hat\varphi$,
$\hat\varphi^2(\bar R)=d(\varphi)^{-2}\bar R$ so $\bar R$ transforms
like a character under the modular group. Thus the corresponding
Jones projection, which is a Jones projection for $\hat\varphi$ lies
in $\Cal M^{(0)}_\omega$. The equivalent conditions are therefore
fulfilled whenever $\hat\varphi$ is
selfconjugate.\hfill$\square$\medskip

\proclaim{Theorem 5.25}Let $\rho$ be real or irreducible and
pseudoreal. If 
 $\varphi$ is amenable  and selfconjugate then $d(\rho)=\|m^\rho\|$.
 The same conclusion holds if $\varphi$ is amenable and $\omega$ is
tracial on 
 $\Cal M^{(0)}_\varphi$. \endproclaim \demo{Proof}When $\varphi$ is
amenable, Proposition 5.7 shows that $\Cal M^{(0)}_\omega$ is a
factor and that $\cup_r(\mu^r,\mu^r),\,\, r\in\Bbb N_0$ is dense in
$\Cal M^{(0)}_\omega$. Thus $\mu(\Cal M^{(0)}_\omega)\subset\Cal
M^{(0)}_\omega$ is a strongly amenable inclusion of type
$II_1$--factors by Theorem 4.1.2(iv) of [\ref(Popa2)]. Hence
$\|m^\mu\|=d(\mu)$, see the remark following the proof of Theorem
5.3.2 of [\ref(Popa2)]. But if $\varphi$ is selfconjugate, then  by
Lemma 5.24,  $d(\rho)= d(\mu)=\|m^\mu\|\leq\|m^\rho\|\leq d(\rho)$,
as required. If $\omega$ is tracial on  $\Cal M^{(0)}_\varphi$, one
may either note that the equivalent conditions of Lemma 5.24 are
satisfied, or, for a more elementary proof use the endomorphism
$\tau$ and Proposition  5.22. \hfill$\square$\medskip 

We do not know whether Theorem 5.25 is valid for all amenable
$\varphi$. The results we have obtained justify the term
``amenable'' for $\varphi$ by analogy of Popa's use of strongly
amenable for the graph of an inclusion [\ref(Popa3)].

We now consider  further examples to complement the above
discussion. We begin with the fundamental representation  of
$SU_q(2)$, $\rho_1$,  first taking $q$ to be a root of unity. The
dual of $SU_q (2)$, with $q$ a root of the unity, contains only
finitely many irreducible objects; it is  a rational tensor
$C^*$--category [\ref(FrKer)]. As a consequence of Popa's
classification [\ref(Popa2)] we have\par

\proclaim{Proposition 5.26} If $q^n=1$ then $\hat\rho_1$ is an
irreducible endomorphism of $\Cal M_{\rho_1}.$ The inclusion
$\hat\rho_1(\Cal M^{(0)}_{\rho_1})\subset \Cal M^{(0)}_{\rho_1}$ is
isomorphic to Wenzl's subfactor [\ref(Wen1)].\endproclaim
\demo{Proof} The inclusion $\hat\rho_1(\Cal M^{(0)}_{\rho_1})\subset
\Cal M^{(0)}_{\rho_1}$ and Wenzl's subfactor have the same standard
invariant  (they both generate a tensor category isomorphic to the
dual of $SU_q (2)$), hence they are isomorphic inclusions
[\ref(Popa2)].\hfill$\square$ \medskip

Now suppose that $q\in (0,1)$.  Since the principal graph of
tensoring by $\rho_1$ is $A_\infty$, $\Cal M^{(0)}_{\rho_1},$  the
grade zero part of $\Cal M_{\rho_1},$ is a factor [\ref(Jones)]:
indeed we have 

\proclaim{Proposition 5.27} $\Cal M_{\rho_1}$ is a factor of type
$III_\lambda, \lambda=(q+q^{-1})^2$.\endproclaim \demo {Proof} 
$\Cal M^{(0)}_{\rho_1}$ is a factor, hence $\Cal M_{\rho_1}$ is a
factor by Corollary 5.3. Since in this case
$(\rho_1{}^r,\rho_1{}^{r+1})=0,\,\,r\in \Bbb N_0$  but
$(\iota,\rho_1{}^2)\neq 0$, the statement follows by Corollary 5.3
and the computation of $d(\rho_1)$ in Section 3 or compare e.g.{}
[\ref(Yam)].\hfill$\square$ \medskip  We therefore have an inclusion
of factors $\hat\rho_1(\Cal M_{\rho_1})\subset \Cal M_{\rho_1}$
associated with the fundamental representation of $SU_q (2)$,\par

\proclaim{Proposition 5.28}  $d(\hat\rho_1)=2.$\endproclaim

\demo {Proof}$\hat\rho_1(\Cal M^{(0)}_{\rho_1})\subset\Cal
M^{(0)}_{\rho_1}$  is an inclusion of type $II_1$--factors with
Jones index $d(\rho_1)^2$. The inclusion 
 is known to be reducible, in fact it is just given by the  direct
sum of two automorphisms. 
 Hence if $\tau$ denotes the restriction of $\hat\rho_1$ to $\Cal
M^{(0)}_{\hat\rho_1},\,\, 
 d(\tau)=2.$ However, $d(\tau)=d(\hat\rho_1)$ by Proposition
5.22.\hfill$\square$\medskip

\proclaim{Corollary 5.29} The fundamental representation of
$SU_q(2),\,q\in (-1,1),$ is not amenable.\par  There are uncountably
many inclusions of factors with graph $A_\infty$ and hence
uncountably many Ocneanu connections [\ref(Ocn)] on the graph
$A_\infty.$\endproclaim\par

We note that Corollary 5.29 yields an example of a non--amenable
tensor $C^*$--category whose tensor product is commutative. Now
there is a certain belief that  Abelianness is responsible for the
amenability of the dual of a locally compact group. This example
shows that matters are not quite so simple and that one has to 
distinguish different notions of Abelianness. The notion in the
following result  is sufficient for amenability.

\proclaim{Lemma 5.30} Let $\rho$ satisfy the conditions of Corollary
5.13 and suppose  there is a unitary $\varepsilon\in
(\rho^2,\rho^2)$ such that $$\varepsilon(s,1)\circ X\otimes
1_\rho=1_\rho\otimes X\circ\varepsilon(r,1),\quad X\in
(\rho^r,\rho^s),$$ where $$\varepsilon(0,1)=1_\rho,\quad
\varepsilon(r,1)=\varepsilon\otimes 1_{\rho^{r-1}}\circ
1_\rho\otimes\varepsilon(r-1,1)$$ then on $\Cal O_\rho$ we have 
$$\hat\rho(X)=\lim_{s\to\infty}\varepsilon(s+k,1)X\varepsilon(s,1)^*,\quad
X\in\Cal O^{(k)}_\rho$$ and in particular
$$\varepsilon(s,1)X=\hat\rho(X)\varepsilon(r,1),\quad
X\in(\hat\rho^r,\hat\rho^s).$$ 

The analogous relations hold on $\Cal M_\varphi$ if the limit is
understood in the $s$--topology provided the state $\omega$ defined
by $\varphi$ satisfies  $$\omega(\varepsilon
X\varepsilon^*)=\omega(X) \quad X\in\Cal M^{(0)}_\varphi,$$  thus in
particular if $\omega$ is tracial on $\Cal M^{(0)}_\varphi.$ More
generally, the  condition is fulfilled if and only if $\varepsilon$
commutes with $e^{-H_2}$, where  $e^{-H_2}$ is as in Lemma 5.2. The
$\varphi$ that satisfy this condition are amenable. In particular,
$\rho$ is amenable. \endproclaim \demo{Proof} In the case of $\Cal
O_\rho$ the proof follows Corollary 4.7 and Lemma 4.14 of
[\ref(DR1)]. In the case of $\Cal M_\varphi$, it suffices to show
that $\varepsilon(s+k,1)X\varepsilon(s,1)^*
\xi_\omega\to\hat\rho(X)\xi_\omega$ since $\xi_\omega$ is separating
for $\Cal M_\varphi$. Now if $Y\in (\rho^r,\rho^{r+k})$ and $s\geq
r$,
$$\|(\varepsilon(s+k,1)X\varepsilon(s,1)^*-\hat\rho(X))\xi_\omega\|\leq\|\varepsilon(s+k,1)
(X-Y)\varepsilon
(s,1)^*\xi_\omega\|+\|\hat\rho(X-Y)\xi_\omega\|=2\|(X-Y)\xi_\omega\|,$$
as required.  The rest of the proof follows as for $\Cal O_\rho$
replacing convergence in the norm by convergence in the s--topology.
This is possible since $\hat\varphi$ is continuous in this topology
(cf. Lemma 5.1).  \hfill$\square$\bigskip Of course, we can find an
$\varepsilon\in (\rho^2,\rho^2)$ with the required properties  if
$\Cal T_\rho$ admits a unitary braiding. However, the above
conditions are weaker and  are known to be fulfilled when $\rho$ is
the regular representation of a Kac--system in  the sense of \S 6 of
[\ref(BaaSka)]. We will refer to $\varepsilon(\cdot,1)$ as a unitary
{\it  half--braiding}. Thus a unitary half--braiding suffices for
amenability  but a non-unitary braiding, as in the case of
$SU_q(2),\,\, q\in (-1,1)$ does not.

We end this section by summing up the result of the previous lemma
in the form of a theorem. \proclaim{Theorem 5.31} Let $\rho$ be a
real object or an irreducible pseudoreal object  in a strict tensor
$C^*$--category. Then if $\rho$ admits a unitary half--braiding it
is $C^*$--amenable and  amenable and in fact $\varphi$ is amenable
whenever  the state $\omega$ associated with $\varphi$ is tracial on
$\Cal M^{(0)}_\varphi$.  \endproclaim

\beginsection {6. Inclusions of Factors and Hopf Algebras}

As a further application of our methods, we give  a description of
two classes of inclusions of factors, the standard AFD inclusions
and the standard graded AFD inclusions.  As compared with the last
section, there is a change of objective. There it could be said that
we were studying endomorphisms of von Neumann algebras starting with
objects in a tensor $C^*$--category. Here instead, we will be
studying inclusions of von Neumann algebras and we need to start
from more particular data.  We thus define  an {\it abstract
$Q$--system} to be a pair $(\lambda,S)$ consisting of  an object
$\lambda$ in a strict tensor $C^*$--category $\Cal T$  
 with $d(\lambda) > 1$ and an isometry  $S\in (\lambda,\lambda^2)$
such that\par \medskip \item{$a)$}$$\,  1_\lambda \otimes S\circ S =
S \otimes 1_\lambda\circ S$$ \smallskip

\item{$b)$}$$\,\cases S^* \circ 1_\lambda \otimes T =1_\lambda  \\
T^*\otimes 1_\lambda \circ S=1_\lambda \endcases$$ \smallskip for
some (unique) $T\in (\iota,\lambda)$.\medskip If in addition
\item{$c)$}  $(\iota, \lambda)$ is one--dimensional, the $Q$--system
will be said to be irreducible.\smallskip

Apart from occasional remarks, we shall continue to suppose, as in
\S 5, that the tensor unit $\iota$ is irreducible.  We recall from
Section 2  that if $\iota$ is irreducible we automatically have such
a structure  whenever $\lambda$ is a canonical object.   If $\iota$
is reducible, it is not clear that we can take $S$ to be an isometry
since it is not clear that $R^*\circ R$ is invertible. Invertibility
in fact corresponds to the notion of finite scalar index  in the
theory of inclusions of von Neumann algebras. An isomorphism of
abstract $Q$--systems will be an isomorphism of tensor
$C^*$--categories $\Cal T_{\lambda}$ preserving the distinguished
intertwiners $T,S$.\par Note that the object $\lambda$ is
selfconjugate and even real, indeed the conjugate equation is
satisfied with $\lambda = \bar\lambda$,  $R =\overline R = S\circ
T$. We denote the corresponding normalized left inverse by $\varphi$
so that  $$d(\lambda)\leq d(\varphi)=\|T\|^2\|S\|^2=\|T\|^2.$$  We
shall see however that $d(\lambda)=\|T\|^2$ if the $Q$-system is
irreducible. It should be noted that not all faithful  normalized
left inverses of $\lambda$ can arise in this manner. In particular,
any such  $\varphi$ is  obviously selfconjugate.

When $\Cal T$ is the category of unital  normal endomorphisms of a
von Neumann algebra  $\Cal M$, it was shown in [\ref(Long6)] that
$\lambda$ is the canonical endomorphism  corresponding to a uniquely
determined von Neumann subalgebra $\Cal N$ of $\Cal M$. The main
step  in the proof consisted in writing an explicit formula for a
conditional expectation onto this subalgebra, namely $$A\mapsto
S^*\lambda(A)S.$$ This step has a categorical analogue. We shall not
need a formal definition of a conditional expectation on a
$C^*$--category here. We will only deal with those arising from left
and  right inverses. If $\varphi$ is a left inverse for $\rho$ then
defining for each pair of objects $\sigma,\,\tau$,
$$\delta_{\sigma,\tau}(X):=1_\rho\otimes\varphi_{\sigma,\tau}(X),\quad
X\in(\rho\sigma,\rho\tau)$$ gives us a conditional expectation onto
the image of the functor of tensoring on the left by $1_\rho$.
Similarly, we get a conditional expectation onto the image of the
functor of tensoring on the right by $1_\rho$ by starting with a
right inverse. The idea here in dealing with a canonical object
$\lambda$ is  to imagine that $\lambda=\rho\bar\rho$ and write down
the conditional expectation onto the image of tensoring on the left
by $\rho$. We work on the category $\Cal T_\lambda$ and define for
$r,\,s\in \Bbb N$,  $$\delta_{r,s}(X):= S^*\otimes
1_{\lambda^{s-1}}\circ 1_\lambda\otimes X\circ S\otimes
1_{\lambda^{r-1}},\quad X\in (\lambda^r,\lambda^s).$$  It is easily
checked that $\delta_{r,r}$ is faithful and a completely positive
projection of norm one  and the same holds for a matrix with entries
$\delta_{r,s},\,\,r,s=1,2,\dots,n.$ Since $S$ is in the image of
$\delta$ by $a)$,  $S\circ S^*$ is also in that image. In fact
$\delta_{2,2}(S\circ S^*)\geq \delta_{1,2}(S)\circ \delta_{1,2}(S)^*
= S\circ S^*$ and letting $\delta_{2,2}$ act on the difference 
$\delta_{2,2}(S\circ S^*) - S\circ S^*$  we get $0$, thus
$\delta_{2,2}(S\circ S^*) = S\circ S^*$ because $\delta$ is
faithful. Computing  $X^*\circ X$ for $X:= 1_\lambda\otimes S^*\circ
S\otimes 1_\lambda-S\circ S^*$, we see that  for any abstract
$Q$--system we have \item{$d)$} $$1_\lambda\otimes S^*\circ S\otimes
1_\lambda=S\circ S^*.$$  This condition was included as part of the
definition of a $Q$--system in [\ref(Long6)]. \proclaim{Lemma 6.1}
The image of $\delta_{r,s}$ is the set of $Y\in
(\lambda^r,\lambda^s)$ such that $$1_\lambda\otimes Y\circ
S\otimes1_{\lambda^{r-1}} = S\otimes 1_{\lambda^{s-1}}\circ Y.$$
\endproclaim \demo{Proof} It is a simple computation using $a)$ and
$d)$ to show that $Y=\delta_{r,s}(X)$ satisfies the required
relation. On the other hand if $Y$ satisfies the relation, 
$Y=\delta_{r,s}(Y)$.\hfill$\square$\medskip 

On the basis of this lemma, it is clear that the image of $\delta$
is a $C^*$--subcategory which is closed under tensor products. It
fails to be a tensor $C^*$--category because the tensor unit is not
in the subcategory. There is also another conditional expectation,
obtained conceptually by thinking of a right inverse for $\bar\rho$
and leading to another $C^*$--subcategory closed under tensor
products. The formulae are obtained by interchanging tensoring on
the left by $\lambda$ with tensoring on the right by $\lambda$. This
second possibility will, however, play no role here.

Now an abstract $Q$--system gives rise to an inclusion $\Cal
N_\varphi\subset \Cal  M_\varphi$ of von Neumann algebras, where
$\Cal N_\varphi$ denotes the von Neumann subalgebra generated by the
above subcategory. Furthermore, we have an  endomorphism
$\hat\lambda$ of $\Cal M_\varphi$ and intertwiners $T\in
(\iota,\hat\lambda),\,S\in(\hat\lambda,\hat\lambda^2)$ satisfying
the relations $a)$ and $b)$ above. Thus we have a $Q$--system in the
sense of [\ref(Long6)] except that $\Cal M_{\varphi}$ might fail to
be a factor. It follows by the analysis in [\ref(Long6)] that
$\hat\lambda$ is the canonical endomorphism of $\Cal M_\varphi$ with
respect to the subalgebra $\Cal N_\varphi$ since this is just the
image of $\Cal M_\varphi$ under  the conditional expectation $\Cal
E,\,\,\Cal E(X):= S^*\hat\lambda(X)S$. Thus  $$\Cal N_\varphi
=\{X\in \Cal M_\varphi : SX=\hat\lambda(X)S \}.$$  In particular,
$\Cal N_\varphi$ is stable under gauge transformations. The normal
conditional expectation will be understood as part of the data when
talking about inclusions and we write $\Cal N_\varphi {\overset{\Cal
E}\to\subset}\Cal M_\varphi$ to emphasise this. The $Q$--system not
only singles out the conditional expectation but  also a particular
choice of the Jones tunnel and a set of Jones projections. The Jones
tunnel is $$\quad\cdots\hat\lambda(\Cal
N_\varphi)\subset\hat\lambda(\Cal M_\varphi)\subset\Cal
N_\varphi\subset\Cal M_\varphi.$$ The Jones projections are as
follows: $E_0=d(\varphi)^{-1}TT^*,\,\, E_1=SS^*,$ and $E_{i+2}=
\hat\lambda(E_i),\,\,i\in \Bbb N_0.$ The Jones index of the
inclusion is $d(\varphi)$.

\proclaim{Lemma 6.2} If $X\in (\hat\lambda^r,\hat\lambda^s)$ then 
$$d^{r+s}E_0E_1\cdots E_{2s}XE_{2r}\cdots E_1E_0=\hat\lambda(X)E_0,$$
where $d=d(\varphi)$. \endproclaim
\demo{Proof}$$TT^*SS^*\cdots\hat
\lambda^{s-1}(SS^*)\hat\lambda^s(TT^*)X\hat\lambda^r(TT^*)\cdots
TT^*=T\hat\lambda^s(T^*)X\hat\lambda^r(T)T^*=d^2\hat\lambda(X)E_0,$$
as required.\hfill$\square$\medskip
 
 Passing to the grade zero part we get an inclusion $\Cal
N^{(0)}_\varphi\,\,{\overset{\Cal E_0}\to\subset}\Cal
M^{(0)}_\varphi$   of von Neumann algebras and we denote by $\tau$
the restriction of $\hat\lambda$ to the homogeneous part.  Since the
Jones projections are of grade zero, $$\quad\cdots\tau(\Cal
N^{(0)}_\varphi)\subset\tau(\Cal M^{(0)}_\varphi)\subset\Cal
N^{(0)}_\varphi\subset\Cal M^{(0)}_\varphi$$ is a Jones tunnel.
Furthermore, by Corollary 5.6 
$$(\tau^r,\tau^r)=(\hat\lambda^r,\hat\lambda^r),\quad r\in\Bbb N_0.$$
Now the union of the relative commutants in the tunnel is just
$\cup_r(\tau^r,\tau^r)$  and is hence dense in $\Cal
M^{(0)}_\varphi$. In fact our inclusion is just  derived by the
canonical procedure from the relative commutants in the tunnel. To
see this we need only verify that the unique normal
$\tau$--invariant state $\omega$ satisfies  $$\omega(\Cal
E_i(X))=\omega(X),\quad X\in \Cal N_{i-1},\,\,i\in \Bbb N_0,$$ where
$\Cal E_i$ is the conditional expectation from $\Cal N_{i-1}$ onto
$\Cal N_i$ defined by  $$\Cal E_i(X)E_{i-1}=E_{i-1}XE_{i-1}.$$ 
Since $\Cal E_{i+2}\tau(X)=\tau\Cal E_i(X)$ and $\omega$ is
$\tau$--invariant, it is enough to test $i=1\,,2$ and hence
alternatively to show that $$\omega(\Cal E_1\Cal
E_0(X))=\omega(X),\quad X\in \Cal M^{(0)}_\varphi.$$ However, we
have  $$E_0\Cal
E_0(X)E_0=d^{-2}TT^*S^*\tau(X)STT^*=\tau\hat\varphi(X)E_0,$$  so
that $\Cal E_1\Cal E_0=\tau\hat\varphi$ and the result is now
obvious.  As far as the zero grade part goes, our construction
therefore yields inclusions of a very special type.  In fact, when
$\omega$ is just a trace on $M^{(0)}_\varphi$ so that we have a
$II_1$--factor, we have just verified that our inclusion is strongly
amenable in the sense of [\ref(Popa2)]. 

 Our next goal is  to show that our inclusion $\Cal
N_\varphi{\overset{\Cal E}\to\subset}\Cal M_\varphi$ can be
reconstructed from its zero grade part.
 In fact, it is just the cross product of its grade zero part by an
endomorphism.  To, see this we note that $\tau$ is a canonical
endomorphism in the sense of Choda[\ref(Cho)]  since we have a
faithful normal $\tau$--invariant state, $E_0\Cal M^{(0)}_\varphi
E_0= \tau(\Cal M^{(0)}_\varphi) E_0$ by Lemma 6.2 and the Jones
projections fulfill the Jones  relations. We can therefore define
the cross product $\Cal M$ of $\Cal M^{(0)}_\varphi$  by the
endomorphism $\tau_0,\,\,\tau_0(X):=\tau(X)E_0.$ $\Cal M$ is
generated  as a von Neumann algebra by $\Cal M^{(0)}_\varphi$ and an
isometry $W$ satisfying  $$WXW^*=\tau(X)W,\quad X\in \Cal
M^{(0)}_\varphi.$$ $\tau$ extends uniquely to an endomorphism
$\tilde\tau$ of $\Cal M$ such that 
$$\tilde\tau(W)=d(\varphi)E_2E_1W.$$ Furthermore, it is easy to see
that $\Cal M$ carries a continuous action $\alpha$ of the  circle
with $\alpha_t(W)=tW,\,\,t\in \Bbb T$ and $\Cal M^{(0)}_\varphi$ as
fixed-point algebra. In fact, the restriction of $\omega$ to $\Cal
M^{(0)}_\varphi$ extends to an $\alpha$--invariant and hence
faithful state $\tilde\omega$ on $\Cal M$.

We now compare $\Cal M,\,\tilde\tau$ and $\Cal
M_\varphi,\,\hat\lambda.$ The ${}^*$-- subalgebra generated by $W$
and $\Cal M^{(0)}_\varphi$ is uniquely determined by $\tau_0,$
(compare the discussion following Lemma 3.7 of [\ref(DR2)]). Hence
there is a unique injective morphism $\chi$ of this
${}^*$--subalgebra into $\Cal M_\varphi$  which is the identity on
$\Cal M^{(0)}_\varphi$ and satisfies $\chi(W)=\sqrt{d(\varphi)}T.$
It intertwines the actions of the circle and $\omega\circ\chi$
coincides with the restriction of $\tilde\omega$. Thus $\chi$
extends to an isomorphism of $\Cal M$ and $\Cal M_\varphi$ and
since  $$d(\varphi)E_2E_1T=\hat\lambda(TT^*)SS^*T=\hat\lambda(T),$$
$\chi$ carries $\tilde\tau$ into $\hat\lambda$. Thus we recover the
initial inclusion $\Cal N_\varphi\subset\Cal M_\varphi$ and have
proved the following result. \proclaim{Theorem 6.3} Let
$(S,\lambda)$ be an abstract $Q$--system and $\varphi$ the
associated left inverse of $\lambda$, then $\Cal M_\varphi$ is just
the cross-product of  $\Cal M^{(0)}_\varphi$ by the endomorphism
$X\mapsto\hat\lambda(X)E_0$ of $\Cal M^{(0)}_\varphi$ and
$\hat\lambda$ is the canonical extension of its restriction to 
$\Cal M^{(0)}_\varphi$.\endproclaim

We now want to characterize the two classes of inclusion discussed
in this section and begin with the appropriate definitions. Suppose
we have an inclusion of properly infinite factors $\Cal
N\,{\overset\Cal E\to\subset}\Cal M$ so that there is a $Q$--system 
$(\lambda,S)$, where $\lambda$ is an endomorphism of $\Cal M$ and
$S\in \Cal M$  satisfies $\Cal E(X)=S^*\lambda(X)S,\quad X\in\Cal
M$.  
 Since any other such $Q$--system is of the form $\lambda',S'$,
where $\lambda'={\text ad}(U)\lambda$ and  $S'=US$, $U$ being a
unitary of $\Cal N$, $U$ induces an isomorphism of abstract
$Q$--systems. Hence different associated $Q$--systems yield
isomorphic inclusions $\Cal N_\varphi\,{\overset\Cal
E\to\subset}\Cal M_\varphi.$ If this inclusion is isomorphic to the
original inclusion, we say that $\Cal N\,{\overset\Cal
E\to\subset}\Cal M$  is a {\it standard graded AFD inclusion}. On
the other hand, even if our factors are not properly infinite, we
can choose a Jones tunnel and construct the associated inclusion
obtained from  the relative commutants. Again, as is well known,
different choices of the tunnel lead to  isomorphic inclusions. If
this inclusion is isomorphic to the original inclusion, we say  that
$\Cal N\,{\overset\Cal E\to\subset}\Cal M$  is a {\it standard AFD
inclusion}. 

\proclaim{Theorem 6.4} Every standard AFD inclusion is the grade
zero part of a standard graded AFD  inclusion. \endproclaim
\demo{Proof} If we tensorialize a Jones tunnel with an infinite
factor, the inclusion derived from the relative commutants remains
unchanged. The tensorialized inclusion, however, has an associated 
$Q$--system. This $Q$--system defines a model graded inclusion whose
zero grade part is, as  we have seen above, derived from the
relative commutants in the associated Jones tunnel.
\hfill$\square$\medskip

Of course, the standard graded AFD inclusion defined by an abstract
$Q$--system $(\lambda,S)$ depends only on the $Q$--system
$(\hat\lambda,S)$. $(\lambda,S)$ and $(\hat\lambda,S)$ are
isomorphic if and only if the associated left inverse $\varphi$ of
$\lambda$ is amenable.  In this case, we say that the $Q$--system
$(\lambda,S)$ is {\it amenable}. Hence, we can sum up the above
discussion in the following way. \proclaim{Corollary 6.5} There is a
natural 1-1 correspondence between isomorphism classes of amenable
abstract $Q$--systems, standard graded AFD inclusions  and standard
AFD inclusions.\endproclaim

Of course, it is hardly surprising that the image of the irreducible
amenable $Q$--systems under this correspondence is precisely the set
of irreducible inclusions.  To see this we note that in a
$Q$--system $(\lambda,S)$, there is a 1--1 correspondence between
the elements of $(\iota,\lambda)$ and those elements $Y$ of 
$(\lambda,\lambda)$ that satisfy $Y\otimes 1_\lambda\circ S=S\circ
Y$. $X\in (\iota,\lambda)$ corresponds to $S^*\circ 1_\lambda\otimes
X\in (\lambda,\lambda)$. However, the relative commutant of $\Cal
N_\varphi$ in $\Cal M_\varphi$ is easily seen  to be $\{ Y\in
(\hat\lambda,\hat\lambda)\,\,: YS=SY\}.$ Now, the relative 
commutant of $\Cal N^{(0)}_\varphi$ in $\Cal M^{(0)}_\varphi$ is
just $\{Y\in (\tau,\tau)\,\,: YE_1=E_1Y\,\}$, where $E_1=SS^*$. It
is not so obvious that  this relative commutant is just the grade
zero part of the previous relative commutant but  the question is
closely related to that resolved in Corollary 5.6 and the same
techniques work.  \proclaim{Lemma 6.6} $$\Cal
N^{(0)}_\varphi{}'\cap\Cal M^{(0)}_\varphi=\Cal N_\varphi{}'\cap
\Cal M_\varphi.$$ \endproclaim \demo{Proof}Let $E,\,E'$ be two
projections in $\Cal M^{(0)}_\varphi$. If there is an $X\in\Cal
N_\varphi$ with $EXE'\neq 0$, we see as in the proof of Lemma 5.5
that there is a $Y\in\Cal N^{(0)}_\varphi$ such that $EYE'\neq 0.$
The result now follows by the argument of Corollary
5.6.\hfill$\square$\medskip

Now the Jones index of an irreducible inclusion coincides with the
minimal index. Hence for an irreducible $Q$--system
$d(\lambda)=\|T\|$.

Let us call a $Q$--system {\it rational} if the completion of $\Cal
T_\lambda$ under  subobjects contains only finitely many
inequivalent irreducible objects. This is obviously  corresponds to
the notion of finite depth in the theory of subfactors. As one sees
from the  class of examples discussed following Corollary 5.6, in
our framework, the Jones index does not necessarily coincide with
the minimal index for an inclusion of factors of finite depth,  i.e.
such inclusions are not necessarily extremal in Popa's
terminology[\ref(Popa2)]. This is because we have not imposed the
condition that $\omega$ defines a trace on the grade zero part. The
Jones index, therefore, has no reason to be a Perron--Frobenius 
eigenvalue. For the same reason, we do not know whether a rational
$Q$--system is  necessarily amenable since $\Cal M_\varphi$ might
possibly fail to be a factor. Of  course, if $\omega$ is a trace on
$\Cal M^{(0)}_\varphi$ then the $Q$--system is amenable since finite
depth inclusions are amenable in Popa's setting[\ref(Popa2)].

As an application, we characterize finite--dimensional Hopf algebras
in terms of tensor $C^*$--categories.

Notice that, as a consequence of Popa's classification
[\ref(Popa2)], we also have a bijective correspondence between
amenable (finite depth) irreducible inclusions of $II_1$--factors
and amenable (finite depth) irreducible inclusions of injective
$III_1$--factors.

\proclaim{Theorem 6.7} There is a bijective correspondence between 
\item {$(i)$}  Irreducible abstract Q--systems $(\lambda, S)$ such
that $$\lambda \otimes \lambda \simeq \lambda \oplus \lambda \cdots
\oplus \lambda \qquad (d\,\text{times}\,\lambda)$$ \item {$(ii)$}
Finite--dimensional complex Hopf algebras of dimension
$d$.\endproclaim \par

\demo{Proof} This follows by Corollary 6.2 and the results in
[\ref(Long6)]. \hfill$\square$\medskip

We recall from $\S$ 6 of [\ref(Long6)] that to each $Q$--system
defined concretely in terms of a von Neumann algebra there is a dual
$Q$--system. We intend in a sequel to study this duality directly in
terms of abstract $Q$--systems.

\beginsection{7.  2--$C^*$--Categories}

The formalism we have presented so far based on the notion of tensor 
$C^*$--categories is not sufficiently general to cover all
interesting  applications and we must generalize to
2--$C^*$--categories. In fact,  this does not present the slightest
difficulty and we could have  worked with the more general notion
from the outset but for our wish  to avoid  unfamiliar  and hence
superficially technical concepts.

Let us begin by understanding the typical situation not yet covered
by  our formalism. Suppose we wish to consider the minimal index of
an inclusion of  factors $\Cal N\subset\Cal M$ directly. To fit in
with our formalism to date, we need an  object in a tensor
$C^*$--category. If $\Cal M$ and $\Cal N$ are isomorphic, we can 
compose the inclusion with an isomorphism to obtain an endomorphism
which  is an object in a tensor $C^*$--category. The intrinsic
dimension of the  endomorphism is independent of the choice of
isomorphism and serves to  define the index of the inclusion.
However, there is a simpler way of  looking at things. We should
look at morphisms rather than endomorphisms of $C^*$--algebras.
There is no difficulty in generalizing the concept of conjugate. A
morphism $\rho$ from $\Cal A$ to $\Cal B$ has a morphism $\bar\rho$
from $\Cal B$ to $\Cal A$ as a conjugate if there are
$R\in(\iota_{\Cal A}, \bar\rho\rho)$ and $\bar R\in(\iota_{\Cal
B},\rho\bar\rho)$ such that $$\bar R^*\rho(R)=1_{\Cal B};\quad
R^*\bar\rho(\bar R)=1_{\Cal A}.$$ As before, $\iota$ denotes the
identity automorphism of the algebra in  question. All we need
therefore is the correct way of formalizing this  situation and this
is provided for by the notion of a 2--$C^*$--category.

A 2-category may be conveniently defined as a set with two
composition laws, $\otimes$ and $\circ$, each making it into a
category. In addition \item{$a$)} Every $\otimes$--unit is an
$\circ$--unit,\par \item{$b$)} The $\otimes$--composition of
$\circ$--units is a $\circ$--unit.\par \item{$c$)} $S'\otimes
T'\circ S\otimes T=(S'\circ S)\otimes (T'\circ T)$ whenever the
right hand side is defined.\par

A 2-category is a 2--$C^*$--category if, given any pair of
$\otimes$--units,
 the category whose elements are the subset having that pair as left
and right $\otimes$--units respectively and $\circ$ as composition
is a $C^*$--category and if furthermore the $\otimes$--composition
is bilinear. 

To illustrate this with the motivating example, we get a
2--$C^*$--category by  taking the set of intertwining operators of
morphisms between elements of a given  set of unital
$C^*$--algebras. $\circ$--composition is just the ordinary
composition of intertwining operators. If $S$ is an intertwining
operator between $\sigma$ and  $\sigma'$ and $T$ an intertwiner
between $\tau$ and $\tau'$ then $S\otimes T$ is defined if and only
if the composition $\sigma\tau$ (and hence also  $\sigma'\tau'$) of
morphisms is defined and $S\otimes T := S\sigma(T)$ considered as an
intertwiner from $\sigma\tau$ to $\sigma'\tau'$. The $\circ$--units
are  therefore the unit intertwiners for the morphisms and may be
identified with the  morphisms themselves and the $\otimes$--units
are the unit intertwiners for the identity automorphisms and may be
identified with the algebras themselves.

But this is only one example and not necessarily the most convenient
one, even for studying inclusions of factors, so we give further
examples and stress the  point that the subject matter of this
paper, conjugation, dimension and trace,  is just an elementary part
of abstract theory of 2--$C^*$--categories.

Of course, the classical example of a 2--category is the 2--category
of natural  transformations between functors between categories.
Hence there is an obvious generalization of our motivating example
which involves replacing unital  $C^*$--algebras by
$C^*$--categories and taking the set of bounded natural 
transformations between $^*$--functors between $C^*$--categories. We
refrain from  giving details since we shall not pursue this example
here. The basic definitions  of the terms involved can be found in
[\ref(GLR)].

Of more relevance to the present paper is the fact that the
motivating example  can be generalized by considering the set of
intertwiners  of morphisms from elements of a given set of
$C^*$--algebras with unit into matrix algebras over the given set.
$T$ is an intertwiner between two such morphisms $\rho$ and
$\sigma$, written $T \in (\rho,\sigma)$, if $\rho : \Cal A \mapsto
M_m(\Cal B)$ and  $\sigma :\Cal A \mapsto M_n(\Cal B)$ then $T$ is
an $n\times m$ matrix with entries
 from $\Cal B$ such that $T=T\rho(I)=\sigma(I)T$ and
$$T\rho(A)=\sigma(A)T,\quad A\in \Cal A.$$ If $\tau : \Cal B \mapsto
M_r(\Cal C)$ then $\tau\rho : \Cal A\mapsto M_{rn}(\Cal  C)$ where
$\tau\rho(A)$ is defined by letting $\tau$ act on each component of 
$\rho(A)$. To give an explicit formula, if $\rho_{ij}(A)$ denotes
the $ij$-th entry of $\rho(A)$, then
$$(\tau\rho)_{ij}=\tau_{i_1j_1}\rho_{i_2j_2}\quad \text{where}\quad
i=i_1+ri_2, \quad j=j_1+rj_2.$$ If $R\in (\rho,\rho')$ and $T\in
(\tau,\tau')$ then $$T\otimes R=T\tau(R)=\tau'(R)T$$ is an element
of $(\tau\rho,\tau'\rho')$ This 2--$C^*$--category has the 
advantage of being closed under finite direct sums and subobjects.

We may also modify the above examples to cater for $C^*$--algebras
without unit. To do this, given two morphisms $\rho$ and $\sigma$
from $\Cal A$ to matrix  algebras over $\Cal B$, an intertwiner $T$
must be allowed to be in the appropriate
 matrix algebra over the multiplier algebra $M(\Cal B)$ and must
satisfy the  condition that $T=\text{lim}T\rho(X_\alpha)$ where
$X_\alpha$ is an approximate  unit in $\Cal A$. 

In much the same spirit one could consider Hermitian bimodules over
a set of  $C^*$--algebras, cf. [\ref(Rob1)], \S8. Thus if $\Cal A$
and $\Cal B$ are $C^*$--algebras  we need a right Hermitian $\Cal
B$--module, i.e. a right Banach $\Cal B$--module  $X$ with a $\Cal
B$--valued inner  product ${<}\,,\,{>}$, conjugate--linear in the
first variable and linear in the  second, such that, for all $x,y
\in X$ and $B\in\Cal B$, \smallskip      \item {a)} ${<}x,x{>}\geq
0,$ \quad \item {b)} ${<}x,yB{>}={<}x,y{>}B,$ \item {c)}
${<}x,y{>}^*={<}y,x{>},$ \quad \item {d)} $\| x\|^2=\|{<}x,x{>}\|.$
\medskip $\Cal A$ then acts on the left on $X$; in other words we
have a morphism of $\Cal A$ into  the $C^*$--algebra of $\Cal
B$--module maps of $X$ into itself with adjoint. If $Y$ is  then an
Hermitian $\Cal B$--$\Cal C$--bimodule, we can define the tensor
product  $X\otimes_{\Cal B}Y$ as an Hermitian $\Cal A$--$\Cal
C$--bimodule spanned by  elements $x\otimes y,\, x\in X,\, y\in Y$
with a $\Cal C$--valued inner product such that  $${<}x\otimes
y,x'\otimes y'{>}={<}y,{<}x,x'{>}y'{>}.$$ Since $X\otimes_{\Cal B}Y$
is not defined uniquely by this procedure but only  up to
isomorphism, we will not get a 2--$C^*$--category whose elements
are  bimodule homomorphisms since the $\otimes$--composition will
not be strictly  associative. Rather than generalizing, and hence
complicating, our formalism  to deal with a problem which already
arises at the level of tensor products  of Hilbert spaces, we
content ourselves with a framework where the  $\otimes$--composition
is strictly associative.

Consider a $C^*$--category and pick for each object $A$ a
$C^*$--subalgebra  $\Cal A$ of $(A,A)$. An Hermitian bimodule
associated with this data is a  closed subspace $X$ of some $(A,B)$
such that\smallskip \item {a)} $X\Cal A\subset X$, \item {b)} $\Cal
B X\subset X$ \item  {c)} $x^*x'\in \Cal A, \quad x,x'\in X.$
\medskip  If $X\subset (A,B)$ and $Y\subset (B,C)$ are two such
Hermitian bimodules their  tensor product is the closed subset of
$(A,C)$ generated by the elements of the  form $yx, x\in X, y\in Y$.
The resulting Hermitian bimodule will be denoted  by $YX$. If now
$X'\subset (A,B)$ and $Y'\subset (B,C)$ are two further Hermitian 
bimodules and if $S$ is a bimodule homomorphism from $X$ to $X'$
having an  adjoint $S^*$ in the sense that $${x'}^*Sx=(S^*x')^*x,\quad
x\in X,\quad x'\in X'.$$ and $T$ is a bimodule homomorphism from $Y$
to $Y'$ then $T\otimes S$ is a bimodule homomorphism from $YX$ to
$Y'X'$. Furthermore, if $\circ$ denotes the composition of  bimodule
homomorphisms then we have $$(T'\otimes S')\circ(T\otimes
S)=(T'\circ T)\otimes(S'\circ S),$$  whenever the right hand side is
defined. The set of such bimodule homomorphisms  forms a
2--$C^*$--category.

Note that if we take a $C^*$--category with a single object, i.e. a
unital  $C^*$--algebra, $\Cal A$, and take $\Bbb C I$ as our
distinguished $C^*$--subalgebra  then the resulting
2--$C^*$--category is just the tensor $C^*$--category of Hilbert 
spaces in $\Cal A$.

Another interesting special case is to take the distinguished
$C^*$--subalgebra  of $(A,A)$ to be $(A,A)$ itself, for each object
$A$. The resulting 2--$C^*$--category  then depends intrinsically on
the $C^*$-category we started from and each of  the associated
bimodules $X\in (A,B)$ is bi-Hermitian in the sense that $X$ has a 
second inner product $x'x^*$ taking values in $(B,B)$ and 
$$x''({x'}^*x)=(x''{x'}^*)x,\quad  x,x',x''\in X,$$ where the identity
is regarded as relating the two inner products and the  module
actions.

After these illustrative examples we turn to the concept of
conjugation for finite  dimensional $\circ$--units. If $\Cal T$ is a
2--$C^*$--category then we can define a full 2--$C^*$--subcategory
$\Cal T_f$ by saying that a $\circ$--unit $\rho$ is in  $\Cal T_f$
if there exists a $\circ$--unit $\bar\rho$ and $R\in  (\iota_{\Cal
A},\bar\rho\rho)$ and $\bar R\in (\iota_{\Cal B},\rho\bar\rho)$ such
that 

$$\bar R^*\otimes 1_\rho \circ 1_\rho\otimes R=1_\rho ; \quad
R^*\otimes 1_{\bar \rho} \,\circ \, 1_{\bar\rho}\,\otimes \bar
R=1_{\bar \rho}.$$ Here $\iota_{\Cal A}$ and $\iota_{\Cal B}$ are
the right and left $\otimes$--units of $\rho$.  These equations will
again be referred to as the conjugate equations. With this
definition  there are obvious analogues of the results of Section 2
whose exact formulation can  be left to the reader. There are also
obvious analogues of the results of Section 3  provided $\iota_{\Cal
A}$ and $\iota_{\Cal B}$ are irreducible but the more challenging
task is to provide analogues without making this assumption.

Here we shall content ourselves with showing that the formalism of
2--C$^*$--categories  that has been developed here allows one to
interpret the Jones index of an inclusion of  factors as the square
of the dimension of a suitable $\circ$--unit in a
2--C$^*$--category. The most obvious approach, namely to consider
the inclusion as a $\circ$--unit in our  motivating example of
2--C$^*$--category, viz. the 2--C$^*$--category of intertwining
operators of morphisms between elements of  a given set of unital
C$^*$--algebras fails in general because the conjugate fails to
exist in this  2--C$^*$--category. Instead, we need the
generalization, discussed above, where we allow  intertwiners
between morphisms into matrix algebras over the given set of
C$^*$--algebras. 

Let $\Cal N\subset \Cal M$ be an inclusion of factors with finite
index and let $\Cal M_1=\langle \Cal M,e_{\Cal N} \rangle$ be the
basic construction. Then the existence of a Pimsner--Popa basis
provides  one with an isomorphism of $\Cal M_1$ with a reduced
subalgebra of some $M_n(\Cal N)$ and hence the inclusion of  $\Cal
M$ 
 into  $\Cal M_1$ provides a morphism of $\Cal M$ into $M_n(\Cal
N)$. This is a $\circ$--unit in our 
 category and is a conjugate for the $\circ$--unit coming from the
original inclusion. More 
 specifically, if $\xi_i\in\Cal M,\,\,i=1,2,\ldots n$ denotes the
Pimsner--Popa basis[\ref(PiPo)] and 
 $\rho$  denotes the inclusion of $\Cal N$ into $\Cal M$, then
$\bar\rho$, a morphism of $\Cal M$ 
 into $M_n(\Cal N)$ is given by: $$\bar\rho_{ij}(m)=E_{\Cal
N}(\xi_i{}^*m\xi_j),$$ where $E_{\Cal N}$ is the conditional
expectation of $\Cal M$ onto $\Cal N$ appearing in the definition 
of the Pimsner--Popa basis. We get a solution $R,\bar R$ of the
conjugate equations for  $\rho$ by setting $$\bar
R_i=\xi_i{}^*;\quad R_i=E(\xi_i{}^*).$$ Computing we see that
$R^*\circ R=1_{\Cal N}$ and $\bar R^*\circ \bar
R=\sum_i\xi_i\xi_i{}^*$. Furthermore the image of the Jones
projection $e_{\Cal N}$ in $M_n(\Cal N)$ is just $R\circ R^*$. Hence
bearing in mind that the inequality of Lemma 2.7 is best possible
and the  Pimsner--Popa characterization of the Jones
index[\ref(PiPo)], we see how the index defined by  Jones for an
inclusion of type $II_1$--factors fits into the framework considered
here. When the inclusion $\rho$ is irreducible, in particular when
the index is less than 4, the  Jones index is just $d(\rho)^2$.
\beginsection {Appendix}

In this appendix we give a canonical procedure for completing a
2--$C^*$--category $\Cal T$  to give a 2--$C^*$--category $C(\Cal
T)$  having (finite) direct sums and subobjects. Before describing
the construction, we  remark that it may, at times, prove more
convenient to think of completion in other terms. Thus in the von
Neumann algebra context one might prefer to tensor with an infinite
factor $\Cal N$, passing for example from End$\Cal M$ to End$(\Cal
M\otimes\Cal N)$.

We first discuss how to add subobjects. The resulting 
2--$C^*$--category has  as elements triples (F,T,E) of elements  of
$\Cal T$, where $E$ and $F$ are projections and $T\circ E=T=F\circ
T$. The composition  laws are defined as follows: $(F,T,E)\circ
(F',T',E')$ is defined if and only if $F'=E$ when the  composition
is $(F,T\circ T',E')$. $(F,T,E)\otimes (F',T',E')$ is defined if and
only if  $E\otimes E'$, $F\otimes F'$ and hence $T\otimes T'$ are
defined and the composition  is given by $(F,T,E)\otimes
(F',T',E')=(F\otimes F',T\otimes T',E\otimes E').$ The adjoint is
 of course given by $(F,T,E)^*=(E,T^*,F),$ and $\| (F,T,E)\|=\| T\|.$
 
 Obviously, $(E,E,E)$ is a right and $(F,F,F)$ a left $\circ$--unit
for $(F,T,E).$ If $A$ and $B$ 
 are right and left 
 $\otimes$--units for $E$, respectively, then $(A,A,A)$ and
$(B,B,B)$ are right and left 
 $\otimes$--units for $(F,T,E).$ Hence, the set of $\otimes$--units
remains unchanged and the 
 $\otimes$--composition of $\circ$--units is a $\circ$--unit. To
show that we have a 2--$C^*$--category, it only remains to check the
interchange law and that $\otimes$ commutes with taking adjoints and
this will be left to the  reader as an exercise.

To add direct sums, we take as elements of the completed
2--$C^*$--category finite matrices  with entries from the original
2--$C^*$--category, where the $ij$-th entry $X_{ij}$ has a right
 and left $\circ$--unit depending just on $j$ and $i$, respectively.
$\circ$-composition 
 is now given by matrix multiplication, $$(X\circ Y)_{ij}=\sum_j
X_{ij}\circ Y_{jk}$$ and 
 is defined provided (each) $X_{ij}\circ Y_{jk}$ is defined. Note
that the entries of $X$ all 
 have a common left and a common right $\otimes$-unit. Hence it is
obvious when 
 $X\otimes X'$ is to be defined and we set $$(X\otimes
X')_{ii',jj'}=X_{ij}\otimes X_{i'j'},$$ 
 where for formal reasons, we should use lexicographic ordering. The
$\otimes$--units remain the same but considered as $1\times
1$--matrices. The $\circ$--units are matrices whose off-diagonal
elements are zero and whose diagonal elements are $\circ$--units.
Hence the $\otimes$--composition of $\circ$--units is again a
$\circ$--unit. As one would expect for matrices, the adjoint is
defined by $(X^*)_{ij}=(X_{ji})^*$ and the norm is defined in the
same way as for a matrix of bounded operators. It will again be left
to the reader to check the interchange law and that $\otimes$
commutes with taking adjoints.

It is easy to check that every $\circ$--unit in the completed
category is, in an  obvious sense, a finite direct sum of
$\circ$--units from the original 2--$C^*$--category. Furthermore,
one can again check easily, that if a 2--$C^*$--category has (finite)
direct sums, then it continues to have direct sums after adding
subobjects as above. Thus we have a canonical procedure for
completing a 2--C$^*$--category with respect to finite direct sums
and subobjects if we first add direct sums  and then add subobjects,
as above. We denote the resulting 2--C$^*$--category by  $C(\Cal
T)$. It is obvious that a 2--$^*$--functor $F$ between
2--C$^*$--categories extends  in a canonical manner to a
2--$^*$--functor $C(F)$ between their respective completions.  Of
course, in the special case that $F$ is a tensor $^*$--functor,
$C(F)$ will also be a  tensor $^*$--functor.

 \medskip 
 \noindent {\it Acknowledgements.} We would like to thank Sorin Popa
for helping us to understand  his work. \bigskip \centerline{\bf
References} \bigskip \references

\end